\documentclass[11pt]{amsart}
  \usepackage[margin=3.3cm]{geometry}
\usepackage{amsfonts,amsmath,amssymb,color,esint}
\usepackage{enumerate}
\usepackage{graphicx}
\usepackage{tikz}
\usetikzlibrary{decorations.markings}
\usetikzlibrary{plotmarks}
\usetikzlibrary{patterns}
\numberwithin{equation}{section}

\usepackage{amsthm}

\usepackage[babel=true,kerning=true]{microtype}

\usepackage{amsthm,leqno}
\usepackage{hyperref}

\newtheorem{theo}{Theorem}[section]

\newtheorem{rem}{Remark}[section]
\newtheorem{defi}{Definition}[section]
\newtheorem{cl}{Claim}[section]
\newtheorem{prop}{Proposition}[section]
\newtheorem{cor}{Corollary}[section]

\newcommand{\eps}{\varepsilon}

\newcommand{\R}{\mathbb{R}}

\begin{document}

\title[]{A variational method for functionals depending on eigenvalues}
\author{Romain Petrides} 
\address{Romain Petrides, Universit\'e Paris Cit\'e, Institut de Math\'ematiques de Jussieu - Paris Rive Gauche, b\^atiment Sophie Germain, 75205 PARIS Cedex 13, France}
\email{romain.petrides@imj-prg.fr}

\begin{abstract} 
We perform a systematic variational method for functionals depending on eigenvalues of Riemannian manifolds. It is based on a new concept of Palais Smale sequences that can be constructed thanks to a generalization of classical min-max methods on $\mathcal{C}^1$ functionals to locally-Lipschitz functionals. We prove convergence results on these Palais-Smale sequences emerging from combinations of Laplace eigenvalues or combinations of Steklov eigenvalues in dimension 2.
\end{abstract}

\maketitle

Optimization of eigenvalues of operators (Laplacian with Dirichlet or Neumann boundary conditions, Dirichlet-to-Neumann operator, bi-laplacian, magnetic Laplacian etc) is a common field of spectral geometry.  We consider the eigenvalues as functionals depending on the shape and topology of the domain, on the operator, and/or on the geometric structure (Riemannian metrics, CR structure, sub-Riemannian metrics, etc). One old and celebrated problem was independently solved by Faber \cite{faber} and Krahn \cite{krahn} in 1923: the domains minimizing the first Laplace eigenvalue with Dirichlet boundary conditions among domains of same volume in $\mathbb{R}^n$ are Euclidean balls. This problem is very similar to the classical problem of isoperimetry, and the proof of this result uses the isoperimetric inequality, so that even when the perimeter is not involved in the renormalization (by a prescribed area/perimeter/diameter or Cheeger constant etc) of an eigenvalue functional, shape optimization on it is often called an isoperimetric problem on the eigenvalue.

\medskip

We can distinguish two main families of optimization of eigenvalues. In the first one, the ambiant geometry is prescribed (for instance, the Euclidean space $\R^n$, sphere, hyperbolic space, etc) and there is an optimization with respect to the shape and topology of a domain in this ambiant space. Emblematic results are the Faber-Krahn inequality \cite{faber}\cite{krahn} and the Szeg\"o-Weinberger \cite{szego}\cite{weinberger} inequality. In the second one, the ambiant topology is prescribed (on a fixed manifold) but the optimization holds with respect to the metric on the manifold, or potentials involved in the eigenvalue operator. An emblematic result is Hersch inequality \cite{hersch}: the round sphere is the maximizer of the first Laplace eigenvalue among metrics of same area on the $2$-sphere. In both problems, we look for bounds on eigenvalues, optimal inequalities and critical domains/metrics/potentials realizing these bounds. 

The current paper is devoted to the second family of problems. In principle, the bigger the space of variations is, the richer the critical points of the functional are. For instance, \textit{critical metrics} for combinations of Laplace eigenvalues over Riemannian metrics with prescribed volume are associated to minimal surfaces into ellipsoids (see \cite{petrides-4}), while \textit{critical metrics} for Steklov eigenvalues with prescribed perimeter are associated to free boundary minimal surfaces into ellipsoids (see \cite{petrides-5}). If only one eigenvalue appears in the functional, the target ellipsoids are spheres/balls as was primarily noticed by Nadirashvili \cite{nadirashvili} for Laplace eigenvalues and Fraser and Schoen for Steklov eigenvalues \cite{fs2}\cite{fs}. This gives an elegant connexion with the theory of minimal surfaces. If we look for critical metrics with respect to variations in a conformal class, we only obtain harmonic maps instead of minimal immersions \cite{EI4}\cite{EI5}\cite{fs2}\cite{petrides-4}\cite{petrides-5}. Other examples of critical metrics will are given in \cite{pt} thanks to a unified approach based on computations of subdifferentials (see e.g. \cite{clarke} and discussions below). Noticing that the harmonic maps enjoy a regularity theory (see e.g \cite{helein}\cite{riviere}), we can start a long story of investigations for variational aspects of eigenvalue functionals.

\medskip

In the past decades, many variational methods have been proposed since the seminal works by Nadirashvili \cite{nadirashvili} for the maximization of the first Laplace eigenvalues on tori and Fraser and Schoen \cite{fs} for the maximization of the first Steklov eigenvalues on surfaces with boundary of genus $0$. We briefly explain the idea with the example of maximization of one eigenvalue in a conformal class $[g] = \{ e^{2u}g ; u\in \mathcal{C}^\infty\left(M \right) \}$
$$ \Lambda_k(M,[g]) = \sup_{\tilde{g}\in [g]} \bar\lambda_k(\tilde{g})$$
where $\bar\lambda_k$ is a renormalized eigenvalue. Notice that conformal classes are convenient not only because the space of variation is a space of functions, but also because there are upper bounds on eigenvalues in this space \cite{Korevaar}\cite{hassann}. The main idea was to build a specific maximizing sequence of conformal factors that emerge from a \textit{regularized} variational problem. 
\begin{itemize}
\item In \cite{nadirashvili}, (Laplacian, dimension 2) the author maximizes the first eigenvalue $\bar\lambda_1$ on the smaller admissible space $E_N$ of conformal factors $f\in \mathcal{C}^{\infty}(M)$ such that $ 0 \leq f \leq N$ for $N\in \mathbb{N}$, giving a maximizing sequence as $N \to +\infty$ of $L^{\infty}$ factors $f_N \in E_N$ for $ \Lambda_1(\Sigma,[g]) = \sup_{\tilde{g}\in [g]} \lambda_1(\tilde{g})$. 
\item In \cite{fs}, \cite{petrides}, (Laplacian, dimension 2) the authors maximize a relaxed functional $f \mapsto \bar\lambda_1\left( K_\eps( f )g \right)$, where $K_\eps( f )$ is the solution at time $\eps>0$ of the heat equation with respect to $g$ at time $\eps>0$ with initial data $f$, obtaining a maximizing sequence $K_\eps(\nu_\eps)$ of smooth positive conformal factors as $\eps\to 0$, for some maximal probability measure $\nu_\eps$ of the relaxed functional $\nu \mapsto \bar\lambda_1\left( K_\eps( \nu )g \right)$.
\item In \cite{GurskyPerez}, (Conformal Laplacian, dimension $n\geq 3$) the authors proposed to modify both the functional and the space of admissible variations.
\end{itemize}
Whatever the choice, the main difficulty is to obtain convergence of this maximizing sequence of conformal factors to a regular conformal factor. Since these maximizing sequences come from the maximization of a regularized variational problem, we obtain Euler-Lagrange equations expected to bring regularity estimates on the sequence, in order to pass to the limit. Of course, these expectations are only possible if sequences of critical metrics already \textit{a priori} satisfy regularity estimates and compactness properties. This is the case for conformal factors associated to harmonic maps \cite{helein}\cite{riviere}.

The second method (see \cite{petrides}), improved in \cite{petrides-2}  and \cite{petrides-3} (Laplace and Steklov eigenvalues with higher index) is now performed for combinations of eigenvalues \cite{petrides-4} \cite{petrides-5}. The first method (see \cite{nadirashvili} \cite{nadirashvili-sire}) was improved in \cite{knpp19} for Laplace eigenvalues of higher index. It is also worth mentioning that there is an indirect method to maximize first and second conformal Laplace eigenvalues \cite{KS} \cite{KS2} based on min-max methods to build harmonic maps. While it is difficult to generalize it to higher eigenvalues or combinations, this gives a nice characterization of the maximizers, also leading to quantified inequalities on first and second eigenvalues \cite{knps21}.

\medskip

In the current paper, we simplify, unify and generalize the previous variational methods by defining a notion of \textit{Palais-Smale} (PS) sequences of conformal factors. It is a significative improvement, e.g for the following reasons:
\begin{itemize}
\item We can observe that maximizing sequences extracted by the maximization of relaxed functionals by the Heat kernel $e^{2u_\eps} = K_\eps[\nu_\eps]$ (in \cite{petrides}\cite{petrides-2}\cite{petrides-4}\cite{petrides-5}) satisfy the properties of (PS) sequences as $\eps\to 0$. Notice that these sequences $(e^{2u_\eps})_{\eps>0}$ are canonical in the sense that they satisfy even more regularity properties (for instance, there are $\mathcal{C}^0$ a priori estimates on eigenfunctions) than a random (PS) sequence.  However, working on these sequences requires an overly high technicality. 
\item All the previous methods are \textit{ad hoc} methods while the concept of (PS) sequences gives a systematic approach.
\item (PS) sequences can be extracted from min-max problems on combinations of eigenvalues, while the previous methods seem specific to maximizations, and for some of them specific to the maximization of only one eigenvalue.
\item With the extraction of (PS) sequences by the Ekeland variational principle (explained in the current paper), we can prove that all the minimizing sequences converge in some sense to a smooth optimizer.
\item This new method easily adapts to equivariant optimization problems with applications in \cite{petrides-7} and \cite{petrides-8}.
\item It is also used in \cite{petrides-9} to prove existence of a minimizer for combinations of eigenvalues of the Laplacian with respect to all the metrics for any topology (and in particular the existence of a maximizer for the first eigenvalue that was left open in general since the seminal papers by \cite{hersch} on spheres and \cite{nadirashvili} on tori)
\item It is also developped in \cite{petrides-6} for eigenvalues of the Laplacian in higher dimensions, with all the specificities due to higher dimensions. 
\end{itemize}

Classically, \textit{Palais-Smale} sequences on a $\mathcal{C}^1$ functional $E : X \to \mathbb{R}$ are sequences  $(x_n)$ such that $E(x_n) \to c$ and $\left\vert DE(x_n) \right\vert \to 0$. The main problem is that a functional involving eigenvalues (depending on a space $X$ of metrics, conformal factors, potentials, etc) is not a $\mathcal{C}^1$ functional. Of course, it is a $\mathcal{C}^1$ functional at any point in which the involved eigenvalues are simple, but we often have multiplicity of eigenvalues at the critical points, corresponding to intersection of smooth branches of eigenvalues. However, thanks to F. Clarke (see e.g \cite{clarke}), the subdifferential $\partial E(x)$ plays the role of the differential for locally Lipschitz functionals. Roughly speaking, it is a space of subgradients containing all the informations on the first variation of the functional, and in particular on the derivatives corresponding to the smooth branches of eigenvalues at points of multiplicity (see \cite{pt} for more details). Then, criticality of $E$ at $x$ can be defined by $0 \in \partial E(x)$. The current paper is devoted to quantify this property by asking a property that can be roughly written as $\left\vert \partial E(x_n) \right\vert \to 0$ for minimizing sequences, for instance thanks to the Ekeland variational principle (see for instance the nice book \cite{struwe}). We emphasize that this systematic approach is promising to solve many other variational problems on eigenvalues.


This method is explained in Section \ref{sec1}. In particular, we develop a new variational framework that is well adapted to eigenvalue functionals : we choose spaces of admissible variables that allow us to define eigenvalues, and their derivatives in order to apply the Ekeland variational principle.
As in the previous methods, the main difficulty is then to prove convergence of \textit{Palais-Smale} sequences. A wide part of the current paper is devoted to prove the convergence of minimizing sequences in a conformal class for functionals depending on combinations of Laplace eigenvalues (proof or Proposition \ref{palaissmalegeneral} in Section \ref{sec2}) or combinations of Steklov eigenvalues (proof of Proposition \ref{palaissmalegeneralsteklov} in Section \ref{sec3}) in dimension 2, that lead to Theorem \ref{theo:main}. In Section \ref{sec4}, we list $\eps$-regularity results on harmonic maps and free boundary harmonic maps into ellipsoids that are independent of the dimension of the target ellipsoids up to control there excentricity: the proof of these quite new results, first oberved in \cite{KS} or \cite{kkms} in the case of the sphere is given in \cite{petrides-10}. All along the paper, we then rewrite a proof of the main theorems in \cite{petrides-4} and \cite{petrides-5} to simplify and enlighten the techniques used there, and we prove that this convergence holds for any maximizing sequence.

\section{The variational approach} 
\label{sec1}

\subsection{The variational problem and notations}
Let $\Sigma$ be a compact surface. If $\partial \Sigma = \emptyset$, we consider for a Riemannian metric $g$ the k-th renormalized eigenvalue of the Laplacian
$$ \bar{\lambda}_k(g):= \inf_{E\in \mathcal{G}_{k+1}\left( H^1(\Sigma,g) \right)} \max_{\varphi \in E\setminus \{0\}} \frac{\int_\Sigma \vert \nabla \varphi \vert^2_g dA_g}{\int_\Sigma \varphi^2 dA_g} A_g(\Sigma) $$
where $ \mathcal{G}_{k+1}\left( H^1(\Sigma) \right)$ is the set of subspaces of $H^1(\Sigma,g)$ of dimension $k+1$, $dA_g$ is the area measure associated to $g$ and $A_g(\Sigma):= \int_\Sigma dA_g$ is the total area with respect to $g$
and if $\partial \Sigma \neq \emptyset$, we consider for a Riemannian metric $g$ the k-th renormalized Steklov eigenvalue
$$ \bar{\sigma}_k(g):= \inf_{E\in \mathcal{G}_{k+1}\left( H^1(\Sigma,g) \right)} \max_{\varphi \in E\setminus \{0\}} \frac{\int_\Sigma \vert \nabla \varphi \vert^2_g dA_g}{\int_{\partial\Sigma} \varphi^2 dL_g} L_g(\partial\Sigma) $$
where $dL_g$ is the length measure associated to the induced metric of $g$ on $\partial\Sigma$ and $L_g(\partial\Sigma):= \int_{\partial\Sigma} dL_g$ is the total length of $\partial\Sigma$ with respect to $g$.

We let $F : \left(\R_+^\star\right)^m \to \R$ be a $\mathcal{C}^1$ function such that 
$$ \forall i \in \{1,\cdots,m\}, \partial_i F(\lambda_1,\cdots,\lambda_m) \leq 0. $$
since $F$ is non-increasing with respect to all the coordinates, we can extend $F$ by continuity to $\R_+^m$ allowing the value $+\infty$ on $\R_+^m \setminus \left(\R_+^\star\right)^m$.

We set if $\partial\Sigma = \emptyset$ ("closed case")
$$E(g):= F(\bar{\lambda}_1(g),\cdots,\bar{\lambda}_m(g))$$
and if $\partial\Sigma \neq \emptyset$ ("boundary case")
$$E(g):= F(\bar{\sigma}_1(g),\cdots,\bar{\sigma}_m(g))$$
and
$$ I_F(\Sigma,[g]) := \inf_{\tilde{g} \in [g] } E(\tilde{g})  $$
where the infimum is taken on the conformal class of a metric $g$ 
$$[g]:= \{ \tilde{g} ; \exists u\in \mathcal{C}^\infty\left(\Sigma\right), \tilde{g}=e^{2u}g  \}$$

We denote $a_F$ the minimal integer such that
$$ I_F(\Sigma,[g]) < I_{F,a_F}(\Sigma,[g]) $$
where for $a \geq 1$
$$I_{F,a}(\Sigma,[g]) := \inf_{\tilde{g} \in [g] } F(0,\cdots,0,\bar{\lambda}_{a+1}(\tilde{g}),\cdots, \bar{\lambda}_m(\tilde{g}))$$
$$ \text{ or } \inf_{\tilde{g} \in [g] } F(0,\cdots,0,\bar{\sigma}_{a+1}(\tilde{g}),\cdots, \bar{\sigma}_m(\tilde{g})).$$
For instance, $a_F=1$ if $F = +\infty$ on $\{0\} \times \R_+^{m-1}$.

In the following theorem, $(\mathbb{B}^2, [ h])$ is a \textit{bubble}: the round sphere $(\mathbb{S}^2, g_{\mathbb{S}^2})$ endowed with its conformal class in the closed case and the Euclidean disk $(\mathbb{D}, g_{\mathbb{D}})$ endowed with its conformal class  in the boundary case. We will denote $dA_{\mathbb{S}^2}$ and $dA_{\mathbb{D}}$ the area measure with respect to these metrics.

\begin{theo} \label{theo:main}
For any minimizing sequence $e^{2u_n} dA_g$ (resp $e^{u_n}dL_g$ if $\partial \Sigma \neq \emptyset$) for $I_F(\Sigma,[g])$, we have up to the extraction of a subsequence that $\left(e^{2u_n} dA_g \right)$ (resp $\left(e^{u_n}dL_g\right)$) MW$\star$-bubble-tree-converges to a minimizer of $I_F(\tilde{\Sigma},[\tilde{g}])$  where
$$ (\tilde{\Sigma},[\tilde{g}]) := (\Sigma, [ g] ) \sqcup \bigsqcup_{j=1}^l (\mathbb{B}^2, [ h]) \text{ or } \bigsqcup_{j=1}^l (\mathbb{B}^2, [ h]).$$
The conformal factors of the minimizer are positive (except on a finite number of conical singularities in the closed case) and smooth.

In addition denoting $s$ the number of connected components of $\tilde{\Sigma}$ and $I$ the maximal integer such that $\bar{\lambda}_I(e^{2u_n}g) \to 0$ as $n\to +\infty$, we have $s \leq I+1 \leq a_F$.


In particular, if 
$$a_F = 1 \text{ and } I_F(\Sigma,[g]) < I_F(\mathbb{B}^2, [ h]) $$ 
then up to the extraction of a subsequence, any minimizing sequence MW$\star$-converges to a minimizer of to a positive (except on a finite number of conical singularities in the closed case) and smooth conformal factor on $\Sigma$.
\end{theo}

The definition of MW$\star$ bubble tree convergence is given in Definition \ref{defi:mwstar}. Beyond the numerous oportunities and simplifications given by the techniques that lead to this result, if we compare it to the main result of \cite{petrides-4} and \cite{petrides-5}, this result is new in the sense that the convergence holds for any minimizing sequence. This is a first step to establish stability results discussed in \cite{knps21}.

\subsection{Extension to the complete functional space of continuous bilinear maps on $H^1$}

We let $B$ be the Banach space of symmetric continuous bilinear forms $\beta : H^1(\Sigma)\times H^1(\Sigma) \to \R$ endowed with the norm
$$ \left\Vert \beta \right\Vert_g = \sup_{\varphi,\psi \in H^1 \setminus\{0\}} \frac{\left\vert \beta(\varphi,\psi) \right\vert}{\Vert \varphi \Vert_{H^1(g)} \Vert \psi \Vert_{H^1(g)}}  $$
where in the context of Laplace eigenvalues in closed surfaces
$$ \Vert \varphi \Vert_{H^1(g)}^2  := \int_{\Sigma} \vert \nabla \varphi \vert^2_g dA_g + \int_{\Sigma} \varphi^2 dA_g $$
and in the context of Steklov eigenvalues on compact surfaces with boundary,
$$ \Vert \varphi \Vert_{H^1(g)}^2  := \int_{\Sigma} \vert \nabla \varphi \vert^2_g dA_g + \int_{\partial} \varphi^2 dL_g. $$
Notice that pairs of the norms $\Vert \cdot \Vert_{H^1(g)}$ or $\Vert \cdot \Vert_{g}$ with different metrics are equivalent and that the space $B$ is independent of the choice of the metric. We denote $B_+$ the subspace of non-negative bilinear forms of $B$. Let $\beta \in B_+$. We set the $k$-th generalized eigenvalue
$$ \lambda_k(\beta) = \inf_{V \in \mathcal{G}_k(V_\beta)} \max_{\varphi \in V\setminus\{0\}} \frac{\int_{\Sigma} \left\vert \nabla \varphi \right\vert^2_g dA_g }{ \beta(\varphi,\varphi) } $$
where $\mathcal{G}_k(V_\beta)$ is the set of $k$-dimensional vector subspace of
$$ V_\beta = \{ \varphi \in \mathcal{C}^\infty(\Sigma), \beta(1,\varphi) = 0 \} $$
Notice that we can replace $V_\beta$ by its closure in $H^1$:
$$ \overline{V_\beta} =  \{ \varphi \in H^1(\Sigma), \beta(1,\varphi) = 0 \}$$
in the definition of $\lambda_k(\beta)$. Notice also that $\left[0,+\infty \right]$ is the set of admissible values of $\lambda_k$ on $B_+ $.

Finally, we set the $k$-th renormalized eigenvalue
$$ \bar{\lambda}_k(\beta) = \lambda_k(\beta) \beta(1,1). $$
and by convention $\bar{\lambda}_k = 0$ if $\beta(1,1)=0$.

\begin{prop} \label{prop:continuous}
$\lambda_k$ is an upper semi-continuous functional on
$$ G=  \{ \beta \in B_+ ; \beta(1,1)\neq 0  \} $$
and 
$\lambda_k$ and $\bar{\lambda}_k$ are locally Lipschitz maps on the open set 
$$ F =  \{\beta \in B_+ ; \beta(1,1)\neq 0 \text{ and } \lambda_k(\beta) < +\infty \} $$
Moreover, for any $\Lambda >0$,
$$  F_\Lambda =  \{\beta \in B_+ ; \bar{\lambda}_k(\beta) \leq \Lambda \}  $$
is a closed set in $B$.
\end{prop}

\begin{proof}

\textbf{Step 1: $\lambda_k$ is upper semi-continuous on $G$.}

\medskip

Let $\beta , \beta_n \in G$ such that $\beta_n \to \beta$ in $B$. If $\lambda_k(\beta)=+\infty$, then, there is nothing to prove. We assume that $\lambda_k(\beta) < +\infty$. Let $V \in \mathcal{G}_k(V_\beta)$ such that
$$ \max_{\varphi \in V\setminus\{0\}} \frac{\int_{\Sigma} \left\vert \nabla \varphi \right\vert^2_g dA_g }{ \beta(\varphi,\varphi) } \leq \lambda_k(\beta)+\delta $$
Then
$$ \lambda_k(\beta_n) \leq \max_{\varphi \in V \setminus\{0\}} \frac{\int_{\Sigma} \left\vert \nabla \varphi \right\vert^2_g dA_g}{\beta_n\left( \varphi - \frac{\beta_n\left(1,\varphi\right)}{\beta_n(1,1)},\varphi - \frac{\beta_n\left(1,\varphi\right)}{\beta_n(1,1)}  \right)} = \max_{\varphi \in V \setminus\{0\}} \frac{\int_{\Sigma} \left\vert \nabla \varphi \right\vert^2_g dA_g}{\beta_n\left( \varphi , \varphi \right) - \frac{\beta_n\left(1,\varphi\right)^2}{\beta_n(1,1)}}  $$
Let $\varphi \in V$ be such that $\Vert \varphi \Vert = 1$
$$ \beta_n\left( \varphi , \varphi \right) - \frac{\beta_n\left(1,\varphi\right)^2}{\beta_n(1,1)} \geq \beta(\varphi,\varphi) - \Vert \beta_n - \beta \Vert - \frac{\Vert \beta_n - \beta \Vert^2}{\beta(1,1) - \Vert \beta_n - \beta \Vert } . $$
Since $\lambda_k(\beta) <+\infty$, we know that $\beta(\varphi,\varphi)>0$, and that $V$ is a finite-dimensional set,
$$ \inf_{\varphi \in V, \Vert \varphi \Vert = 1} \beta(\varphi,\varphi) > 0 $$
and since $\beta(1,1)\neq 0$, and $\beta_n \to \beta$, we obtain that
$$ \lambda_k(\beta_n) \leq \lambda_k(\beta) + \delta +o(1) $$
as $n\to +\infty$. Letting $n \to +\infty$ and then $\delta \to 0$, we obtain the property.

\medskip

\textbf{Step 2: $\lambda_k$ is continous on $F$ and $F_{\Lambda}$ is closed}

\medskip

Let $\beta \in B , \beta_n \in F$ be such that $\beta_n \to \beta$ in $B$. We assume that
$$ \Lambda := \limsup_{n\to +\infty} \lambda_k(\beta_n) < + \infty. $$
Let $V_n \in \mathcal{G}_k(V_{\beta_n})$ be such that
$$ \max_{\varphi \in V_n\setminus\{0\}} \frac{\int_{\Sigma} \left\vert \nabla \varphi \right\vert^2_g dA_g }{ \beta_n(\varphi,\varphi) } \leq \lambda_k(\beta_n)+\delta \leq \Lambda + 2\delta $$
where the last inequality holds for $n$ large enough. Then
$$ \lambda_k(\beta) \leq \max_{\varphi \in V_n \setminus\{0\}} \frac{\int_{\Sigma} \left\vert \nabla \varphi \right\vert^2_g dA_g}{\beta\left( \varphi - \frac{\beta\left(1,\varphi\right)}{\beta(1,1)},\varphi - \frac{\beta\left(1,\varphi\right)}{\beta(1,1)}  \right)} = \max_{\varphi \in V_n \setminus\{0\}} \frac{\int_{\Sigma} \left\vert \nabla \varphi \right\vert^2_g dA_g}{\beta\left( \varphi , \varphi \right) - \frac{\beta\left(1,\varphi\right)^2}{\beta(1,1)}} . $$
Let $\varphi \in V_n$, then
$$ \beta\left( \varphi , \varphi \right) - \frac{\beta\left(1,\varphi\right)^2}{\beta(1,1)} \geq \beta_n(\varphi,\varphi) - \left( \Vert \beta_n - \beta \Vert - \frac{\Vert \beta_n - \beta \Vert^2}{\beta(1,1)  } \right) \left\Vert \varphi \right\Vert_{H^1}^2 . $$
 We have the following general Poincar\'e inequality (see e.g \cite{Zie89}, lemma 4.1.3]):
$$ \int_{\Sigma}\left(\varphi - \frac{\beta_n(1,\varphi)}{\beta_n(1,1)} \right)^2dA_g \leq C \left\Vert \frac{\beta_n(1,.)}{\beta_n(1,1)} \right\Vert_{H^{-1}}^2 \int_{\Sigma} \left\vert \nabla \varphi \right\vert^2 dA_g $$
so that knowing that $\varphi \in V_n$, 
\begin{equation*}
\begin{split}
 \Vert \varphi \Vert_{H^1}^2 \leq & \left( C \left\Vert \frac{\beta_n(1,.)}{\beta_n(1,1)} \right\Vert_{H^{-1}}^2 + 1  \right)\left( \lambda_k(\beta_n)+ \delta \right) \beta_n(\varphi,\varphi) \\
 \leq & \left( C \left( \frac{\left\Vert\beta\right\Vert + \Vert \beta_n - \beta \Vert}{\beta(1,1) - \Vert \beta_n - \beta \Vert } \right)^2 + 1  \right)\left(\Lambda + 2 \delta \right)\beta_n(\varphi,\varphi)
\end{split} 
  \end{equation*}
and we obtain that
$$ \lambda_k(\beta) \leq \left(\lambda_k(\beta_n) + \delta\right)(1+o(1)) $$
so that letting $n\to+\infty$ and then $\delta\to 0$, we obtain the expected result.

\medskip

\textbf{Step 3: $\lambda_k$ is locally Lipschitz on $F$}

\medskip

Let $\beta\in F$. We set $\Lambda = \lambda_k(\beta)+1$. Let $\eps_0$ and let $\beta_1,\beta_2 \in F_\Lambda \cap B(\beta,\eps_0)$ be such that
$$ \left\Vert \beta_1 - \beta_2 \right\Vert =: \eps \leq 2\eps_0 \text{ and } \sup_{B(\beta,\eps_0)} \lambda_k \leq \Lambda. $$
$\eps_0$ exists by continuity of $\lambda_k$. Let $0<\delta<1$ we shall fix later and let $V \in \mathcal{G}_k(V_{\beta_1})$ be such that
$$ \max_{\varphi \in V\setminus\{0\}} \frac{\int_{\Sigma} \left\vert \nabla \varphi \right\vert^2_g dA_g }{ \beta_1(\varphi,\varphi) } \leq \lambda_k(\beta_1)+\delta $$
Then, we test the space 
$$ \tilde{V} := \left\{ \varphi - \frac{\beta_2(1,\varphi)}{\beta_2(1,1)}  ; \varphi \in V \right\} \in \mathcal{G}_k(V_{\beta_2}) $$
in the variational characterization of $\lambda_k(\beta_2)$:
$$ \lambda_k(\beta_2) \leq \max_{\varphi \in V \setminus \{0\}}  \frac{\int_{\Sigma} \left\vert \nabla \varphi \right\vert^2_g dA_g }{ \beta_2\left(\varphi - \frac{\beta_2(1,\varphi)}{\beta_2(1,1)},\varphi - \frac{\beta_2(1,\varphi)}{\beta_2(1,1)}\right) } $$
for $\varphi\in V$, we have
\begin{equation*} 
\begin{split} \beta_2\left(\varphi - \frac{\beta_2(1,\varphi)}{\beta_2(1,1)},\varphi - \frac{\beta_2(1,\varphi)}{\beta_2(1,1)}\right) = & \beta_1(\varphi,\varphi) + \left(\beta_1-\beta_2\right)(\varphi,\varphi) - \frac{\left(\beta_2-\beta_1\right)(1,\varphi)^2}{\beta_2(1,1)} \\
\geq & \beta_1(\varphi,\varphi) - \left( \Vert \beta_1-\beta_2 \Vert + \frac{ \Vert \beta_1-\beta_2 \Vert^2 }{\beta(1,1) - 2\eps_0} \right) \Vert \varphi \Vert_{H^1}^2
\end{split}
 \end{equation*}
 We have the following general Poincar\'e inequality (see e.g \cite{Zie89}, lemma 4.1.3]):
$$ \int_{\Sigma}\left(\varphi - \frac{\beta_1(1,\varphi)}{\beta_1(1,1)} \right)^2dA_g \leq C \left\Vert \frac{\beta_1(1,.)}{\beta_1(1,1)} \right\Vert_{H^{-1}}^2 \int_{\Sigma} \left\vert \nabla \varphi \right\vert^2 dA_g $$
so that knowing that $\varphi \in V$, 
\begin{equation*}
\begin{split}
 \Vert \varphi \Vert_{H^1}^2 \leq & \left( C \left\Vert \frac{\beta_1(1,.)}{\beta_1(1,1)} \right\Vert_{H^{-1}}^2 + 1  \right)\left( \lambda_k(\beta_1)+ \delta \right) \beta_1(\varphi,\varphi) \\
 \leq & \left( C \left( \frac{\left\Vert\beta\right\Vert + 2\eps_0}{\beta(1,1) - 2\eps_0} \right)^2 + 1  \right)\left(\Lambda + \delta \right)\beta_1(\varphi,\varphi)
\end{split} 
  \end{equation*}
and gathering all the previous inequalities, we obtain
\begin{equation*}
\begin{split} \lambda_k(\beta_2) \leq & \lambda_k(\beta_1) \left(1 - \left(\eps + \frac{\eps^2}{\beta(1,1)-2\eps_0}\right)\left( C \left( \frac{\left\Vert\beta\right\Vert + 2\eps_0}{\beta(1,1) - 2\eps_0} \right)^2 + 1  \right)\left(\Lambda + \delta \right) \right)^{-1} \\
\leq & \left(\lambda_k(\beta_1)+\delta\right) ( 1 -  C_{\Lambda}(\eps_0) \eps)^{-1}
\end{split}\end{equation*}
where 
$$ C_\Lambda = \left(1 + \frac{2\eps_0}{\beta(1,1)-2\eps_0}\right)\left( C \left( \frac{\left\Vert\beta\right\Vert + 2\eps_0}{\beta(1,1) - 2\eps_0} \right)^2 + 1  \right)\left(\Lambda + 1 \right). $$
Choosing $2\eps_0 < \beta(1,1)$ such that $C_\Lambda(\eps_0)\eps_0 \leq \frac{1}{2}$, we obtain
$$ \lambda_k(\beta_2) \leq \left(\lambda_k(\beta_1)+\delta\right)\left(1 + 2 C_\Lambda(\eps_0) \eps \right) $$
Now, letting $\delta \to 0$, we obtain
$$ \lambda_k(\beta_2)- \lambda_k(\beta_1) \leq 2 \Lambda C_\Lambda(\eps_0) \Vert \beta_1-\beta_2 \Vert $$
Exchanging $\beta_1$ and $\beta_2$, the same argument leads to
$$ \left\vert \lambda_k(\beta_2)- \lambda_k(\beta_1) \right\vert \leq 2 \Lambda C_\Lambda(\eps_0) \Vert \beta_1-\beta_2 \Vert. $$
\end{proof}

We set $\overline{X}$ the closure of $X$ in $B$ where
$$ X = \left\{ (\varphi,\psi)\in H^1\times H^1 \mapsto \int_{\Sigma} e^{2u} \varphi\psi dA_g ; u \in \mathcal{C}^\infty\left( \Sigma \right) \right\} $$
if we consider the problem of Laplace eigenvalues and
$$ X = \left\{ (\varphi,\psi)\in H^1\times H^1 \mapsto \int_{\partial\Sigma}e^{u}\varphi\psi dL_g ; u\in \mathcal{C}^\infty\left(\partial\Sigma\right) \right\} $$
if we consider the problem of Steklov eigenvalues.

\medskip

We denote $Q_+$ the set of squares of $H^1$ functions and $Q = Span(Q_+)$. 
One immediate property of $\beta \in \overline{X}$ is that $\beta$ acts as a linear map on $Q$.

\begin{prop} For any $\beta \in \overline{X}$, there is a unique linear map $L_\beta : Q \to \R$ such that
$$ \forall \phi,\psi \in H^1\left(\Sigma\right), L_{\beta}\left( \phi\psi \right) = \beta\left(\phi,\psi\right) $$
and in particular
$$ \forall \phi \in H^1\left(\Sigma\right), L_{\beta}\left( \phi^2 \right) = \beta\left( \phi , \phi \right) \geq 0. $$
In addition in the closed case, $L_\beta : Q \cap \mathcal{C}^0(\Sigma) \to \R $ has a unique extension $L_\beta : \mathcal{C}^0(\Sigma) \to \R$ ($L_\beta$ is a non-negative Radon measure on $\Sigma$). In the case of compact surfaces with boundary, $L_\beta : Q \cap \mathcal{C}^0(\partial\Sigma) \to \R $ has a unique extension $L_\beta : \mathcal{C}^0(\partial\Sigma) \to \R$ ($L_\beta$ is a non-negative Radon measure on $\partial\Sigma$)
\end{prop}

\begin{proof} Let $\theta \in Q$. Let $\{ \phi_i \}_{i\in I}$ and $\{ \psi_j \}_{j\in J}$ two finite families of $H^1$ functions and $\{ t_{i} \}_{i\in I}$ and $\{ s_j \}_{j\in J}$ associated families of real numbers such that
$$\theta = \sum_{i\in I} t_i \phi_i^2 = \sum_{j\in J} s_j \psi_j^2  $$
Then it is clear that
\begin{equation} \label{eqlinearity} \sum_{i\in I} t_i \beta\left(\phi_i,\phi_i\right) = \sum_{j\in J} s_j \beta\left(\psi_j,\psi_j\right). \end{equation}
Indeed, if $e^{2u_k}$ converges to $\beta$ in $B$. 
$$ \sum_{i\in I} t_i \int_\Sigma e^{2u_k} \phi_i^2 = \sum_{j\in J} s_j \int_\Sigma e^{2u_k} \psi_j^2  $$
and letting $k\to +\infty$, we easily deduce \eqref{eqlinearity}. Then we can set a unique linear map $L_\beta : Q \to \R$ such that
$$ \forall \phi \in H^1(\Sigma), L_{\beta}(\phi^2) = \beta(\phi,\phi) $$
More generality, we compute that
$$ L_{\beta}( 4 \phi \psi) = L_{\beta}( (\phi+\psi)^2 - (\phi - \psi)^2 ) = \beta(\phi+\psi,\phi+\psi) - \beta(\phi-\psi,\phi-\psi) = 4\beta(\phi,\psi)  .$$
Finally, we have that for $\varphi \in \mathcal{C}^\infty(\Sigma)$,
$$ L_\beta(\varphi) = \beta(1,\varphi) = \lim_{k\to +\infty} \left\vert \int_\Sigma e^{2u_k} \varphi dA_g \right\vert \leq \lim_{k\to +\infty}  \int_\Sigma e^{2u_k} dA_g  \Vert \varphi \Vert_{\mathcal{C}^0}  \leq \Vert \beta \Vert \Vert \varphi \Vert_{\mathcal{C}^0} $$
and we complete the claim by the theorem of unique extension of continuous linear forms. The  case of surfaces with boundary is similar.
\end{proof}

We also obtain the immediate corollary for eigenvalues by \cite{Korevaar} and \cite{hassann}
\begin{cor}
$$ \sup_{\beta\in \overline{X} \setminus \{0\} } \bar{\lambda}_k(\beta) = \sup_{\beta\in X \setminus \{0\} } \bar{\lambda}_k(\beta) < +\infty $$
\end{cor}

We also have the very useful compactness property of bilinear forms in $\overline{X}$
\begin{prop}
Let $c,c'>0$. Let $\beta \in \overline{X}$ be such that $\beta(1,1)\neq 0$, then the image of 
$$ S_{c,c'} = \{ (\phi,\psi) \in H^1 \times H^1 ; \Vert \phi \Vert_{H^1}^2 \leq c \text{ and } \Vert \psi \Vert_{H^1}^2 \leq c'   \} $$
and of
$$ \tilde{S}_{\beta,c,c'} = \{ (\phi,\psi) \in \overline{V_{\beta}} \times \overline{V_{\beta}} ; \int_\Sigma \vert \nabla \phi \vert_g^2 \leq c \text{ and } \int_\Sigma \vert \nabla \psi \vert_g^2 \leq c'   \} $$
by $\beta$ is a compact set. More generally if $(\beta_n) \in \overline{X}$ satisfies $\beta_n \to \beta$ in $\overline{X}$ and if $(\phi_n,\psi_n) \in \tilde{S}_{\beta_n ,c,c'}$, then there is a subsequence $(\phi_{j(n)},\psi_{j(n)})$ that converges weakly to $ (\phi,\psi) \in \tilde{S}_{\beta ,c,c'}$ in $H^1 \times H^1$ and such that
$$ \beta_{j(n)}(\phi_{j(n)},\psi_{j(n)}) \to \beta(\phi,\psi) $$
as $n\to +\infty$
\end{prop}

\begin{proof} We only prove the proposition in the context of closed surfaces. The case of surfaces with boundary is similar.
We first notice that if $\phi \in \overline{V_{\beta_n}}$, then by the Poincaré inequality,
$$ \Vert \phi \Vert_{L^2(g)}^2 \leq C \left\Vert \frac{\beta_n(1,.)}{\beta_n(1,1)} \right\Vert_{H^{-1}}^2 \int_{\Sigma}\left\vert \nabla \phi \right\vert^2_g dA_g $$
so that setting $a = \left(1 + C\left(\left\Vert \frac{\beta(1,.)}{\beta(1,1)} \right\Vert_{H^{-1}}^2 +1 \right) \right) c $ and $b = \left(1 + C\left(\left\Vert \frac{\beta(1,.)}{\beta(1,1)} \right\Vert_{H^{-1}}^2 +1\right)\right) c'  $, we obtain that $\tilde{S}_{\beta_n,c,c'} \subset S_{a,b}$ for $n$ large enough.

Let $(\phi_n,\psi_n) \in H^1 \times H^1$ be such that $ \Vert \phi_n \Vert_{H^1} \leq c$ and $\Vert \psi_n\Vert_{H^1} \leq c'$. By the weak compactness of the ball of $H^1$, up to the extraction of a subsequence, we have that $\phi_n$ and $\psi_n$ weakly converge to $\phi$ and $\psi$ in $H^1$. We aim at proving that
$$ \beta_n(\phi_n,\psi_n) \to \beta(\phi,\psi) $$
as $n\to +\infty$. Let $\delta>0$. Since $\beta \in \overline{X}$, there is a smooth positive function $e^{2u}$ such that $ \left\Vert \beta - e^{2u} \right\Vert \leq \delta $. By the compact injection of $W^{1,2} \subset L^{2}(e^{2u}g)$, we have up to the extraction of a subsequence that $\psi_n \to \psi$ and $\phi_n \to \phi$ in $L^2(e^{2u}g)$ so that
$$ \int_{\Sigma} \phi_n\psi_n e^{2u} dA_g \to \int_{\Sigma} \phi\psi e^{2u}dA_g. $$
We obtain that
$$ \left\vert \beta_n(\phi_n,\psi_n) - \beta(\phi,\psi) \right\vert \leq \left\vert \int_{\Sigma} \phi_n\psi_n e^{2u} dA_g - \int_{\Sigma} \phi\psi e^{2u}dA_g \right\vert + \left(  \Vert \beta_n - \beta \Vert + 2 \Vert \beta - e^{2u} \Vert \right) c c' $$
so that passing to the limit as $n \to +\infty$,
$$ \limsup_{n\to +\infty} \left\vert \beta(\phi_n,\psi_n) - \beta(\phi,\psi) \right\vert \leq \delta c c' $$
and letting $\delta\to 0$, we obtain the expected result.
\end{proof}

Notice also that the norms $N_{\beta}(\phi)^2 := \int_\Sigma \left\vert \nabla \phi \right\vert^2_g +\beta(\phi,\phi)$ satisfy for $\beta \in \overline X$ the existence of an open neighborhood $U_{\beta}$ and a constant $C_{\beta}$ such that
$$ \forall \beta \in U_{\beta}, \forall \phi \in H^1, C_{\beta}^{-1} N_{\tilde{\beta}}(\phi)^2 \leq N_{\beta}(\phi)^2 \leq C_{\beta} N_{\tilde{\beta}}(\phi)^2 $$

By \cite{pt}, we obtain from this compactness property that the spectrum associated to $\beta\in \overline{X}$ is discrete, that is
$$ 0 = \lambda_0 < \lambda_1(\beta) \leq \lambda_2(\beta) \leq \cdots \leq \lambda_k(\beta) \to +\infty \text{ as } k\to +\infty $$ 
and in particular that the multiplicity of eigenvalues is finite and that there are eigenfunctions. Notice that if $\Sigma$ is connected, $\lambda_0 = 0$ is a simple eigenvalue associated to the constant functions. We denote the equations on eigenfunctions $\Delta_g \varphi = \lambda \beta\left(\varphi,\cdot \right)$ if we consider $\Delta_g \varphi$ as a linear form on $H^{1}$. This notation also holds in the Steklov case.

As soon as $\beta$ belongs to the interior of $\overline{X}$, we can also compute the directional derivatives, the generalized directional derivatives, the subdifferential and the Clarke subdifferential of 
$$ \beta \mapsto F(\bar{\lambda}_1(\beta),\cdots,\bar{\lambda}_m(\beta)) $$
where $F : \left(\mathbb{R}_{+}^\star \right)^m \to \mathbb{R}_{+}^\star $ such that $\partial_i F \leq 0$ for any $i$. 
$$ \partial E(\beta) \subset co \left\{ \sum_{i=1}^m \partial_i F(\lambda_1(\beta),\cdots,\lambda_m(\beta)) \lambda_i(\beta) \left((\phi_i,\phi_i) - (1,1) \right) ; (\phi_1,\cdots,\phi_m)\in \mathcal{O}_m(\beta) \right\}$$
where $\mathcal{O}_m(\beta)$ is the set of orthonormal families with respect to $\beta$ $(\phi_1,\cdots,\phi_m)$ such that $\phi_i$ is an eigenfunction associated to $\lambda_i(\beta)$. 

In our case, we will compute right directional derivatives on points $\beta \in \overline{X}$ that do not belong to the interior of $\overline{X}$ but it is not a problem if the variation $\beta+ t b$ still belongs to the admissible set as soon as $t \searrow 0$. For the sake of completeness, we write the method of \cite{pt} in our context:

We denote by 
$$ i(k) := \min\{ i \in \mathbb{N}^* ; \lambda_i = \lambda_k \}$$
$$ I(k) := \max\{ i \in \mathbb{N}^* ; \lambda_i = \lambda_k \}$$

\begin{prop} \label{prop:firstderivative}
For $\beta \in  \bar{X}$, and $b \in  \bar{X}$,
\begin{equation} \label{eq:firstderivativeminmax}
\begin{split} \lim_{t \searrow 0 } \frac{\bar{\lambda}_k(\beta+tb) - \bar{\lambda}_k(\beta)}{t} = & \bar{\lambda}_k(g,\beta) \left( b(1,1) - \min_{ V \in \mathcal{G}_{k-i(k)+1}(E_k(\beta)) } \max_{\phi \in V\setminus \{0\}} \frac{b(\phi,\phi)}{\beta(\phi,\phi)} \right) \\ 
= & \bar{\lambda}_k(g,\beta) \left( b(1,1) - \max_{ V \in \mathcal{G}_{I(k)-k+1}(E_k(\beta)) } \min_{\phi \in V\setminus \{0\}}   \frac{b(\phi,\phi)}{\beta(\phi,\phi)} \right)
\end{split}
\end{equation}
\end{prop}

\begin{proof} The right-hand terms are equal as a consequence of the min-max formula for the quotients of a quadratic form by a positive definite quadratic form on finite-dimensional spaces. Notice that from Proposition \ref{prop:continuous}, we have that $\lambda_k(\beta+tb) \to \lambda_k(\beta)$ as $t \searrow 0$. 

We denote by 
$$\phi_{i(k)}^t,\cdots, \phi_{I(k)}^t$$
a family of $\beta$-orthonormal eigenfunctions associated to the eigenvalues 
$$\lambda_{i(k)}(\beta+tb) \leq \cdots \leq \lambda_{I(k)}(\beta+tb)$$
 we rename $\lambda_{i(k)}^t \leq \cdots \leq \lambda_{I(k)}^t$ that all converge to $\lambda_k := \lambda_k(\beta)$ as $t\to 0$.
Up to the extraction of a subsequence as $t \to 0$, $\phi_i^t$ converges to $\phi_i$ weakly in $H^1$, and 
$$ \left(\beta+tb\right)(\phi-\phi_i^t,\phi-\phi_i^t) \to 0$$
 as $t\to 0$. Passing to the weak limit on the equation satisfied by $\phi_i^t$ and to the strong limit on $\left(\beta+tb\right)(\phi_i^t,\phi_j^t) = \delta_{i,j}$, we obtain
$$ \Delta_g \phi_i = \lambda_k \beta( \phi_i,\cdot) \text{ and } \beta(\phi_i,\phi_j) = \delta_{i,j} $$
for $i(k)\leq i,j\leq I(k)$. Integrating the equation with respect to $\phi_i$ proves that
\begin{equation*}
\begin{split} \int_\Sigma \vert \nabla \phi_i \vert_g^2 dA_g & =  \lambda_k \beta(\phi_i,\phi_i) 
=  \lim_{t\to 0}  \lambda_i^t \beta(\phi_i^t,\phi_i^t) = \lim_{t\to 0} \int_\Sigma \vert \nabla \phi_i^t \vert_g^2 dA_g 
\end{split} \end{equation*}
so that $\phi_i^t$ converges strongly in $H^1$.

\medskip 

For $i(k) \leq i \leq I(k)$. We set $R_i^t := \phi_i^t - \pi_k(\phi_i^t) $ where for $v\in H^1$
$$ \pi_k(v) := v - \sum_{i=i(k)}^{I(k)} \beta(v,\phi_i) \phi_i $$
is the orthogonal projection of $v$ on $E_k(g,\beta)$ with respect to $\beta$. We have
$$ \Delta_g R_i^t  - \lambda_k \beta( R_i^t,\cdot) = \lambda_i^t (\beta + tb)( \phi_i^t,\cdot) - \lambda_k \beta( \phi_i^t,\cdot) = (\lambda_i^t - \lambda_k) \beta( \phi_i^t,\cdot) + \lambda_i^t tb( \phi_i^t,\cdot).  $$
We set
\begin{equation} \label{eq:defialphait} \alpha_i^t := \left\vert \lambda_i^t - \lambda_k \right\vert + t + \sqrt{\beta( R_i^t ,R_i^t)} \end{equation}
and
\begin{equation} \label{eq:defitilderit} \tilde{R}_i^t = \frac{R_i^t}{\alpha_i^t} \hspace{5mm} \tau_i^t = \frac{t}{\alpha_i^t} \hspace{5mm} \delta_i^t := \frac{\lambda_i^t - \lambda_k}{\alpha_i^t}. \end{equation}
Let's prove that $\tilde{R}_i^t $ is bounded in $H^1$. Let $v \in H^1$, we have that
$$ \int_M \nabla \tilde{R}_i^t \nabla v dA_g = \lambda_k \beta(\tilde{R}_i^t,v) + \delta_i^t \beta(\phi_i^t,v) + \lambda_i^t b(\phi_i^t,v) $$
so that
$$ \left\vert \int_\Sigma \nabla \tilde{R}_i^t \nabla v dA_g + \beta( \tilde{R}_i^t, v) \right\vert \leq \left( \left(\lambda_k + 1\right) \sqrt{\beta(\tilde{R}_i^t,\tilde{R}_i^t )} \Vert \beta \Vert  + \left(\delta_i^t \Vert \beta \Vert + \lambda_i^t \Vert b \Vert \right) \Vert  \phi_i^t \Vert_{H^1}  \right) \Vert v \Vert_{H^1}$$
so that by the Riesz theorem associated to the Hilbert norm $N_\beta$, and the equivalence of the $H^1$ norm and the $N_\beta$ norm, we obtain that
$ \tilde{R}_i^t $ is bounded with respect to $N_\beta$ as $t\to 0$.
By equivalence between the $H^1$ norm and the norm $N_\beta$, again, $ \tilde{R}_i^t $ is bounded in $H^1$. Then, up to the extraction of a subsequence as $t\to 0$,
$$ \tilde{R}_i^t \to \tilde{R}_i \text{ weakly in } H^1 \hspace{5mm} \tau_i^t \to \tau_i \hspace{5mm} \delta_i^t \to \delta_i.$$
Passing to the weak limit in the equation, we obtain
\begin{equation} \label{eq:limritilde} \Delta_g \tilde{R}_i - \lambda_k \beta(  \tilde{R}_i ,\cdot)= \delta_i \beta( \phi_i,\cdot ) + \tau_i \lambda_k b( \phi_i,\cdot). \end{equation}
In addition, up to the extraction of a subsequence, 
$$\beta(\tilde{R}_i^t - \tilde{R}_i,\tilde{R}_i^t - \tilde{R}_i) \to 0$$
as $t\to 0$ and we obtain because of the definitions \eqref{eq:defialphait} and \eqref{eq:defitilderit}
\begin{equation} \label{eq:sumRitildedeltatau} \beta(\tilde{R}_i,\tilde{R}_i) + \vert \delta_i \vert + \tau_i = 1 \end{equation}
Now, we integrate \eqref{eq:limritilde} against $\phi_i$ and we obtain that
\begin{equation} \label{eq:eqriintegrated}
\delta_i \beta(\phi_i\phi_i) + \tau_i \lambda_k  b( \phi_i,\phi_i ) = 0.
\end{equation}
Now, we assume by contradiction that $\tau_i = 0$, then by \eqref{eq:eqriintegrated}, $\delta_i=0$ and by \eqref{eq:limritilde}, $\tilde{R}_i \in E_k(\beta) \cap E_k(\beta)^{\perp_{Q(\beta,\cdot)}} = \{0\}$. This contradicts \eqref{eq:sumRitildedeltatau}. Therefore $\tau_i \neq 0$ and
$$ \frac{\delta_i}{\tau_i} = \frac{-\lambda_k(\beta) b(\phi_i,\phi_i) }{\beta(\phi_i,\phi_i)} $$
Integrating \eqref{eq:limritilde} against $\phi_j$ for $j\neq i$, we obtain that
$$ \lambda_k(\beta) b( \phi_i, \phi_j ) = 0 $$
so that $\phi_{i(k)},\cdots,\phi_{I(k)}$ are nothing but an orthonormal basis with respect to $\beta$ that is orthogonal with respect to $ -\lambda_k(\beta) b $. Since in addition we have that $\delta_{i(k)}\leq \cdots \leq \delta_{I(k)}$, classical min-max formulae for orthonormal diagonalization give
$$ \frac{\delta_i}{\tau_i} = \min_{ V \in \mathcal{G}_{i-i(k)+1}(E_k(\beta)) } \max_{v \in V\setminus \{0\}} \frac{-\lambda_k(\beta) b(v,v)}{\beta(v,v)} $$
Since the right-hand term is independent of the choice of the subsequence as $t\to 0$, we obtain that the directional derivative exists and
$$ \lim_{t\searrow 0} \frac{\lambda_i^t-\lambda_i}{t} = \lim_{t\searrow 0} \frac{\delta_i^t}{\tau_i^t} = \frac{\delta_i}{\tau_i} $$
and a straightforward chain rule using the directional derivative of $(\beta + t b)(1,1)$  completes the proof of the proposition.
\end{proof}

\subsection{Regularization of minimizing sequences by Ekeland's variational principle}
The familly of functionals $E$ depending on $F$ given in the beginning of the section can be extended to $\overline{X}$. We obtain the following proposition for the extraction of Palais-Smale sequences $(PS)_K$ (see Definition \ref{def:PS})

\begin{prop} \label{prop:constructPSK}
For any $\eps>0$, we let $e^{2u_\eps}$ be a conformal factor, and $g_\eps:= e^{2u_\eps}g$ such that
$$ E(e^{2u_\eps}dA_g) \leq \inf_{\beta\in \overline{X}} E(\beta) + \eps^2. $$ 
(or replace by $E(e^{u_\eps}dL_g)$ in the Steklov case). Then, there is $K \leq m$ and a $(PS)_K$ sequence $(\beta_\eps,\Phi_\eps,g_\eps)$ as $\eps\to 0$.
\end{prop}

\begin{defi} \label{def:PS} Let $(\Sigma,g)$ be a compact Riemannian surface. Let $\beta_\eps \in \overline{X}$ (the definition of $\overline{X}$ depends if $\partial\Sigma =\emptyset$ or not), $\Phi_{\eps} : \Sigma \to \mathbb{R}^{m_\eps}$ be a sequence of maps with $(m_\eps)_{\eps>0} \in \left(\mathbb{N}^*\right)^{\R_+^*}$, $g_\eps := e^{2u_\eps} g$ a family of metrics conformal to $g$ and $K \in \mathbb{N}^\star$. We say that $(\beta_\eps,\Phi_\eps,g_\eps)$ satisfies the Palais-Smale assumption (with eigenvalue indices bounded by $K$) $(PS)_{K}$ as $\eps\to 0$, if
\begin{itemize}
\item The diagonal terms of $\Lambda_{\eps} := diag(\lambda_{1}^{\eps},\cdots, \lambda_{m_{\eps}}^{\eps})$ are the $m_\eps$ first (Laplace if $\partial\Sigma=\emptyset$, Steklov if $\partial\Sigma\neq \emptyset$) eigenvalues associated to $\beta_\eps$ such that $\lambda_{1}^{\eps}\leq \cdots \leq \lambda_{m_{\eps}}^{\eps} = \lambda_K^\eps$ where $\lambda_K^\eps$ is the $K$-th eigenvalue.
\item $ \Delta_g \Phi_{\eps} = \beta_\eps\left( \Lambda_{\eps} \Phi_{\eps} , .\right)$, where $\beta_\eps\left( \Lambda_{\eps} \Phi_{\eps} , . \right) : H^1(\Sigma)^{m_\eps} \to \R$
\item $L_\eps(1) = L_\eps\left( \left\vert \Phi_{\eps} \right\vert_{\Lambda_{\eps}}^2 \right) = \int_{\Sigma}\left\vert \nabla\Phi_{\eps} \right\vert_g^2 dA_g =  1$ where we denote $L_\eps$ the linear form associated to $\beta_\eps$
\item For $i\in \{1,\cdots,m_\eps\}$, $t_i^{\eps} = L_\eps \left( \left(\phi_i^\eps\right)^2 \right)$ and $\sum_{i=1}^{m_\eps} \lambda_{i}^\eps t_i^\eps  = 1$ and we have that $\lambda_i^\eps t_i^\eps \to 0$ as $\eps \to 0$ for any $i \in \{1,\cdots,m_\eps\}$ such that $\lambda_i^\eps \to 0$ as $\eps \to 0$.
\item $\left\vert \Phi_\eps \right\vert_{\Lambda_{\eps}}^2 \geq_{a.e} 1 - \theta_\eps^2$ in $\Sigma$ if $\partial\Sigma =\emptyset$ and  $\left\vert \Phi_\eps \right\vert_{\Lambda_{\eps}}^2 \geq_{a.e} 1 - \theta_\eps^2$ in $\partial\Sigma$ if $\partial\Sigma\neq \emptyset$ where $\Vert \theta_{\eps} \Vert_{H^1(g_\eps)}^2 \leq \eps$
\item $\Vert \beta_\eps - e^{2u_\eps}dA_g \Vert_{g_\eps} \leq \eps$ if $\partial\Sigma=\emptyset$, $\Vert \beta_\eps - e^{u_\eps}dL_g \Vert_{g_\eps} \leq \eps$ if $\partial\Sigma\neq\emptyset$.
\end{itemize}
\end{defi}

\begin{rem} We proved in \cite{pt} (lemma 2.1) that up to transformations coming from linear algebra, $\Phi_\eps$ can be chosen as an orthogonal family with respect to $\beta_\eps$. These transformations do not affect the other properties of $(PS)_K$ sequences. However, this extra property is not necessary in the current paper.
\end{rem}

\begin{proof}
We assume up to a renormalization that
$$ \int_\Sigma e^{2u_\eps}dA_g=1 \text{ if } \partial\Sigma = \emptyset \text{ and } \int_\Sigma e^{u_\eps}dL_g=1 \text{ if } \partial\Sigma \neq \emptyset. $$
By Ekeland's variational principle, knowing that $\{ \beta \in \overline{X} ; \beta(1,1)\geq 1 \}$ endowed with the distance $d_\eps(\beta_1,\beta_2) = \Vert \beta_1 - \beta_2 \Vert_{g_\eps}$ where $g_\eps = e^{2u_\eps} g$ and if $\partial \Sigma = \emptyset$,
$$ \Vert b \Vert_{g_\eps} := \sup_{\varphi,\psi \in H^1\setminus\{0\}} \frac{\vert \beta(\varphi,\psi) \vert}{\Vert \varphi\Vert_{H^1(g_\eps)} \Vert \psi \Vert_{H^1(g_\eps)}} $$
is a complete space as a closed subset of $\overline{X}$, we obtain the existence of $\beta_\eps \in  \overline{X}$ with $1 \leq \beta_\eps(1,1)\leq 1+ \eps$ such that 
$$ E(\beta_\eps) \leq \inf_{\beta\in \overline{X}} E(\beta) + \eps^2 $$
and
$$ \Vert \beta_\eps - e^{2u_\eps} dA_g \Vert_{g_\eps} \leq \eps \text{ if } \partial \Sigma = \emptyset \text{ and } \Vert \beta_\eps - e^{u_\eps} dL_g \Vert_{g_\eps} \leq \eps \text{ if } \partial \Sigma \neq \emptyset $$
and
$$ \forall \beta \in \bar{X}, E(\beta_\eps) - E(\beta) \leq \eps \Vert \beta_\eps - \beta \Vert_{g_\eps}. $$
In particular, we have that for any $b \in \overline{X}$ 
$$ \lim_{t\downarrow 0} \frac{E(\beta_\eps) - E(\beta_\eps + t b)}{t} \leq \eps \Vert b \Vert_{g_\eps} $$
where we know that this limit exists by the previous subsection. Without loss of generality, we can assume that all the previous inequalities hold with $\beta_\eps(1,1)=1$.
Let $V \in {L}^{2}(\Sigma)$ such that $V\geq_{a.e} 0$. Then there is $\left(\phi_1,\cdots,\phi_m\right) \in \mathcal{O}_m(\beta_\eps)$ such that
$$ \int_{\Sigma}  \left( \sum_{i=1}^m t_i^\eps \lambda_i(\beta_\eps) \left( \phi_i^2 - 1 \right) \right) V dA_{g_\eps} \geq -\eps \Vert V dA_{g_\eps} \Vert_{\eps}$$
where $t_i^\eps = - \partial_i F(\lambda_1(\beta_\eps),\cdots,\lambda_m(\beta_\eps)) \geq 0$. We also have the existence of $\theta_\eps \in W^{1,2}$ such that $\Vert \theta_\eps \Vert_{H^1(g_\eps)} = 1$ and 
$$ \Vert V dA_{g_\eps} \Vert_{g_\eps} = \int_{\Sigma} V \theta_\eps^2 dA_{g_\eps}  = \max_{\Vert \phi \Vert_{H^1(g_\eps)} = 1 } \int_{\Sigma} V \phi^2 dA_{g_\eps} 
$$
Indeed the supremum in the definition of the norm of $V$ in the vector space $B$ is realized because of the compact embedding $H^1 \subset L^p$ for any $1\leq p < +\infty$. 
We obtain that
for any $V \in L^2$ such that $V\geq_{a.e} 0$, there is $\left(\phi_1,\cdots,\phi_m\right) \in \mathcal{O}_m(\beta_\eps)$ and $\theta \in H^1$ with $\Vert \theta_\eps \Vert_{H^1(g_\eps)} \leq 1$ such that
\begin{equation} \label{eqconsequenceekeland} \int_{\Sigma} V \left( \left( \sum_{i=1}^m t_i^\eps \lambda_i(\beta_\eps) \left( \phi_i^2 - 1 \right) \right) + \eps\theta^2 \right) dA_{g_\eps} \geq 0 \end{equation}
Now, let's give a Hahn-Banach separation argument. We first notice that the set
$$ \{ \theta^2 ; \theta \in H^1 \text{ and } \Vert \theta \Vert_{H^1(g_\eps)} \leq  1 \} $$
is a compact convex subset of $L^p(\Sigma)$ for any $1\leq p < +\infty$.
Indeed, we just have to prove that it is a convex set by the compact embedding $H^1 \subset L^p$. Let $\theta_1, \theta_2 \in H^1$ such that $\Vert \theta_i \Vert_{H^1(g_\eps)} \leq 1$ for $i=1,2$. Let $t \in [0,1]$. We aim at proving that $\theta:= \sqrt{(1-t)\theta_1^2 + t \theta_2^2} \in H^1(g_\eps)$ and satisfies $\Vert \theta \Vert_{H^1(g_\eps)} \leq 1$:
$$ \int_{\Sigma} \theta^2 dA_{g_\eps} = (1-t)\int_{\Sigma} \left(\theta_1\right)^2dA_{g_\eps} + t \int_{\Sigma} \left(\theta_2\right)^2dA_{g_\eps} $$
and since $(x_1,x_2)\mapsto \sqrt{(1-t)x_1^2 + t x_2^2}$ is a Lipschitz map, $\theta \in H^1(g_\eps)$ and by the computation
\begin{equation*} 
\begin{split}\left\vert \nabla \theta \right\vert^2_{g_\eps} =_{a.e} & \left\vert \frac{(1-t)\theta_1 \nabla \theta_1 + t \theta_2 \nabla \theta_2 }{\theta} \right\vert^2_{g_\eps} \\
= & \frac{ (1-t)^2 \theta_1^2 \left\vert \nabla \theta_1 \right\vert^2_{g_\eps} + 2t(1-t) \theta_1 •\theta_2 \left\langle \nabla \theta_1\nabla \theta_2 \right\rangle_{g_\eps} + t^2 \theta_2^2 \left\vert \nabla \theta_2 \right\vert^2_{g_\eps}  }{ (1-t)\theta_1^2 + t \theta_2^2  } \\
\leq & \frac{ (1-t)^2 \theta_1^2 \left\vert \nabla \theta_1 \right\vert^2_{g_\eps} + t(1-t)( \theta_2^2 \left\vert \nabla \theta_1 \right\vert^2_{g_\eps} + \theta_1^2 \left\vert \nabla \theta_2 \right\vert^2_{g_\eps}) + t^2 \theta_2^2 \left\vert \nabla \theta_2 \right\vert^2_{g_\eps} }{(1-t)\theta_1^2 + t \theta_2^2} \\
 = & (1-t) \left\vert \nabla \theta_1 \right\vert^2_{g_\eps} + t \left\vert \nabla \theta_2 \right\vert^2_{g_\eps}
\end{split}\end{equation*}
we obtain
$$ \Vert \theta \Vert_{H^1(g_\eps)}^2 \leq (1-t)\Vert \theta_1 \Vert_{H^1(g_\eps)}^2 + t \Vert \theta_2 \Vert_{H^1(g_\eps)}^2 \leq 1 $$
which is the expected result.

Therefore, the set
$$ K = co \left\{  \left( \sum_{i=1}^m t_i \lambda_i(\beta) \left( \phi_i^2 - 1 \right) \right) + \eps\theta^2  ;  \theta \in H^1 , \Vert \theta \Vert_{H^1(g_\eps)} \leq  1 , (\phi_1,\cdots,\phi_m)\in \mathcal{O}_m(\beta_\eps) \right\} $$
is a compact convex subset of $L^p$ for any $1\leq p < +\infty$ and
$$ F = \{ f \in L^2 ; f\geq_{a.e} 0 \} $$
is a closed cone in $L^2$. We assume by contradiction that $F \cap K = \emptyset$.
Then, there is $V \in L^2$ such that
$$ \forall \psi \in K ; \int_{\Sigma} V \psi dA_{g_\eps} \leq - \alpha < 0 $$
$$ \forall f \in F ; \int_{\Sigma} V f dA_{g_\eps} \geq 0 $$
and we deduce from the second property that $V \geq_{a.e} 0$. 
From the first property is then a contradiction with \eqref{eqconsequenceekeland}. Then $F\cap K \neq \emptyset$
and there is $J_\eps$, $s_j$ for $j\in \{1,\cdots,J_\eps\}$ such that $\sum_{j=1}^{J_\eps} s_j = 1$, $(\phi^\eps_{1,j},\cdots,\phi^\eps_{m,J_\eps}) \in \mathcal{O}_m(\beta_\eps)$ and $\theta_\eps \in H^1$ such that $\Vert \theta_\eps \Vert_{H^1(g_\eps)} \leq 1$ and
$$ \sum_{j=1}^{J_\eps} s_j^\eps \sum_{i=1}^{m} t_i^\eps \lambda_i(\beta_\eps) \left( \left(\phi^\eps_{i,j}\right)^2 - 1 \right)  + \eps\theta_\eps^2 \geq_{a.e} 0 $$

We now prove that $\lambda_i^\eps t_i^\eps \to 0$ as $\eps \to 0$ for any $i \leq I$, that is $i$ such that $\lambda_i^\eps \to 0$ as $\eps \to 0$. For any $i \leq I$, we have
$$ I_{F,i}(\Sigma,[g]) < +\infty  $$
since $(\lambda_1^\eps,\cdots,\lambda_m^\eps)$ corresponds to a minimizing sequence. Then for $x_l,,\cdots,x_m >0$
$$ F(0,\cdots,0,x_{i+1},\cdots,x_m) = F(0,\cdots 0,x_i,\cdots,x_m) - \int_{0}^{x_i} \frac{t \partial_i F(0,\cdots 0,t,x_{i+1},\cdots,x_K)}{t}dt $$
implies that $\lim_{t\to 0} t \partial_i F(0,\cdots 0,t,x_{i+1},\cdots,x_m) = 0$. This implies that $t_i^\eps \lambda_i^\eps \to 0$ as $\eps \to 0$.

Finally as noticed in the remark after the proposition \cite{pt} (lemma 2.1) allows us to conclude the proof of the proposition up to replace $(t_1^\eps,\cdots,t_m^\eps)$ by an element of $Mix(t_1^\eps,\cdots,t_m^\eps)$  (see notations in \cite{pt} (lemma 2.1)).
\end{proof}

\section{Convergence of regularized minimizing sequences in the closed case}
\label{sec2}

We aim at proving the following proposition (see definition \ref{defi:mwstar} for the MW$\star$ bubble tree convergence)

\begin{prop} \label{palaissmalegeneral}
Let $(\Sigma,g)$ be a Riemannian surface without boundary and $(\beta_\eps, \Phi_\eps,g_\eps)$, be a $(PS)_{K}$ sequence as $\eps\to 0$.
Then, up to the extraction of a subsequence 
$e^{2u_\eps}dA_g$ and $\beta_\eps(1,.)$ MW$\star$--bubble tree converge to the measures $V_0 dA_g$ (possibly $0$ if $l\geq 1$) on $\Sigma$ and $V_j dA_{\mathbb{S}^2}$ on $\left(\mathbb{S}^2\right)_j$ where $V_0,V_1,\cdots,V_l$ are $L^{\infty}$ densities.

If in addition $(\beta_\eps)$ $(g_\eps)$ are minimizing sequences for $E$, then denoting 
$$\Lambda := diag\left(\bar{\lambda}_1(\widetilde{\Sigma},VdA_{\tilde{g}}),\cdots,\bar{\lambda}_n(\widetilde{\Sigma},VdA_{\tilde{g}}) \right)$$
where $\widetilde{\Sigma} = \Sigma \sqcup \bigsqcup_{j=1}^l (\mathbb{S}^2)_j$ endowed with $\widetilde{g}$ equal to $g$ on $\Sigma$ and the round metric on the copies of $\mathbb{S}^2$ and $V = V_0$ in $\Sigma$ and $V = V_j$ in $(\mathbb{S}^2)_j$, we have that
$$ V_0 = \frac{\left\vert \nabla \Phi_0 \right\vert_{\Lambda,g}^2}{\left\vert \Lambda \Phi_0 \right\vert^2}  \text{ and } V_j = \frac{\left\vert \nabla \Phi_j \right\vert_{\Lambda,g_{\mathbb{S}^2}}^2}{\left\vert \Lambda \Phi_j \right\vert^2}  $$  
where $\Phi_0 : \Sigma \to \mathcal{E}_\Lambda$ and $\Phi_j : \left(\mathbb{S}^2\right)_j \to \mathcal{E}_\Lambda$ are harmonic maps into $\mathcal{E}_\Lambda:= \{\vert x \vert^2_\Lambda = 1\}$ and we have that
$$ I_F(\Sigma,[g]) = I_F(\widetilde{\Sigma},[\widetilde{g}]) $$
\end{prop}

\begin{rem}
Notice that by a glueing method similar to \cite{ces}, we always have 
$$ I_F(\Sigma,[g]) \leq  I_F(\widetilde{\Sigma},[\widetilde{g}]) $$
and if we know that the inequality is strict, then we automatically deduce that all the minimizing sequences for $I_F(\Sigma,[g])$ MW$\star$ converge to a measure absolutely continuous with respect to $dA_g$ with a smooth density ($l=0$ in the proposition)
\end{rem}

This proposition and the remark proves Theorem \ref{theo:main} in the case of Laplace eigenvalues, noticing that $V \tilde{g}$ is a smooth metric up to conical singularities which correspond to the zeros of $V$ or of the energy densities of harmonic maps.

\subsection{Tree of concentration points }

We define MW$\star$--bubble tree convergence of sequences of $\overline{X}$ by weak-star convergence in the sense of measures in multiple scales. Here we say that a measure MW$\star$ converges if there is a weak $\star$ convergence in the set of Radon measures standing as the dual $\left(\mathcal{C}^0(\Sigma)\right)^\star$.

\begin{defi} \label{defi:mwstar}
Let $(\Sigma,g)$ be a Riemannian surface without boundary. We say that a sequence $(\mu_n)$ of positive Radon measures MW$\star$-bubble-tree converges if there is $l\in \mathbb{N}$ such that for $1\leq j \leq l$ there are sequences of points  $x_j^n \in \Sigma$ and of scales $\alpha_j^n >0$ satisfying for all $0 \leq i\neq j\leq l$
$$ \frac{d_g(x_i^n,x_j^n)}{\alpha_i^n +\alpha_j^n} + \frac{\alpha_i^n}{\alpha_j^n} + \frac{\alpha_j^n}{\alpha_i^n} \to +\infty \text{ and } \alpha_i^n \to 0 \text{ and } \alpha_j^n \to 0 $$
as $n\to +\infty$ such that 
\begin{itemize}
\item $\mu_0^n:=\mu_n$ MW$\star$ converges to $\nu_0$ in $\Sigma$.
\item for $\varphi \in \mathcal{C}^0_c(\mathbb{R}^2)$, we set
$$ \mu_j^n(\varphi) = \mu_n\left( \varphi\left( \frac{x-x_{j}^n}{\alpha_j^n}\right)  \right) $$
and $\mu_j^n$ MW$\star$ converges to $\nu_j$ in $\mathbb{R}^2$
\end{itemize}
In addition, letting $Z_j$ be the set of concentration points of $\mu_j^n$, the sets $Z_j$ are finite and
$$ \lim_{n\to +\infty} \mu_n(\Sigma) = \int_{\Sigma\setminus Z_0} d\nu_0 + \sum_{i=1}^l \int_{\mathbb{R}^2\setminus Z_i} d\nu_i  \text{ and } \forall i\in\{1,\cdots,l\}, \int_{\mathbb{R}^2\setminus Z_i} d\nu_i \neq 0.$$
Denoting 
$$\mu_0 := \nu_0 - \sum_{x\in Z_0} \nu_0(\{x\})\delta_x $$
a Radon measure of $\Sigma$ and for $1\leq j \leq l$
$$ \mu_j := \pi_{\mathbb{S}^2}^{\star} \left( \nu_j - \sum_{x\in Z_j} \nu_j(\{x\})\delta_x \right) $$
where $\pi_{\mathbb{S}^2} : \mathbb{S}^2 \to \mathbb{R}^2$ is the stereographic projection, we say that $(\mu_n)$ MW$\star$ bubble tree converges to the measure $\mu$ on $\widetilde{\Sigma} := \Sigma \sqcup \bigsqcup_{j=1}^l \left(\mathbb{S}^2\right)_j$ such that $\mu$ is equal to $\mu_0$ on $\Sigma$ and $\mu_j$ on $\left(\mathbb{S}^2\right)_j$ for $1\leq j \leq l$.
%
\end{defi}

Notice that in this definition, there is a slight abuse of notations with the use of sums of points on a manifold. This ambiguity is solved by the use of an atlas of conformal charts. More generally, in the following, all disk $\mathbb{D}_r(p)$ that appears in our analysis will correspond to a flat disk in a conformal chart of the manifold.

Before going deeply in the analysis of Palais-Smale sequences, we prove the following upper semi-continuity of eigenvalues with respect to MW$\star$ bubble tree convergence of measures. This generalizes the proof of Kokarev \cite{kokarev} for semi-continuity with respect to MW$\star$ convergence.

\begin{cl}
We assume that $(\Sigma,g,\mu_n)$ MW$\star$ bubble tree converges to $(\widetilde{\Sigma},\widetilde{g},\mu)$ as $n \to +\infty$. Then
$$ \limsup_{n\to +\infty}\lambda_k(\Sigma,g,\mu_n) \leq \lambda_k(\widetilde{\Sigma},\tilde{g},\mu) $$
\end{cl}

\begin{proof} If $\lambda_k(\widetilde{\Sigma},\tilde{g},\mu) = +\infty$, there is nothing to prove. We assume that $\lambda_k(\widetilde{\Sigma},\tilde{g},\mu) < +\infty$.
Let $\delta >0$ we will let go to $0$ at the end of the proof. We let $\phi_0,\cdots,\phi_k$ be a set of smooth functions, which we can assume to be orthogonal with respect to $\mu$ such that 
$$ \max_{\phi \in \langle \phi_0,\cdots,\phi_k \rangle} \frac{\int_{\widetilde{\Sigma}} \vert \nabla \phi \vert^2_g dA_g}{\int_{\widetilde{\Sigma}} \phi^2 d\mu} \leq \lambda_k(\widetilde{\Sigma},\tilde{g},\mu) + \delta $$
We will use these functions as test functions for the sequence $(\Sigma,g,\mu_n)$. For $1\leq j \leq l$, we let $\eta_j^n$ be smooth functions such that 
\begin{itemize}
\item $\eta_0^n \in \mathcal{C}^\infty_c(\Sigma \setminus \mathbb{D}_{\rho}(Z_0))$ with $0 \leq \eta_0^n \leq 1$, $\eta_0^n = 1$ in $\Sigma \setminus \mathbb{D}_{\sqrt{\rho}}(Z_0)$ and 
$$\int_{\Sigma} \vert \nabla \eta_0^n \vert^2_g dA_g \leq \frac{C}{\ln \frac{1}{\rho}}$$.
\item for $1\leq j \leq l$, $\eta_j^n \in \mathcal{C}^\infty_c(\mathbb{D}_{\frac{1}{\rho}} \setminus \mathbb{D}_{\delta}(Z_j))$  with $0 \leq \eta_j^n \leq 1$, $\eta_j^n = 1$ in $\mathbb{D}_{\frac{1}{\sqrt{\rho}}} \setminus \mathbb{D}_{\sqrt{\rho}}(Z_j)$ and 
$$\int_{\mathbb{R}^2} \vert \nabla \eta_j^n \vert^2 dxdy \leq \frac{C}{\ln \frac{1}{\rho}}.$$
\end{itemize}
We set for $0\leq i \leq k$ and for $x\in \Sigma$ (with the abuse of notation corresponding to local computation in conformal charts of an atlas)
\begin{equation} \label{eqdefphiintest} \phi_i^n(x) := \left(\eta_0^n \cdot (\phi_i)_{\vert \Sigma}\right)(x) + \sum_{j=1}^l \left(  \eta_j^n\cdot (\phi_i)_{\vert \mathbb{S}^2} \circ \pi_{\mathbb{S}^2}^{-1} \right)\left( \frac{x - x_j^\eps}{\alpha_j^\eps} \right) . \end{equation}

Now we aim at testing the functions $\phi_0^n,\cdots,\phi_k^n$ in the variational characterization of $\lambda_k(\Sigma,g,\mu_n)$. Let $\phi \in  \langle \phi_0^n,\cdots,\phi_k^n\rangle$ written as $\phi =  \sum_{j=0}^l a_j^n \phi_j^n$ and let $\psi = \sum_{j=0}^l a_j^n \phi_j$. We renormalize $(a_j^n)$ so that $\int_{\widetilde{\Sigma}} \psi^2 d\mu = 1$. Since every term in the sum of \eqref{eqdefphiintest} have disjoint support for $n$ large enough, and using the conformal invariance of the Dirichlet energy, we obtain
\begin{equation*}
\begin{split}\int_{\Sigma} \vert \nabla \phi \vert^2_g dA_g = \int_\Sigma \left\vert \nabla \left(\eta_0^n \cdot \psi_{\vert \Sigma}\right) \right\vert^2_g dA_g + \sum_{j=1}^l \int_{\mathbb{S}^2} \left\vert \nabla \left(\eta_j^n \circ \pi_{\mathbb{S}^2} \cdot \psi_{\vert (\mathbb{S}^2)_j}\right) \right\vert^2_{\mathbb{S}^2} dA_{\mathbb{S}^2} \\
\leq \int_{\widetilde{\Sigma}} \vert \nabla \psi \vert_{\widetilde{g}}^2 dA_{\widetilde{g}} + 2(l+1) \Vert \psi \Vert_{L^\infty} \sqrt{\frac{C}{\ln \frac{1}{\rho}}} \sqrt{\int_{\widetilde{\Sigma}} \vert \nabla \psi \vert_{\widetilde{g}}^2 dA_{\widetilde{g}}} + (l+1)\frac{C}{\ln \frac{1}{\rho}} \Vert \psi \Vert_{L^\infty}^2
\end{split} \end{equation*}
We also have that
\begin{equation*}
\begin{split}
\int_{\Sigma}  \phi^2 d\mu_n & = \int_\Sigma  \left(\eta_0^n \cdot \psi_{\vert \Sigma}\right)^2 d\mu_n + \sum_{j=1}^l \int_{\mathbb{S}^2} \left(\eta_j^n \circ \pi_{\mathbb{S}^2} \cdot \psi_{\vert (\mathbb{S}^2)_j}\right)^2  d\mu_j^n \\
& = \int_{\widetilde{\Sigma}}  \psi^2 d\mu + \int_{\widetilde{\Sigma}} \psi^2 \left( \sum_{j=0}^l \left(\eta_j^n\right)^2 d\mu_j^n - d\mu\right) 
\end{split} 
\end{equation*}
so that
$$ \left\vert \int_{\Sigma}  \phi^2 d\mu_n -  \int_{\Sigma}  \psi^2 d\mu  \right\vert \leq \Vert \psi \Vert_{L^{\infty}}^2 o_{n,\rho}(1)   $$
where $o_{n,\rho}(1)$ converges to $0$ when $n\to +\infty$ and then $\rho\to 0$. By equivalence of the $L^\infty$ norm and the Euclidean norm associated to $\left(\psi_1,\psi_2\right)\mapsto \int_{\widetilde{\Sigma}} \psi_1\psi_2 d\mu$, on the finite dimensional set $\langle \phi_0,\cdots,\phi_n \rangle$ and the assumption $\int_{\widetilde{\Sigma}} \psi^2 d\mu = 1$ we obtain that $\Vert \psi \Vert_{L^{\infty}}$ is bounded by a constant independent of $n$ and $\rho$. Then
$$ \frac{\int_{\Sigma} \vert \nabla \phi \vert^2_g dA_g}{\int_{\Sigma} \phi^2 d\mu_n} \leq \frac{\int_{\widetilde{\Sigma}} \vert \nabla \psi \vert_{\widetilde{g}}^2 dA_{\widetilde{g}} + O\left(\frac{1}{\sqrt{\ln \frac{1}{\rho}}}\right)}{\int_{\widetilde{\Sigma}} \psi^2 d\mu + o_{n,\rho}(1) } \leq \lambda_k(\widetilde{\Sigma},\widetilde{g},\mu) + \delta + o_{n,\rho}(1). $$
In addition, similarly to the previous computations, we have
$$ \left\vert \int_{\Sigma}  \phi_i^n \phi_j^n d\mu_n -  \int_{\Sigma}  \phi_i\phi_j d\mu  \right\vert \leq \Vert \phi_i \Vert_{L^{\infty}} \Vert \phi_j \Vert_{L^{\infty}} o_{n,\rho}(1)   $$
so that for $n$ large enough and $\rho$ small enough, the family $(\phi_0^n,\cdots,\phi_l^n)$ is independent on $supp(\mu_n)$ since $(\phi_0,\cdots,\phi_n)$ is orthonormal with respect $\mu$. The variational characterization of $\lambda_k(\Sigma,g,\mu_n)$ then yields
$$ \lambda_k(\Sigma,\mu_n,g) \leq \max_{\phi \in \langle \phi_0^n,\cdots,\phi_k^n\rangle } \frac{\int_{\Sigma} \vert \nabla \phi \vert^2_g dA_g}{\int_{\Sigma} \phi^2 d\mu_n} \leq \lambda_k(\widetilde{\Sigma},\widetilde{g},\mu) + \delta + o_{n,\rho}(1). $$
Letting $n\to +\infty$, then $\rho\to 0$ and then $\delta\to 0$, we obtain the expected result.
\end{proof}

We currently have the following general property. The notation for $I$ is kept in all the paper.
\begin{prop} \label{prop:bubbletree}
Let $(\beta_\eps ,\Phi_\eps,e^{2u_\eps}g)$ be a $(PS)_K$ sequence. We assume in addition that $\lambda_\eps^I \to 0$ and that $\left(\lambda_\eps^{I+1}\right)$ is uniformly lower bounded as $\eps\to 0$. Then up to the extraction of a subsequence, $\beta_\eps$ and $e^{2u_\eps}dA_g$ MW$\star$-bubble tree converge to the same measures $\mu_0$ on $\Sigma$ and $\mu_j$ in $\left(\mathbb{S}^2\right)_j$ for $1\leq j\leq l$ where $l\leq I+1$. 
\end{prop}

We denote for $1\leq j\leq l$, $(x_j^\eps,\alpha_j^\eps)$ the associated points and scales. We denote for $0\leq j \leq l$, $\mu_j$ the pullback of the continuous part of $\nu_j$ (having the set of atoms $Z_j$) with repect to $\pi_{\mathbb{S}^2}$. The functions $f_\eps : \Sigma \to \R$ we consider are seen at the scale $(x_j^\eps,\alpha_j^\eps)$ with the formula
$$ f_j^\eps := f_\eps( x_j^\eps + \alpha_j^\eps \pi_{\mathbb{S}^2}^{-1}\left(\cdot\right) ) $$
and in particular, we denote $\widetilde{\Phi_\eps}^j := \left(\Phi_\eps\right)_j$
while linear forms on continuous functions (measures) $\mu_\eps$ or bilinear forms on $H^1$ functions $\beta_\eps$ satisfy at the scale $(x_j^\eps,\alpha_j^\eps)$ for $\varphi,\psi \in \mathcal{C}_c^\infty\left( \mathbb{\pi}_{\mathbb{S}^2}^{-1}\left( \R^2 \setminus Z_j \right) \right)$
$$ \langle \mu_j^\eps , \varphi \rangle := \left\langle \mu_\eps , \varphi\left( \frac{ \pi_{\mathbb{S}^2}\left(\cdot\right) - x_j^\eps }{\alpha_j^\eps} \right) \right\rangle$$
$$ \beta_j^\eps(\varphi,\psi):= \beta_\eps \left(\varphi\left( \frac{\pi_{\mathbb{S}^2}\left(\cdot\right) - x_j^\eps}{\alpha_j^\eps} \right),\psi\left( \frac{ \pi_{\mathbb{S}^2}\left(\cdot\right) - x_j^\eps }{\alpha_j^\eps} \right)\right)$$
and in particular, we denote $e^{2u_j^\eps}dA_{\mathbb{S}^2} := \left( e^{2u_\eps} dA_g \right)_j $. 

We say that the analysis in $\Sigma$ if $\mu_0 \neq 0$ and in $\left(\mathbb{S}^2\right)_j$ for $1\leq j \leq l$ (in this case $\mu_j \neq 0$) at the scale $(x_j^\eps,\alpha_j^\eps)$ of functions, measures, bilinear forms and sets we consider is an analysis in a "thick part" since the measure has a positive mass at this scale.

We also localize the space $\bar{X}$. We denote for an open set $\Omega$ of a smooth Riemannian surface $(\Sigma,g)$, $\overline{X}\left(\Omega,g\right)$ the closure of 
$$ X(\Omega,g) := \{ \left(\varphi,\psi\right)\mapsto \int_\Sigma \varphi \psi e^{2u}dA_g ; u\in \mathcal{C}^\infty(\Omega) \} $$
in the set of symmetric bilinear forms on $H^{1}_0(\Omega)$ endowed with the norm
$$ \Vert \beta \Vert_{\overline{X}\left(\Omega,g\right)} := \sup_{\varphi,\psi \in H^1_0(\Omega) } \frac{\vert \beta(\varphi,\psi)\vert}{\Vert \varphi \Vert_{H^1_0(\Omega,g)} \Vert \psi \Vert_{H^1_0(\Omega,g)}}. $$

\begin{proof}[Proof of Proposition \ref{prop:bubbletree}]
It was proved for sequences of smooth metrics $g_\eps = e^{2u_\eps}g$ in \cite{petrides-2} and \cite{petrides-3}. Let us proove that $\beta_\eps$ MW$\star$ converges to the same limits as $e^{2u_\eps}dA_g$ in the same scales $\nu_0,\cdots,\nu_j$. In $\Sigma$, we have that for $\varphi \in \mathcal{C}_c^{\infty}(\Sigma)$,
$$ \left\vert \beta_\eps(1,\varphi) - \int_{\Sigma} \varphi d\nu_0 \right\vert \leq \left\vert \int_\Sigma \varphi(e^{2u_\eps}dA_g - d\nu_0) \right\vert + \Vert \varphi \Vert_{H^1(g_\eps)} \Vert e^{2u_\eps} - \beta_\eps \Vert_{g_\eps} = o(1) $$
as $\eps \to 0$ since by conformal invariance
$$ \Vert \varphi \Vert_{H^1(g_\eps)}^2 = \int_\Sigma \vert \nabla \varphi \vert^2_g dA_g + \int_{\Sigma} \varphi^2 dA_{g_\eps} \leq \Vert \nabla \varphi \Vert_{H^1(g)}^2 + \Vert  \varphi \Vert_{L^{\infty}}^2 $$
is bounded by a constant independent of $\eps$ and $\Vert e^{2u_\eps} - \beta_\eps \Vert_{g_\eps} = O(\eps)$. Since $(\beta_\eps)$ can be seen as a sequence of measures, its weak$\star$ limit has to be $\nu_0$ by uniqueness of the limit in the sense of distributions.

Similarly, we have at the scale $(x_j^\eps,\alpha_j^\eps)$ that for $\varphi \in \mathcal{C}^{\infty}_c(\pi_{\mathbb{S}^2}^{-1}\left(\mathbb{R}^2\right))$,
\begin{equation*}\begin{split} \left\vert \beta_j^\eps(1,\varphi) - \int_{\mathbb{S}^2} \varphi d\pi_{\mathbb{S}^2}^\star(\nu_j) \right\vert \leq & \left\vert \int_{\mathbb{S}^2} \varphi(e^{2u_j^\eps}dA_{\mathbb{S}^2} - d\pi_{\mathbb{S}^2}^\star(\nu_j)) \right\vert \\ & + \Vert \varphi \Vert_{H^1\left(\Omega, e^{2u_j^\eps }\right)} \Vert e^{2u_j^\eps}dA_{\mathbb{S}^2} - \beta_j^\eps \Vert_{\Omega, e^{2u_j^\eps}}
\end{split} \end{equation*}
where $\Omega$ is an open set that contains the support of $\varphi$ and $\Vert \varphi \Vert_{H^1\left(\Omega, e^{2u_j^\eps }\right)}$ is again uniformly bounded by the use of the conformal invariance of the Dirichlet energy since we used conformal charts and by definition,
$$ \Vert e^{2u_j^\eps}dA_{\mathbb{S}^2} - \beta_j^\eps \Vert_{\Omega, e^{2u_j^\eps}} = \Vert e^{2u_\eps}dA_g - \beta_\eps \Vert_{g_\eps} = O(\eps) $$
as $\eps\to 0$.
\end{proof}

The goal in all the section is to prove that the limiting measures $\mu_0,\cdots,\mu_l$ are absolutely continous with respect to $dA_g$ or  $dA_{\mathbb{S}^2}$ with densities satisfying the conclusions of Proposition \ref{palaissmalegeneral}.

\subsection{Some convergence of $\omega_{\eps}$ to $1$ in thick parts and first replacement of $\Phi_{\eps}$}
We set 
$$\omega_{\eps} = \sqrt{\left\vert \Phi_{\eps} \right\vert_{\Lambda_{\eps}}^2 + \theta_\eps^2 } 
$$ 
We first have that $\nabla \omega_{\eps}$ converges to $0$ in $L^2$ and that $\sqrt{\Lambda_\eps}\cdot\nabla\Phi_{\eps}$ has a similar $L^2$ behaviour as $\sqrt{\Lambda_\eps}\cdot\nabla \frac{\Phi_{\eps}}{\omega_{\eps}}$.
\begin{cl} We have that
\begin{equation}  \label{eqomegaepsto1}   \int_{\Sigma} \left\vert \nabla \omega_\eps \right\vert^2 + \int_\Sigma \left\vert \nabla\left( \Phi_\eps - \frac{\Phi_\eps}{\omega_\eps} \right) \right\vert_{\Lambda_\eps}^2  + \int_\Sigma(\omega_\eps^2-1)\left\vert \nabla \frac{\Phi_\eps}{\omega_\eps} \right\vert^2_{\Lambda_\eps} \leq O(\eps) \end{equation}
as $\eps\to 0$.
\end{cl}

\begin{proof} We first prove
\begin{equation} \label{eqomega_epsminus1} L_\eps\left( \left\vert \Lambda_\eps \Phi_\eps \right\vert^2 \left(1- \frac{1 }{\omega_\eps}\right) \right) \leq O(\eps) \end{equation}
as $\eps\to 0$. Since $\omega_\eps \geq 1$,  
and $\left\vert \Phi_\eps \right\vert_{\Lambda_\eps}^2 \leq \omega_\eps^2$, we have that
\begin{equation*}
\begin{split}
 L_\eps\left( \left\vert \Lambda_\eps \Phi_\eps \right\vert^2 \left(1- \frac{1 }{\omega_\eps} \right) \right) \leq & \lambda_K^\eps L_\eps\left(  \left(\omega_\eps^2-\omega_\eps\right) \right) \\
\leq & \lambda_K^\eps \left( L_\eps\left(  \left\vert \Phi_\eps \right\vert_{\Lambda_\eps}^2 \right) + L_\eps \left(  \theta_\eps^2 \right) - L_{\eps}(1) \right) \\
 \end{split} \end{equation*}
 so that
$$ L_\eps\left( \left\vert \Lambda_\eps \Phi_\eps \right\vert^2 \left(1- \frac{1 }{\omega_\eps}\right) \right) \leq  \lambda_K^\eps  L_\eps( \theta_\eps^2) \leq  \lambda_K^\eps  \left\Vert \beta_\eps \right\Vert_{g_\eps} \left\Vert \theta_\eps \right\Vert^2_{H^1(g_\eps)} \leq O(\eps) $$
as $\eps \to 0$ and we obtain \eqref{eqomega_epsminus1}.

We now prove \eqref{eqomegaepsto1}:
\begin{equation*} 
\begin{split}
\int_{\Sigma} & \left\vert \nabla \frac{\Phi_\eps}{\omega_\eps} \right\vert_{\Lambda_\eps}^2 - \int_\Sigma  \left\vert \nabla \Phi_\eps \right\vert_{\Lambda_\eps}^2 - \int_\Sigma  \left\vert \nabla\left( \Phi_\eps - \frac{\Phi_\eps}{\omega_\eps} \right) \right\vert_{\Lambda_\eps}^2 
\\ = & - 2 \int_\Sigma \left\langle \nabla \Phi_\eps, \nabla\left( \Phi_\eps - \frac{ \Phi_\eps}{\omega_\eps}\right) \right\rangle_{\Lambda_\eps} 
= - 2 \int_\Sigma \Delta \Phi_\eps  {\Lambda_\eps}.\left( \Phi_\eps - \frac{ \Phi_\eps}{\omega_\eps}\right)  \\
= & - 2 \beta_\eps\left( \Lambda_\eps.\Phi_\eps , {\Lambda_\eps}.\left( \Phi_\eps - \frac{ \Phi_\eps}{\omega_\eps}\right) \right)
= - 2 L_\eps\left(   \left\vert \Lambda_\eps \Phi_\eps \right\vert^2 \left(1- \frac{1 }{\omega_\eps} \right) \right) = O(\eps)
\end{split}
\end{equation*}
where we tested $ \Delta \Phi_{\eps} =  \beta_\eps( \Lambda_\eps\Phi_\eps,.) $ in $\Sigma$ against $ \Lambda_\eps. \left( \Phi_\eps-\frac{ \Phi_\eps}{\omega_\eps}\right)$, and we used \eqref{eqomega_epsminus1}.

In particular, we have
$$ 0 \leq \int_\Sigma \left\vert \nabla\left( \Phi_\eps - \frac{\Phi_\eps}{\omega_\eps} \right) \right\vert_{\Lambda_\eps}^2 \leq \int_\Sigma  \left( \left\vert \nabla \frac{\Phi_\eps}{\omega_\eps} \right\vert_{\Lambda_\eps}^2 - \left\vert \nabla \Phi_\eps \right\vert_{\Lambda_\eps}^2\right) + O(\eps) $$
as $\eps\to 0$ and knowing that with the straightforward computations we have
\begin{equation*}
\begin{split} \left\vert \nabla \frac{\Phi_\eps}{\omega_\eps} \right\vert_{\Lambda_\eps}^2 - \left\vert \nabla \Phi_\eps \right\vert_{\Lambda_\eps}^2 = & \left( 1 - \omega_\eps^2  \right)\left\vert \nabla \frac{\Phi_\eps}{\omega_\eps} \right\vert_{\Lambda_\eps}^2 - \left\vert \nabla \omega_\eps \right\vert^2 \frac{\omega_\eps^2+\theta_\eps^2}{\omega_\eps^2}  + 2 \frac{\theta_\eps}{\omega_\eps} \nabla \omega_\eps \nabla \theta_\eps \\
=  & \left( 1 - \omega_\eps^2  \right)\left\vert \nabla \frac{\Phi_\eps}{\omega_\eps} \right\vert_{\Lambda_\eps}^2 - \left\vert \nabla \omega_\eps \right\vert^2 - \left\vert \frac{\theta_\eps}{\omega_\eps} \nabla \omega_\eps - \nabla \theta_\eps \right\vert^2 + \left\vert \nabla \theta_\eps \right\vert^2
\end{split} \end{equation*}
where 
$$ \left\vert \frac{\theta_\eps}{\omega_\eps} \nabla \omega_\eps - \nabla \theta_\eps \right\vert^2 = \omega_\eps^2 \left\vert \nabla \frac{\theta_\eps}{\omega_\eps} \right\vert^2  $$
we obtain that 
\begin{equation*}
\begin{split} \int_{\Sigma} \left( \omega_\eps^2 - 1 \right)\left\vert \nabla \frac{\Phi_\eps}{\omega_\eps} \right\vert_{\Lambda_\eps}^2 + \int_{\Sigma} \left\vert \nabla \omega_\eps \right\vert^2 +  \int_{\Sigma} \omega_\eps^2 \left\vert \nabla \frac{\theta_\eps}{\omega_\eps} \right\vert^2 + \int_\Sigma \left\vert \nabla\left( \Phi_\eps - \frac{\Phi_\eps}{\omega_\eps} \right) \right\vert_{\Lambda_\eps}^2 \\
 \leq \int_\Sigma \left\vert \nabla \theta_\eps \right\vert^2 + O\left(\eps \right) 
\end{split} \end{equation*}
as $\eps \to 0$. 
\end{proof}

\subsection{Good/bad points in thick parts and immediate consequences}

\subsubsection{Construction of a finite number of bad points}

In the following, we perform local regularity estimates on $(\Phi_\eps)$. These estimates can only be done far from "bad points" we select in Claim \ref{clbadpoints}.
For $\Omega \subset \Sigma$ a domain of $\Sigma$, we recall that
$$ \lambda_\star(\Omega, g, \beta_\eps) = \inf_{\varphi \in \mathcal{C}_c^{\infty}\left(\Omega\right) } \frac{\int_\Omega  \left\vert \nabla \varphi \right\vert_g^2 dA_g }{\beta_\eps(\varphi,\varphi) }. $$
We recall that $\lambda_K^\eps := \max_{i\in\{1,\cdots,m_\eps\}}\lambda_{i}^\eps$ where $\lambda_i^\eps$ is a $i$-th eigenvalue on $(\Sigma,g,\beta_\eps)$.  We also let $g_j = g$ if $j =0$ and $g_j = g_{\mathbb{S}^2}$ if $j \geq 1$. We have:
\begin{cl} \label{clbadpoints} Up to a subsequence, there is $0< r_{\star} <1$ and a set of at most $K+1$ bad points $P_j \subset \Sigma$ and  such that for any $p \in \Sigma \setminus P_j$ and any $r < \min\left( r_\star, d_g(p,P_j ) \right) $, then for $\eps$ small enough,
$$ \lambda_\star\left( \mathbb{D}_r(p), g_j, \beta_\eps^j \right) \geq \lambda_K^\eps. $$
\end{cl}

\begin{proof}
We only prove the result for $j=0$ since it is similar for $j\geq 1$. Just notice that for $j\geq 1$, the north pole (used for the stereographic projection) is automatically considered as a bad point. We set 
$$ r_\eps^1 =  \inf\{ r>0 ; \exists p\in \Sigma, \lambda_\star\left( \mathbb{D}_r(p),g, \beta_\eps \right) < \lambda_K^\eps \}. $$
If $r_\eps^1$ does not converge to $0$, then up to a subsequence, there is $r_{\star\star}$ such that $r_\eps^1 > r_{\star\star}$ and Claim \ref{clbadpoints} is proved for this $r_{\star\star}$ and $P = \emptyset$.
If $r_\eps^1 \to 0$, then, we let $p_1^\eps$ be such that $ \lambda_\star\left( \mathbb{D}_{r_\eps^1}(p^1_\eps),g, \beta_\eps \right) < \lambda_\eps $  (up to take $r_\eps^1 + \eps$ instead of $r_\eps^1$ in order to have the strict inequality).
By induction assume that for $j \in \mathbb{N}$ we constructed $r_\eps^1 \leq r_\eps^2 \leq \cdots \leq r_\eps^{j-1}$ such that $r_\eps^{j-1}\to 0$ and points $p^1_\eps,\cdots,p^{j-1}_\eps$ such that 
$$ \forall i \neq l, \mathbb{D}_{r_\eps^i}(p_\eps^i) \cap \mathbb{D}_{r_\eps^l}(p^l_\eps) = \emptyset  \text{ and } \forall i, \lambda_\star\left( \mathbb{D}_{r_\eps^i}(p^i_\eps),g,\beta_\eps \right) <\lambda_K^\eps$$
then we let $r_\eps^j$ be the following infimum
$$ \inf\{ r>0 ; \exists p\in \Sigma,  \forall i, \mathbb{D}_r(p) \cap  \mathbb{D}_{r^i_\eps}(p^i_\eps) = \emptyset  \text{ and } \lambda_\star\left( \mathbb{D}_r(p), g,\beta_\eps \right) < \lambda_K^\eps \} $$
Then if $r_\eps^j $ does not converge to $0$ and up to a subsequence, there is $r_{\star\star}$ such that $r_\eps^j > r_{\star\star} $ and Claim  \ref{clbadpoints} is proved for this $r_{\star\star}$ and $P = \{p_1,\cdots, p_{j-1}\}$ where up to a subsequence we took $p_1,\cdots, p_{j-1}$ as limits of $p_1^\eps,\cdots, p_{j-1}^\eps$ as $\eps \to 0$. 

If $r_\eps^j \to 0$, then let $p_j^\eps$ be such that  $ \lambda_\star\left( \mathbb{D}_{r_\eps^j}(p^1_\eps), g,\beta_\eps \right) < \lambda_K^\eps$ and $\mathbb{D}_{r_\eps^j}(p_\eps^j) \cap \mathbb{D}_{r_\eps^j}(p_\eps^i)=\emptyset$ for $i<j$ (up to take $r_\eps^j + \eps$ again).

This induction process has to stop because if we have we constructed $r_\eps^1 \leq r_\eps^2 \leq \cdots \leq r_\eps^{k+1}$ such that $r_\eps^{k+1}\to 0$ and points $p^1_\eps,\cdots,p^{k+1}_\eps$ such that 
$$ \forall i \neq l, \mathbb{D}_{r^i_\eps}(p_\eps^i) \cap \mathbb{D}_{r^l_\eps}(p^l_\eps) = \emptyset  \text{ and } \forall i, \lambda_\star\left( \mathbb{D}_{r_\eps^i}(p^i_\eps), g,\beta_\eps \right) < \lambda_K^\eps$$
Let $\varphi_i$ be the first eigenfunction associated to $\lambda_\star\left( \mathbb{D}_{r_\eps^i}(p^i_\eps),g, \beta_\eps \right)$ extended by $0$ in $\Sigma \setminus \mathbb{D}_{r_\eps^i}(p^i_\eps)$. We have by the min-max characterization of the $K$-th eigenvalue on $M$, $\lambda_\eps$ and since $\varphi_i$ are orthogonal functions that
$$ \lambda_K^\eps \leq \max_{i=1,\cdots, K+1} \frac{\int_\Sigma \left\vert \nabla \varphi_i \right\vert^2_g  dA_g}{\beta_\eps\left(\varphi_i,\varphi_i\right)} < \lambda_K^\eps$$
which is a contradiction.
\end{proof}

In the following, for $\rho >0$, we denote 
$$ \Omega_\rho^0 = \Sigma \setminus \bigcup_{p\in P_0} \mathbb{D}_{\rho}(p) \text{ and } \Omega_\rho^j = \mathbb{S}^2 \setminus \bigcup_{p\in P_j} \mathbb{D}_{\rho}(p). $$

\subsubsection{Smallness of $\omega_\eps - 1$ and $\theta_\eps$ near good points of thick parts}

We have the following convergence of $\omega_\eps$ to $1$ and $\theta_\eps$ to $0$ in thick parts. It also gives that if $\lambda_i^\eps \to 0$, then $\int_{\Omega_\rho^j} \left(\sqrt{\lambda_i^\eps}  \widetilde{\phi_i^\eps}^j  \right)^2 \to 0$ as $\eps \to 0$. 
\begin{cl} \label{clsmallthetaomega-1} We have for any $0 < \rho \leq \rho_0$ that for $1\leq j\leq l$ and for $j = 0$ if $\mu_0\neq 0$ that
\begin{equation}  \label{eqomegaepsto1general2}   \int_{\Omega_\rho^j}(\omega_\eps^j -1)^2 +\int_{\Omega_\rho^j}\left(\theta_\eps^j\right)^2 \leq O\left( \eps\right) \end{equation}
\begin{equation} \label{eq:poincareeigenfunctions}
\int_{\Omega_\rho^j} \left(\sqrt{\lambda_i^\eps}  \widetilde{\phi_i^\eps}^j  \right)^2  \leq O\left( \lambda_i^\eps t_i^\eps \right) 
\end{equation}
as $\eps\to 0$ and then $\rho \to 0$ (where the integrals are computed with respect to $dA_g$ if $j=0$ and the measure $dA_{\mathbb{S}^2}$ if $j\geq 1$).
\end{cl}

\begin{proof}
We will use Poincaré inequalities. Let $\eta \in \mathcal{C}^{\infty}_c(\Omega_{\rho}^0)$ be such that $0 \leq \eta \leq 1$ and $\eta = 1$ in $\Omega_{2\rho}^0$. In particular, for $\rho$ small enough, since $\mu_0 \neq 0$, $\beta_\eps(1,\eta)$ is uniformly lower bounded. Since $ \frac{\beta_\eps(.,\eta)}{\beta_\eps(1,\eta)}$ is a projection on $H^1 \to H^1$ that
$$ \int_{\Sigma} \left( \omega_\eps - \frac{\beta_\eps(\omega_\eps,\eta)}{\beta_\eps(1,\eta)} \right)^2 dA_g \leq C \left\Vert \frac{\beta_\eps(.,\eta)}{\beta_\eps(1,\eta)} \right\Vert_{H^{-1}(g)}^2 \int_{\Sigma} \left\vert \nabla \omega_\eps \right\vert^2_g dA_g $$
Since $ \frac{\beta_\eps(.,\eta)}{\beta_\eps(1,\eta)}$ is bounded in $H^{-1}$, the left-hand term is bounded by $O(\eps)$ as $\eps\to 0$. Now, similarly to the proof of \eqref{eqomega_epsminus1}, we have that
$$ \left\vert \frac{\beta_\eps(\omega_\eps,\eta)}{\beta_\eps(1,\eta)} - 1 \right\vert \leq  \frac{1}{\beta_\eps(1,\eta)} L_\eps\left(\left(\omega_\eps-1\right)^2\right)^{\frac{1}{2}} L_\eps\left(\eta^2\right)^{\frac{1}{2}}$$
where using that $\omega_\eps \geq 1$, $\omega_\eps^2 = \vert \Phi_\eps\vert_{\Lambda_\eps}^2 + \theta_\eps^2 $ and that $L_\eps(1) = L_\eps(\vert \Phi_\eps\vert_{\Lambda_\eps}^2)$,
$$  L_\eps\left(\left(\omega_\eps-1\right)^2\right) = L_\eps(\omega_\eps^2 - 2\omega_\eps +1) \leq  L_\eps(\omega_\eps^2 - 1) \leq L_\eps(\theta_\eps^2) \leq O(\eps)$$
as $\eps \to 0$ and that $\beta_\eps(1,\eta)$ is uniformly lower bounded so that
$$ \int_{\Sigma} \left( \omega_\eps - 1 \right)^2 dA_g = O(\eps) $$
as $\eps\to 0$ and doing the same with $\theta_\eps$ and $\phi_i^\eps$ completes the proof of estimate \eqref{eqomegaepsto1general2} and \eqref{eq:poincareeigenfunctions} for $j=0$ if $\mu_0 \neq 0$. Notice that the proof is analogous for $j\geq 1$.
\end{proof}

\subsubsection{Good annuli close to bad points}
We denote for a point $p$ and $r_2 < r_1$. 
$$ \mathbb{A}_{r_1,r_2}(p) := \mathbb{D}_{r_1}(p) \setminus \mathbb{D}_{r_2}(p) $$
\begin{cl} \label{clannulusbadpoint}
Let $j \in \{0,\cdots,l\}$ and let $p \in P_j$, then, up to the extraction of a subsequence there is $r>0$ and $s_\eps \to 0$ as $\eps \to 0$ such that
$$ \lambda_\star\left( \mathbb{A}_{r,s_\eps}(p),g_j, \beta_\eps^j \right) \geq \lambda_K^\eps $$
\end{cl}

\begin{proof}
We assume that the claim does not hold. 

\medskip

\noindent \textbf{Step 1:} Up to the extraction of a subsequence as $\eps \to 0$, we build by induction points $r_{K+1}< r_K < \cdots < r_1 < r_0$ such that for $i \in \{0,\cdots,K\}$
$$  \lambda_\star\left( \mathbb{A}_{r_i,r_{i+1}}(p),g_j, \beta_\eps^j \right) < \lambda_K^\eps $$

\medskip

\noindent \textbf{Proof of Step 1:} Let $r_0 >0$. Then, for $\eps$ small enough the set 
$$ \{ 0 < s < r ; \lambda_\star\left( \mathbb{A}_{r_0,s}(p),g_j, \beta_\eps^j \right) < \lambda_K^\eps \} $$
is not empty because if not, the Claim holds. Therefore, we can set
$$ s_\eps := \sup \{ 0 < s < r ; \lambda_\star\left( \mathbb{A}_{r_0,s}(p),g_j, \beta_\eps^j \right) < \lambda_K^\eps \} $$
We have that $s_\eps$ is lower bounded by a constant $c_0 > 0$ as $\eps \to 0$  (because if not, there is a subsequence such that $s_\eps \to 0$ take $r = r_0$ and $s_\eps + \eps$ instead of $s_\eps$ the claim holds). We set $ r_1 = \frac{c_0}{2} $. We now assume that $r_0 > r_1 > \cdots > r_k$ are built for some $k$ and we build $r_{k+1}$. As before, we set
$$ s_\eps := \sup \{ 0 < s < r ; \lambda_\star\left( \mathbb{A}_{r_k,s}(p),g_j, \beta_\eps^j \right) < \lambda_K^\eps \} $$
which satisfies $0 < s_\eps \leq r_\eps$. $s_\eps$ is lower bounded by a constant $c_k$ as $\eps\to 0$ because if not, the claim holds. We set $r_{k+1} = \frac{c_k}{2}$.
The proof of Step 1 is complete.

\medskip

\noindent \textbf{Step 2:} We obtain a contradiction: for $0\leq i\leq K$ we let $\varphi_i$ be the first eigenfunction associated to $\lambda_\star\left( \mathbb{A}_{r_i,r_{i+1}}(p),g_j, \beta_\eps^j \right) $ extended by $0$ outside $\mathbb{A}_{r_i,r_{i+1}}(p)$ and we test $\left\langle \varphi_i \right\rangle_{0 \leq i\leq K}$ (if $j= 0$) or $\left\langle \varphi_i\left( \frac{\cdot - x_j^\eps}{\alpha_j^\eps} \right) \right\rangle_{0 \leq i\leq K}$ (if $j\geq 1$) that belongs to $ \mathcal{G}_{K+1}(H^1(\Sigma))$  in the variational characterization of $\lambda_K^\eps$. Since $ \varphi_i  $ are orthogonal, we obtain that 
$$\lambda_K^\eps \leq \max_{i\in \{0,\cdots K\}} \lambda_\star\left( \mathbb{A}_{r_i,r_{i+1}}(p),g_j, \beta_\eps^j \right)  <  \lambda_K^\eps $$
 and this is a contradiction.
\end{proof}

\subsubsection{Non concentration of energies near good points and arbitrarily close to bad points}

\begin{cl} \label{propnonconcentgeneral}
Let $p\in \Sigma \setminus P_0 $ or $\mathbb{S}^2 \setminus P_j$, be a good point
then for any $r$ such that $\sqrt{r} < r_\star(p):= \min\left(r_\star,\frac{d(p,P_j)}{2}\right)$ and any function $\zeta \in \mathcal{C}_c^\infty(\mathbb{D}_r(p))$ such that $0 \leq \zeta \leq 1$
\begin{equation} \label{noconcent}
\lim_{r\to 0} \limsup_{\eps\to 0} L_j^\eps( \zeta ) = \lim_{r\to 0} \limsup_{\eps\to 0}\int_{\mathbb{D}_r(p)} \left\vert \nabla \widetilde{\Phi_{\eps}}^j \right\vert_{\Lambda_\eps}^2 =  0
\end{equation}
In addition, we have that for a bad point $p \in P_j$ and $r \leq r_\star$, and any function $\zeta \in \mathcal{C}_c^\infty(\mathbb{A}_{r,\sqrt{s_\eps}}(p))$ such that $0 \leq \zeta \leq 1$
\begin{equation} \label{noconcentannulus}
\lim_{r\to 0} \limsup_{\eps\to 0} L_j^\eps( \zeta ) = \lim_{r\to 0} \limsup_{\eps\to 0}\int_{\mathbb{A}_{r,\sqrt{s_\eps}}(p)} \left\vert \nabla \widetilde{\Phi_{\eps}}^j \right\vert_{\Lambda_\eps}^2 =  0
\end{equation}
\end{cl}

\begin{proof}
Let $\eta \in \mathcal{C}_c^{\infty}(\mathbb{D}_{\sqrt{r}}(p))$ with $0\leq \eta \leq 1$, $\eta = 1$ in $\mathbb{D}_{r}(p)$ and $\int_{\Sigma}\left\vert \nabla \eta\right\vert^2_{g} \leq \frac{C}{\ln\left(\frac{1}{r}\right)}$
$$ L_j^\eps(\zeta) \leq L_j^\eps(\eta^2) \leq \frac{1}{\lambda_{\star}(\mathbb{D}_{r_\star(x)}(p),g_j, \beta_j^\eps)} \int_{\Sigma}\left\vert \nabla \eta \right\vert^2 \leq  \frac{C}{\lambda_K^\eps\ln\left(\frac{1}{r}\right)} $$
Letting $\eps \to 0$ and then $r\to 0$ we obtain the first non-concentration property.
Now, we drop the index/exponent $j$ on the function $\widetilde{\Phi_\eps}^j$ and on $\beta_j^\eps$, $L_j^\eps$. We test $\Delta_g \Phi_\eps = \beta_\eps\left(\Lambda_\eps \Phi_\eps,.\right)$ against $\eta \frac{\Lambda_\eps\Phi_\eps}{\omega_\eps}$ and we obtain
$$ \int_{\Sigma} \eta \nabla \Phi_\eps \nabla  \frac{\Lambda_\eps \Phi_\eps}{\omega_\eps} = - \int_\Sigma \frac{\Lambda_\eps \Phi_\eps}{\omega_\eps} \nabla \Phi_\eps \nabla \eta + \beta_\eps\left( \frac{\vert \Lambda_\eps \Phi_\eps \vert^2}{\omega_\eps} , \eta\right) $$
so that
\begin{equation*} 
\begin{split}
\int_{\mathbb{D}_r(p)} \left\vert \nabla \Phi_\eps \right\vert_{\Lambda_\eps}^2 \leq &  \int_{\Sigma} \eta \left\vert \nabla \Phi_\eps \right\vert_{\Lambda_\eps}^2 \leq \left\vert \int_{\Sigma}\eta \left\langle \nabla  \Phi_\eps\nabla\left(\Phi_\eps-\frac{\Phi_\eps}{\omega_\eps}\right) \right\rangle_{\Lambda_\eps} \right\vert + \left\vert\int_{\Sigma} \eta \nabla \Phi_\eps \nabla \frac{\Lambda_\eps \Phi_\eps}{\omega_\eps} \right\vert \\
\leq & \left(\int_\Sigma \left\vert \nabla \Phi_\eps \right\vert^2_{\Lambda_\eps} \int_\Sigma \left\vert \nabla \left( \Phi_\eps - \frac{\Phi_\eps}{\omega_\eps} \right) \right\vert^2_{\Lambda_\eps}  \right)^{\frac{1}{2}} \\
& + \left(\int_{\Sigma}\left\vert \nabla \eta\right\vert^2_{g}\right)^{\frac{1}{2}}  \left(\int_\Sigma \left\vert \nabla \Phi_\eps \right\vert^2\right)^{\frac{1}{2}} \left\Vert \frac{\vert \Lambda_\eps \Phi_\eps \vert}{\omega_\eps} \right\Vert_\infty 
  + L_\eps(\eta^2)^{\frac{1}{2}}L_\eps\left(\frac{\vert \Lambda_\eps \Phi_\eps \vert^4}{\omega_\eps^2}\right)^{\frac{1}{2}} \\
\leq & O(\eps^{\frac{1}{2}}) + \frac{C^{\frac{1}{2}}}{\ln\left(\frac{1}{r}\right)^{\frac{1}{2}}} \left( \left\Vert \frac{\vert \Lambda_\eps \Phi_\eps \vert}{\omega_\eps} \right\Vert_\infty + \frac{1}{\sqrt{\lambda_K^\eps}} \left(\lambda_K^\eps\right)^{\frac{3}{2}} \right) \leq O(\eps^{\frac{1}{2}}) + \frac{2 \lambda_K^\eps C^{\frac{1}{2}}}{\ln\left(\frac{1}{r}\right)^{\frac{1}{2}}}
\end{split}
\end{equation*}
so that letting $\eps\to 0$ and then $r\to 0$, we obtain the second expected non-concentration property in \eqref{noconcent}.

The proof of \eqref{noconcentannulus} is similar with the choice of $\eta_\eps \in \mathcal{C}_c^{\infty}(\mathbb{A}_{\sqrt{r},s_\eps}(p))$ with $0\leq \eta_\eps \leq 1$, $\eta_\eps = 1$ in $\mathbb{A}_{r,\sqrt{s_\eps}}(p)$ and $\int_{\Sigma}\left\vert \nabla \eta_\eps\right\vert^2_{g} \leq \frac{C}{\ln\left(\frac{1}{r}\right)}$ and the use of Claim \ref{clannulusbadpoint}.
\end{proof}

\subsection{Construction of local harmonic replacements}
We set 
\begin{equation}\label{def:varphitauhatlambda}\hat{\theta}_\eps^j := (\theta_\eps^j , \sqrt{\lambda_1^\eps} \widetilde{\phi_1^\eps}^{j} ,\cdots,\sqrt{\lambda_I^\eps}\widetilde{\phi_I^\eps}^{j}) \text{ and } \varphi_\eps^j := \left(\widetilde{\phi_1^\eps}^j,\cdots,\widetilde{\phi_{m_\eps}^\eps}^{j}\right) \text{ and } \hat{\Lambda}_\eps := ( \lambda_{I+1}^\eps,\cdots,\lambda_{m_\eps}^\eps ). \end{equation} First we build a local replacement of $\widetilde{\Phi_\eps}^j$ which will be written $ \sqrt{\omega_\eps^2 - \vert\tau_\eps\vert^2} \Psi_\eps $ where $\tau_\eps$ is a local harmonic replacement into $\R^{I+1}$ of $\hat{\theta}_\eps^j$ and $\Psi_\eps$ is a local harmonic replacement into an Euclidean ellipsoid of parameter $\hat{\Lambda}_\eps$ of $\frac{\varphi_\eps^j}{\left\vert \varphi_\eps^j \right\vert_{\hat{\Lambda}_\eps}}$. In particular, in the following claim, we give a sense to the replacement $\Psi_\eps$ and prove that it can have an arbitrary small energy. We choose $\eps_0 := \eps_{\alpha}'$  in order to have \ref{eq:energyconvexity} with $\alpha$ an upper bound for $\max\left\{ \lambda_K^\eps, \left(\lambda_{I+1}^\eps\right)^{-1} \right\}$. This implies the uniqueness of the harmonic replacement.

\begin{cl} \label{clharmreplacement} There is $\eta>0$ such that for any $p\in \mathbb{S}^2$ (or $\Sigma$ if $j=0$ and $\mu_0\neq 0$) there is $r(p)>0$ such and $r(p)^2 \leq  r_\eps(p) \leq r(p)$ such that there are unique maps $\tau_\eps$ and $\Psi_\eps$ satisfying
$$ \tau_\eps = \hat{\theta}_\eps^j \text{ and } \left\vert \varphi_\eps^j \right\vert_{\hat{\Lambda}_\eps} \geq \frac{1}{2} \text{ and } \Psi_\eps = \frac{\varphi_\eps^j}{\left\vert \varphi_\eps^j \right\vert_{\hat{\Lambda}_\eps}} $$
almost everywhere on $ \partial \mathbb{D}_{r_\eps}(p)$, $\vert\Psi_\eps\vert_{\hat{\Lambda}_\eps} = 1$ and
$$ \int_{\mathbb{D}_{r_\eps(p)}(p)} \left\vert \nabla \Psi_\eps \right\vert^2 = \inf \left\{ \int_{\mathbb{D}_{r_\eps(p)}(p)} \left\vert \nabla \Psi \right\vert^2 ; \Psi \in H^1 \begin{cases} \vert \Psi \vert^2_{\hat{\Lambda}_\eps} =_{a.e} 1 \text{ in } \mathbb{D}_{r_\eps(p)}(p) \\ \Psi =_{a.e} \frac{\varphi_\eps^j}{\left\vert \varphi_\eps^j \right\vert_{\hat{\Lambda}_\eps}} \text{ on } \partial \mathbb{D}_{r_\eps(p)}(p) \end{cases} \right\} \leq \eps_\alpha' $$
and in particular $\Psi_\eps$ is a harmonic map into the ellipsoid $\{ \left\vert x \right\vert_{\hat{\Lambda}_\eps} = 1 \}$ and satisfies
$$ \Delta \Psi_\eps = \frac{\left\vert \nabla \Psi_\eps \right\vert^2_{\hat{\Lambda}_\eps}}{\left\vert \hat{\Lambda}_\eps \Psi_\eps \right\vert^2} \hat{\Lambda}_\eps \Psi_\eps $$
$$ \Delta_g \tau_\eps = 0 $$
and $\vert \tau_\eps \vert^2 \leq \frac{1}{4}$.
\end{cl}

\begin{proof}
During all the proof, we drop the indices or exponents $j$ of all the functions because the argument is similar in all the thick parts. Thanks to \eqref{noconcent}, let $p \in \Sigma \setminus P_0$ or $\mathbb{S}^2 \setminus P_j$, let $r_0(p) \leq r_{\star}$ be such that any small $\eps$,
$$ \int_{\mathbb{D}_{r_0}(p)} \left\vert \nabla \varphi_\eps \right\vert^2 \leq \delta \eps_0. $$
for a constant $0<\delta\leq 1$ we will choose later. If $p \in P_j$, with the use of \eqref{noconcentannulus}, we choose $r_0(p)$ such that , 
$$ \int_{\mathbb{A}_{r_0(p),\frac{r_0(p)^2}{4}}(p)} \left\vert \nabla \varphi_\eps \right\vert^2 \leq \delta \eps_0. $$
Let $\frac{r_0(p)}{2} < r < r_0(p)$. By the Courant-Lebesgue lemma, let $r^2 \leq r_{\eps} \leq r$ be a radius such that  
\begin{equation} \label{eqcourantlebesgue}
\begin{split} 
\int_{\partial\mathbb{D}_{r_\eps(p)}(p)} & \left\vert \partial_\theta \hat{\theta}_\eps \right\vert ^2 d\theta +
 \int_{\partial\mathbb{D}_{r_\eps(p)}(p)} \left\vert \partial_\theta \varphi_\eps \right\vert ^2 d\theta \\ \leq & 
\frac{1}{\ln 2} 
\left( \int_{\mathbb{A}_{r,r^2}(p)}  \left\vert \nabla \hat{\theta}_\eps \right\vert^2 + 
\int_{\mathbb{A}_{r,r^2}(p)}  \left\vert \nabla \varphi_\eps \right\vert^2 
\right) 
 \leq  \frac{2}{\ln 2} \delta \eps_0.
\end{split}  \end{equation}
A vector-valued Morrey embedding theorem yields
\begin{equation} \label{eqconsequencecourantlebesgue} 
\max_{q,q' \in \partial\mathbb{D}_{r_\eps(p)}(p)} \left\vert \tau_\eps(q) - \tau_\eps(q') \right\vert^2 + 
\max_{q,q' \in \partial\mathbb{D}_{r_\eps(p)}(p)} \sum_{i=1}^{n_\eps}  \left\vert \varphi_i^\eps(q) - \varphi_i^\eps(q') \right\vert^2 \leq \frac{2\pi}{\ln 2}  \delta \eps_0. \end{equation}
By the classical trace $L^2$ embedding into $H^1$ and the estimates \eqref{eqomegaepsto1general2} and \eqref{eq:poincareeigenfunctions}, we have that
$$ \int_{\partial \mathbb{D}_{r_\eps(p)}(p)} \left( \omega_\eps - 1  \right)^2 + \int_{\partial \mathbb{D}_{r_\eps(p)}(p)} \vert \hat{\theta}_\eps\vert ^2  \leq o(1) $$
as $\eps\to 0$. Knowing that $\left\vert \varphi_\eps \right\vert_{\hat{\Lambda}_\eps}^2 - 1 = \omega_\eps^2 -1 - \vert \hat{\theta}_\eps \vert^2$, we obtain that
$$ \left\vert \int_{\partial\mathbb{D}_{r_\eps(p)}(p)} ( \left\vert \varphi_\eps \right\vert_{\hat{\Lambda}_\eps} - 1) \right\vert \leq \left\vert \int_{\partial\mathbb{D}_{r_\eps(p)}(p)} ( \left\vert \varphi_\eps \right\vert_{\hat{\Lambda}_\eps}^2 - 1) \right\vert =  o(1) $$
as $\eps \to 0$ 
and since with \eqref{eqconsequencecourantlebesgue} we have
$$ -  \sqrt{ \frac{2\pi}{\ln 2}  \delta \eps_0 \lambda_K^\eps} + \left\vert \varphi_\eps \right\vert_{\hat{\Lambda}_\eps}(q') \leq \left\vert \varphi_\eps \right\vert_{\hat{\Lambda}_\eps}(q) \leq \left\vert \varphi_\eps \right\vert_{\hat{\Lambda}_\eps}(q') + \sqrt{ \frac{2\pi}{\ln 2} \delta \eps_0 \lambda_K^\eps }, $$
taking the mean value on $\partial\mathbb{D}_{r_\eps}(p)$ with respect to $q'$ gives
$$  \left\vert \left\vert \varphi_\eps \right\vert_{\hat{\Lambda}_\eps}(q) - 1 \right\vert \leq  o(1) +  \sqrt{\frac{2\pi}{\ln 2} \delta \eps_0 \lambda_K^\eps}$$
We choose $\eps$ small enough and $\eta \leq \frac{1}{64} \left(\frac{\pi}{\ln 2} \eps_0 \lambda_K^\eps\right)^{-1} $ we obtain that $\left\vert \varphi_\eps \right\vert_{\hat{\Lambda}_\eps}(q) \geq \frac{3}{4}$ and $\vert \hat{\theta}_\eps \vert^2 \leq \frac{1}{4}$ for any $q\in \partial\mathbb{D}_{r_\eps(p)}(p)$ and $\eps$ small enough. By the maximum principle $\vert \tau_\eps \vert^2 \leq \frac{1}{4}$ in $\mathbb{D}_{r_\eps(p)}(p)$

We let $\Psi_\eps : \mathbb{D}_{r_\eps(p)}(p) \to \mathcal{E}_{\hat{\Lambda}_\eps}$ be a harmonic extension of $\frac{\varphi_\eps}{\left\vert \varphi_\eps \right\vert_{\hat{\Lambda}_\eps}}$ (that is a minimizer of the energy on maps $\Psi$ satisfying $\left\vert \Psi \right\vert_{\hat{\Lambda}_\eps} = 1$).
In order to prove uniqueness of $\Psi_\eps$, we have to prove that its energy is small enough. 

Let $\eta \in \mathcal{C}^\infty_c\left( \mathbb{D}_{r^2}(p) \right)$ be a cut-off function such that $\eta \geq 1$ in $\mathbb{D}_{\frac{r^2}{2}}(p)$ and $\vert \nabla \eta \vert \leq \frac{1}{r}$. We set $T_\eps(x) := \left( 1 - \eta \right) \varphi_\eps\left( r_\eps \frac{x}{\vert x \vert} \right) + \eta \varphi_\eps(q_\eps)$ and we compute the energy of $\frac{T_\eps}{\vert T_\eps \vert}$ knowing that
$$ \int_{\mathbb{D}_{r_\eps(p)}(p)} \left\vert \nabla \Psi_\eps \right\vert^2_g dA_g \leq \int_{\mathbb{D}_{r_\eps(p)}(p)} \left\vert \nabla \frac{T_\eps}{\vert T_\eps \vert} \right\vert^2_g dA_g $$
We have that 
$$\left\vert \nabla \frac{T_\eps}{\vert T_\eps \vert} \right\vert^2 \leq  \frac{\vert \nabla T_\eps \vert^2}{\vert T_\eps \vert^2} \leq \frac{2\left( 1-\eta \right)^2 \frac{\vert \nabla_\tau \varphi_\eps \vert^2}{r^2} + 2\vert \nabla \eta \vert^2 \max_{q \in \partial\mathbb{D}_{r_\eps}(p)}\vert \varphi_\eps(q)-\varphi_\eps(q_\eps)  \vert^2}{\left( \vert \varphi_\eps(q_\eps) \vert - \max_{q \in \partial\mathbb{D}_{r_\eps}(p)}\vert \varphi_\eps(q)-\varphi_\eps(q_\eps)  \vert \right)^2} $$
so that using the previous smallness estimates coming from the Courant-Lebesgue property \eqref{eqconsequencecourantlebesgue},
and up to reduce $\delta$, we complete the proof of the Claim.
\end{proof}

\subsection{Local $H^1$ comparison of eigenfunctions to the harmonic replacements}

\begin{cl} \label{clcomparegeneral} We have for all $p\in \Sigma$ and $r_\eps(p)$ given by Claim \ref{clharmreplacement}
$$  \int_{\mathbb{D}_{r_\eps(p)}(p)}   \left\vert \nabla \left( \Psi_\eps- \hat{\varphi}_\eps^j \right) \right\vert^2 
= o\left(1\right) $$
as $\eps \to 0$ where with the notations of Claim \ref{clharmreplacement}
\begin{equation*}
 \hat{\varphi}_\eps^j = \begin{cases} \frac{\varphi_\eps^j}{\rho_\eps^j} \text{ if } p \in \Sigma \setminus P_0 \text{ if } j=0 \text{ or } \mathbb{S}^2 \setminus P_j \text{ if } j \geq 1 \\ \left(1-\eta_\eps\right)  \frac{\varphi_\eps^j}{\rho_\eps^j} +  \eta_\eps \Psi_\eps \text{ if } p \in P_j
 \end{cases} \end{equation*}
where $\rho_\eps^j := \sqrt{\left(\omega_\eps^j\right)^2 - \vert \tau_\eps \vert^2}$ and $\eta_\eps \in \mathcal{C}_c^\infty\left(\mathbb{D}_{\sqrt{s_\eps}}(p)\right)$ such that $\eta_\eps = 1$ in $\mathbb{D}_{s_\eps}(p)$, $0\leq \eta_\eps \leq 1$ satisfy as $\eps \to 0$
\begin{equation}\label{eq:smallenergiesetarho} \int_{\mathbb{D}_{r_\eps(p)}(p)} \vert \nabla \eta_\eps \vert^2 = O\left(\frac{1}{\ln \frac{1}{s_\eps}}\right) \text{ and } \int_{\mathbb{D}_{r_\eps(p)}(p)} \vert \nabla \rho_\eps^j \vert^2 = O(\eps)\end{equation}
\end{cl}

\begin{proof}
We only write the proof of the claim for $p \in P_j$ since the other case exactly follows the same proof with $\eta_\eps = 0$ and $\mathbb{D}_{r_\eps(p)}(p)$ instead of $\mathbb{A}_{r_\eps(p),s_\eps}(p)$. We drop the index/exponent $j$ in all the proof since it works the same way in every thick part. We let $r_\eps(p)$, $\Psi_\eps$, $\tau_\eps$ be given by Claim \ref{clharmreplacement}. 

Notice that \eqref{eq:smallenergiesetarho} on $\rho_\eps^j$ is a simple consequence of Claim \eqref{eqomegaepsto1} and the $(PS)_K$ that gives $\int_\Sigma \vert \nabla \tau_\eps \vert^2 = O(\eps)$. Notice that $\rho_\eps$ is chosen so that  $\hat{\varphi}_{\eps}^i - \Psi_{\eps}^i$ is equal to $0$ on $\partial\mathbb{D}_{r_\eps(p)}(p)$.   With the choice of $\eta_\eps$, it is equal to $0$ on $\partial\mathbb{A}_{r_\eps(p),s_\eps}(p)$. We will use this property in Step 1 and Step 2. Using both steps will complete the proof of the Claim.

\medskip
\noindent\textbf{Step 1:}
\medskip

\begin{equation} \label{eqsmalldiffenergy}
 \int_{\mathbb{D}_{r_\eps}(p)}  \left\vert \nabla \hat{\varphi}_\eps \right\vert^2 - \int_{\mathbb{D}_{r_\eps}(p)} \left\vert \nabla  \Psi_\eps \right\vert^2  \leq  o(1)
\end{equation}
as $\eps \to 0$

\medskip
\noindent \textbf{Proof of Step 1:}
\medskip
We test the function $ \hat{\varphi}_{\eps}^i - \Psi_{\eps}^i$ in the variational characterization of $\lambda_{\star}:= \lambda_{\star}\left(\mathbb{A}_{r_\eps(p),s_\eps}(p),\beta_\eps \right)$ knowing Claim \ref{clbadpoints}:
$$ \lambda_{i}^\eps L_\eps \left( \left(\hat{\varphi}_{\eps}^i - \Psi_{\eps}^i \right)^2 \right) \leq  \lambda_{\star} L_\eps \left( \left(\hat{\varphi}_{\eps}^i - \Psi_{\eps}^i \right)^2 \right) \leq \int_{ \mathbb{D}_{r_\eps(p)}(p) } \left\vert \nabla\left(\hat{\varphi}_{\eps}^i - \Psi_{\eps}^i\right)\right\vert^2$$
and we sum on $i$ to get
\begin{equation} \label{eqtestlambdastar}
 L_\eps\left( \left\vert \hat{\varphi}_{\eps} - \Psi_{\eps}\right\vert_{\hat{\Lambda}_\eps}^2\right) \leq   \int_{ \mathbb{D}_{r_\eps(p)}(p) } \left\vert \nabla \hat{\varphi}_{\eps}  \right\vert^2 +  \int_{ \mathbb{D}_{r_\eps(p)}(p) } \left\vert \nabla  \Psi_{\eps}  \right\vert^2 
- 2\int_{ \mathbb{D}_{r_\eps(p)}(p) } \nabla \hat{\varphi}_{\eps} \nabla \Psi_{\eps} 
 \end{equation}
 Now, we test the equation on $\Phi_\eps$: $\Delta_g \Phi_\eps = \beta_\eps(\Lambda_\eps\Phi_\eps,.) $ against $\frac{1-\eta_\eps}{\rho_\eps}\left(\hat{\varphi}_\eps - \Psi_\eps\right)$ and we multiply by $2$:
\begin{equation*}
\begin{split}2 \int_{\mathbb{D}_{r_\eps(p)}(p)} \nabla \varphi_\eps \nabla\left( \frac{1-\eta_\eps}{\rho_\eps}\left(\hat{\varphi}_\eps - \Psi_\eps\right)  \right) & = 2 L_{\eps}\left(  \left\langle \varphi_\eps , \frac{1-\eta_\eps}{\rho_\eps}\left(\hat{\varphi}_\eps - \Psi_\eps\right) \right\rangle_{\hat{\Lambda}_\eps} \right) \\
& = L_\eps\left( \left\vert \hat{\varphi}_\eps - \Psi_\eps  \right\vert_{\hat{\Lambda}_\eps}^2\right) + L_\eps\left((1-\eta_\eps)^2 \frac{ \vert \tau_\eps \vert^2 - \vert \hat{\theta}_\eps \vert^2 }{\omega_\eps^2 - \vert \tau_\eps \vert^2 }  \right)
\end{split}
\end{equation*}
where for the last equality, we used that $\langle X,(X-Y)\rangle_\Lambda = \frac{1}{2} \left\vert X-Y \right\vert^2_\Lambda + \frac{1}{2} \left(\left\vert X \right\vert^2_\Lambda - \left\vert Y \right\vert^2_\Lambda \right)$ with $X = (1-\eta_\eps)  \frac{\varphi_\eps}{\rho_\eps}$, $Y = (1-\eta_\eps) \Psi_\eps$ and the equality
$$ \hat{\varphi}_\eps - \Psi_\eps = (1-\eta_\eps) \left( \frac{\varphi_\eps}{\rho_\eps}-\Psi_\eps \right). $$

We obtain that 
\begin{equation*}
\begin{split}  \int_{\mathbb{D}_{r_\eps(p)}(p)}  \left\vert \nabla \hat{\varphi}_\eps \right\vert^2 - \int_{\mathbb{D}_{r_\eps(p)}(p)} \left\vert \nabla  \Psi_\eps \right\vert^2  \leq  L_\eps\left((1-\eta_\eps)^2 \frac{ \vert \tau_\eps \vert^2 - \vert \hat{\theta}_\eps \vert^2 }{\omega_\eps^2 - \vert \tau_\eps \vert^2 }  \right). \\
+ 2 \int_{\mathbb{D}_{r_\eps(p)}(p)} \left(  \nabla \hat{\varphi}_\eps  \nabla\left( \hat{\varphi}_\eps - \Psi_\eps  \right) -\nabla \varphi_\eps \nabla \left( \frac{1-\eta_\eps}{\rho_\eps}\left( \hat{\varphi}_\eps - \Psi_\eps  \right) \right) \right) = I + II
\end{split}
\end{equation*}
The first right-hand term satisfies by a Cauchy-Schwarz inequality and properties of $\lambda_{\star}:= \lambda_{\star}\left(\mathbb{A}_{r_\eps(p),s_\eps}(p),\beta_\eps \right)$
\begin{equation*}
\begin{split} I^2 \leq & 4 L_\eps\left( \left\vert (1-\eta_\eps)(\tau_\eps - \hat{\theta}_\eps) \right\vert^2 \right) L_\eps\left( \left\vert (1-\eta_\eps)(\tau_\eps + \hat{\theta}_\eps) \right\vert^2 \right) \\
\leq & C \frac{1}{\lambda_\star} \int_{\mathbb{D}_{r_\eps(p)}(p)} \left\vert \nabla \left((1-\eta_\eps)(\tau_\eps - \hat{\theta}_\eps)   \right) \right\vert^2 \leq o(1) 
 \end{split} \end{equation*}
as $\eps \to 0$ since the energies of $\hat{\theta}_\eps$, $\tau_\eps$ and $\eta_\eps$ go to $0$ as $\eps\to 0$. The second right-hand term satisfies
\begin{equation*}
\begin{split}  II = &  2 \int_{\mathbb{D}_{r_\eps(p)}(p)} \nabla\left( \hat{\varphi}_\eps - \psi_\eps\right)\nabla \left( \hat{\varphi}_\eps - \varphi_\eps\frac{(1-\eta_\eps)}{\rho_\eps} \right)  \\
& + 2 \int_{\mathbb{D}_{r_\eps(p)}(p)} \nabla\frac{\eta_\eps}{\rho_\eps} \left( \left(\hat{\varphi}_\eps-\psi_\eps\right) \nabla \varphi_\eps - \varphi_\eps \nabla \left(\hat{\varphi}_\eps - \Psi_\eps \right) \right)  \\
 = & 2\int_{\mathbb{D}_{r_\eps(p)}(p)} \nabla\left( \hat{\varphi}_\eps - \psi_\eps\right)\nabla \left( \eta_\eps \Psi_\eps \right) \\
 & + 2 \int_{\mathbb{D}_{r_\eps(p)}(p)} \left(\nabla \eta_\eps - \eta_\eps \frac{\nabla \rho_\eps}{\rho_\eps}\right) \left( \left(\hat{\varphi}_\eps-\psi_\eps\right) \frac{\nabla \varphi_\eps}{\rho_\eps} - \frac{\varphi_\eps}{\rho_\eps} \nabla \left(\hat{\varphi}_\eps - \Psi_\eps \right) \right)  \\
 \leq & C \left(\int_{\mathbb{D}_{r_\eps(p)}(p)} \eta_\eps^2 \vert \nabla \Psi_\eps \vert^2 + \vert \nabla \eta_\eps \vert^2  \right)^{\frac{1}{2}} + C \left( \int_{\mathbb{D}_{r_\eps(p)}(p)}  \vert \nabla \rho_\eps \vert^2 + \vert \nabla \eta_\eps \vert^2 \right)^{\frac{1}{2}} = o(1)
\end{split}\end{equation*}
as $\eps\to 0$ where we used for the inequality that the energy of $\varphi_\eps$ and $\hat{\varphi}_\eps-\Psi_\eps$ is uniformly bounded, that $\rho_\eps^{-1}$, $\frac{\varphi_\eps}{\rho_\eps}$ and $\hat{\varphi}_\eps-\Psi_\eps$ are uniformly bounded in $L^{\infty}$ as $\eps \to 0$. For the last equality, we use that the energy of $\rho_\eps$ and $\eta_\eps$ converges to $0$, and that the $L^\infty$ norm of $\vert \nabla \Psi_\eps \vert^2 $ is uniformly bounded in $\mathbb{D}_{\frac{r_\eps(p)}{2}}(p)$ by $\eps$-regularity on harmonic maps (see Claim \ref{cl:epsregclosedcase}). Finally we obtain \eqref{eqsmalldiffenergy}

\medskip
\noindent \textbf{Step 2:}
\medskip

$$ \int_{\mathbb{D}_{r_\eps(p)}(p)}   \left\vert \nabla \left( \Psi_\eps- \hat{\varphi}_\eps^j \right) \right\vert^2 \leq  \int_{\mathbb{D}_{r_\eps(p)}(p)}  \left\vert \nabla \hat{\varphi}_\eps \right\vert^2 - \int_{\mathbb{D}_{r_\eps(p)}(p)} \left\vert \nabla  \Psi_\eps \right\vert^2 + o(1) $$
as $\eps \to 0$.

\medskip
\noindent \textbf{Proof of Step 2:}
\medskip

We test the equation on $\Psi_\eps$: $\Delta \Psi_\eps = \frac{\left\vert \nabla \Psi_\eps \right\vert_{\Lambda_\eps}^2}{\left\vert \Lambda_\eps\Psi_\eps \right\vert^2} \Lambda_\eps \Psi_\eps$ against $\Psi_\eps - \hat{\varphi}_\eps$ and we multiply by $2$ to obtain
\begin{equation*}
\begin{split} 2 & \int_{\mathbb{D}_{r_\eps(p)}(p)}  \nabla \Psi_\eps \nabla \left( \Psi_\eps - \hat{\varphi}_\eps  \right) = 2 \int_{\mathbb{D}_{r_\eps(p)}(p)}  \frac{\left\vert \nabla \Psi_\eps \right\vert_{\hat{\Lambda}_\eps}^2}{\left\vert \hat{\Lambda}_\eps\Psi_\eps \right\vert^2} \langle \Psi_\eps, \Psi_\eps - \hat{\varphi}_\eps \rangle_{\hat{\Lambda}_\eps}  \\ 
= &  \int_{\mathbb{D}_{r_\eps(p)}(p)}  \frac{\left\vert \nabla \Psi_\eps \right\vert_{\hat{\Lambda}_\eps}^2}{\left\vert \hat{\Lambda}_\eps\Psi_\eps \right\vert^2} \left( \left\vert \Psi_\eps - \hat{\varphi}_\eps \right\vert_{\Lambda_\eps}^2 + (1-\eta_\eps) \frac{  \vert \hat{\theta}_\eps \vert^2 - \vert \tau_\eps \vert^2  }{\omega_\eps^2 - \vert \tau_\eps \vert^2 }  \right)   \\
\leq & C \left(\frac{ \lambda_K^\eps}{ \lambda_{I+1}^\eps}\right)^2\eps_0 \left(  \int_{\mathbb{D}_{r_\eps(p)}(p)}  \left\vert \nabla \left( \Psi_\eps - \hat{\varphi}_\eps\right) \right\vert^2 + \left(\int_{\mathbb{D}_{r_\eps(p)}(p)} \left\vert \nabla \left( \hat{\theta}_\eps-\tau_\eps \right) \right\vert^2 + \vert \nabla \eta_\eps \vert^2\right)^{\frac{1}{2}}  \right)
\end{split} \end{equation*}
where we used again that $\langle X,(X-Y)\rangle_\Lambda = \frac{1}{2} \left\vert X-Y \right\vert^2_\Lambda + \frac{1}{2} \left(\left\vert X \right\vert^2_\Lambda - \left\vert Y \right\vert^2_\Lambda \right)$ with $X = \Psi_\eps$ and $Y = \Psi_\eps - \hat{\varphi}_\eps$ for the second equality. The first inequality is a consequence of the rescaling on $\mathbb{D}_{r_\eps(p)(p)}$ of the following classical Hardy inequality (see e.g \cite{laurainpetrides}, Theorem 3.1)
$$\forall u \in H_0^1(\mathbb{D}), \frac{1}{4} \int_{\mathbb{D}} \frac{u^2}{\left(1-\vert x \vert\right)^2} \leq  \int_\mathbb{D} \vert \nabla u \vert^2 $$
using the $\eps$-regularity of the energy of harmonic maps coming from Claim \ref{cl:epsregclosedcase}, we have
$$ \vert \nabla \Psi_\eps \vert^2(x) \leq  \frac{C}{\left(r_\eps(p)-\vert x - p \vert\right)^2} \int_{\mathbb{D}_{r_\eps(p)}(p)}  \vert \nabla \Psi_\eps \vert^2.$$
Then, we have that
\begin{equation*}
\begin{split} & \int_{\mathbb{D}_{r_\eps(p)}(p)}    \left\vert \nabla \left( \Psi_\eps- \hat{\varphi}_\eps^j \right) \right\vert^2 =  \int_{\mathbb{D}_{r_\eps(p)}(p)}\left( \left\vert \nabla \hat{\varphi}_\eps \right\vert^2 -  \left\vert \nabla  \Psi_\eps \right\vert^2 + 2    \nabla \Psi_\eps \nabla \left( \Psi_\eps - \hat{\varphi}_\eps  \right)\right) \\ 
& \leq \int_{\mathbb{D}_{r_\eps(p)}(p)}  \left\vert \nabla \hat{\varphi}_\eps \right\vert^2 - \int_{\mathbb{D}_{r_\eps(p)}(p)} \left\vert \nabla  \Psi_\eps \right\vert^2 + C' \eps_0 \int_{\mathbb{D}_{r_\eps(p)}(p)}   \left\vert \nabla \left( \Psi_\eps- \hat{\varphi}_\eps^j \right) \right\vert^2 + o(1)    \end{split} \end{equation*}
as $\eps \to 0$. Choosing $\eps_0 \leq \left(2C'\right)^{-1}$, we obtain Step 2 and the Claim.
\end{proof}

\subsection{Convergence results on the Palais-Smale sequence}

We consider $\widetilde{\Sigma} := \Sigma \sqcup \bigsqcup_{j=1}^l \left(\mathbb{S}^2\right)_j$ endowed with the metric $\tilde{g}$ equal to $g$ on $\Sigma$ and the round metric $g_{\mathbb{S}^2}$ on $(\mathbb{S}^2)_j$ for $1\leq j\leq l$. 
Thanks to the previous claims, we can construct a covering of $\widetilde{\Sigma}$ of disks $\{ \mathbb{D}_{r_\eps(p)}(p) \}_{p \in Q}$ where $Q$ is a finite set independent of $\eps$
such that the conclusions of Claim \ref{clcomparegeneral} hold on any $\mathbb{D}_{r_\eps(p)}(p)$. We use this property to localize and prove the following:

\begin{cl}
There is $V_0 \in L^\infty_+(\Sigma)$ and $V_1,\cdots,V_j \in L^\infty_+(\mathbb{S}^2)$ such that for any $\eta_0 \in \mathcal{C}^\infty_c\left( \Sigma \setminus P_0 \right)$ and $\eta_j \in \mathcal{C}^\infty_c\left( \mathbb{S}^2 \setminus P_j \right)$ 
for $0 \leq j \leq l$,
\begin{equation}\label{eqconvbetaepstoV} \beta_j^\eps(\eta_j, 1) - \int \eta_j V_j \leq o(1) \left(\Vert \nabla \eta \Vert_{L^2} + \Vert \eta \Vert_{L^\infty} \right).\end{equation}
as $\eps\to 0$. In particular $\mu_0 = V_0 dA_g$ and $\mu_j = V_j dA_{\mathbb{S}^2}$ for $1 \leq j\leq l$
\end{cl}

\begin{proof}
We prove the result for a given $0 \leq j\leq l$ and we drop the use of $j$ in the indices/exponents of functions. We localize the result: let $\eta$ be a cut-off function at the neighborhood of a good point such that a harmonic replacement given by Claim \ref{clharmreplacement} is well-defined on $K = supp(\eta)$ for any large $\eps$, and such that for any large $\eps$,
$$ \Vert \left\vert \nabla \Psi_\eps \right\vert^2 \Vert_{L^{\infty}(K)} \leq A $$
for some constant $A$ by $\eps$-regularity of harmonic maps in Claim \ref{cl:epsregclosedcase}. Then, $ \frac{\left\vert \nabla \Psi_\eps \right\vert_{\hat{\Lambda}_\eps}^2}{\left\vert \hat{\Lambda}_\eps \Psi_\eps \right\vert^2} $ converges to some function $V_j \in L^\infty(K)$ strongly in $L^p(K)$ for $1\leq p < +\infty$. 

We test the function $\frac{\eta \varphi_\eps}{\rho_\eps^2}$ against the equation on $\varphi_\eps$: $\Delta \varphi_\eps = \hat{\sigma}_\eps\beta_\eps(\varphi_\eps,\cdot)$. We obtain
\begin{equation*}\begin{split} \beta_\eps(1,\eta)  = & \beta_\eps\left(\frac{\vert\varphi_\eps\vert_{\hat{\sigma}_\eps}^2}{\rho_\eps^2},\eta\right) = \hat{\Lambda}_\eps\beta_\eps\left(\varphi_\eps,\frac{\varphi_\eps \eta}{\rho_\eps^2}\right) = \int_{K} \nabla \varphi_\eps \nabla \frac{\varphi_\eps \eta}{\rho_\eps^2} \\
 = & \int_K  \frac{\varphi_\eps}{\rho_\eps}  \nabla \frac{\varphi_\eps}{\rho_\eps}  \nabla \eta - \int_K  \frac{\vert\varphi_\eps\vert^2}{\rho_\eps} \nabla \frac{1}{\rho_\eps} \nabla \eta  \\ 
& + \int_K \eta \left\vert \nabla \frac{\varphi_\eps}{\rho_\eps} \right\vert^2 +  \int_K \eta \nabla\frac{1}{\rho_\eps}\left( \frac{\varphi_\eps}{\rho_\eps}  \nabla \varphi_\eps - \varphi_\eps \nabla\frac{\varphi_\eps}{\rho_\eps}\right)    \\
 = & \int_K \left( \eta \vert \nabla \Psi_\eps \vert^2 + \Psi_\eps \nabla \Psi_\eps  \nabla \eta\right) + \int_K  \eta\left( \left\vert \nabla \frac{\varphi_\eps}{\rho_\eps} \right\vert^2 - \vert \nabla \Psi_\eps \vert^2\right)  \\
& + \int_K \left( \Psi_\eps \nabla \Psi_\eps - \frac{\varphi_\eps}{\rho_\eps}\nabla \frac{\varphi_\eps}{\rho_\eps}\right)  \nabla \eta  + \int_K \nabla \frac{1}{\rho_\eps}\left( \frac{\vert\varphi_\eps\vert^2}{\rho_\eps}\nabla \eta + \eta\left( \frac{\varphi_\eps}{\rho_\eps}  \nabla \varphi_\eps - \varphi_\eps \nabla\frac{\varphi_\eps}{\rho_\eps}  \right)  \right) \\
= & \int_K \nabla (\eta\Psi_\eps)\nabla \Psi_\eps + o(1) \left(\Vert \nabla \eta \Vert_{L^2} + \Vert \eta \Vert_{L^\infty} \right) \\
= & \int_{K} \eta \frac{\left\vert \nabla \Psi_\eps \right\vert_{\hat{\Lambda}_\eps}^2}{\left\vert \hat{\Lambda}_\eps \Psi_\eps \right\vert^2}  +o(1)\left(\Vert \nabla \eta \Vert_{L^2} + \Vert \eta \Vert_{L^\infty} \right)
\end{split}\end{equation*}
where the penultimate equality comes from Claim \ref{clcomparegeneral} and \eqref{eq:smallenergiesetarho}. We completed the proof.

\end{proof}

We recall that for a Riemannian surface $(\Sigma,g)$,
$$ I_F(\Sigma,g)  = \inf_{\beta \in \bar{X}} F(\bar{\lambda}_1(\Sigma,g,\beta),\cdots,\bar{\lambda}_m(\Sigma,g,\beta)) $$
From the previous claim, we obtain a measure $VdA_{\tilde{g}}$ equal to $V_0dA_g$ on $\Sigma$ and $V_j dA_{\mathbb{S}^2}$ on $(\mathbb{S}^2)_j$ for $1\leq j\leq l$.
By upper semi-continuity of eigenvalues with respect to bubble tree convergence, and then lower semi-continuity of $f(\Sigma,g,\beta) := F(\bar{\lambda}_1(\Sigma,g,\beta),\cdots,\bar{\lambda}_m(\Sigma,g,\beta))$ with respect to bubble tree convergence, we obtain that
$$ I_F(\Sigma,g) = \liminf_{\eps\to 0} E(\Sigma,g,\beta_\eps) \geq E(\widetilde{\Sigma},\widetilde{g}, V dA_{\tilde{g}} )  ) \geq I_{F}(\widetilde{\Sigma}, \widetilde{g})  $$
In addition, we know by glueing methods that 
$I_{F}(\widetilde{\Sigma}, \widetilde{g}) \geq  I_F(\Sigma,g) $ (see \cite{ces}). Therefore, all the inequalities are equalities and $V dA_{\tilde{g}}$ is a minimizer for $E$ on $\left(\widetilde{\Sigma},\widetilde{g}\right)$. 

By Euler-Lagrange equation applied to the minimizer $V dA_{\tilde{g}}$ (see Proposition \ref{prop:constructPSK} for $\eps = 0$), we obtain the existence of $\Phi : \widetilde{\Sigma} \to \mathbb{R}^n$ such that setting $\lambda_k := \lambda_k(\widetilde{\Sigma},\widetilde{g}, VdA_{\widetilde{g}})$, and $\Lambda := (\lambda_1,\cdots,\lambda_n)$
\begin{itemize}
\item $\Delta_{\widetilde{g}} \Phi = \Lambda V \Phi$
\item $\vert \Phi \vert^2_\Lambda \geq 1$ and $\int_{\widetilde{\Sigma}} \vert \Phi \vert^2_\Lambda V dA_{\widetilde{g}} = 1 $.
\end{itemize}
Applying Claim \eqref{eqomegaepsto1} with $\theta_\eps = 0$, we obtain that $\vert \Phi \vert^2_\Lambda = 1$, so that $\Phi : \widetilde{\Sigma} \to \mathcal{E}_{\Lambda}$ is a harmonic map. In addition, we have by the computation of $\frac{1}{2} \Delta_{\widetilde{g}} \vert \Phi \vert^2_\Lambda = 0$, we obtain that
$$ V = \frac{\left\vert \nabla \Phi \right\vert^2_\Lambda}{\vert \Lambda\Phi \vert^2} $$
and since a harmonic map has to be smooth, $V$ is a smooth function and vanishes at most at a finite number of points, that correspond to conical singularities of $V \widetilde{g}$. The proof of Proposition \ref{palaissmalegeneral} is complete.

\section{Convergence of regularized minimizing sequences in the Steklov case}
\label{sec3}
We aim at proving the following proposition (see definition \ref{defi:mwstar} for the MW$\star$ bubble tree convergence where we take measures that have their suport in $\partial \Sigma$ and we replace surfaces $\Sigma$ by curves $\partial \Sigma$ and $\mathbb{R}^2$ by $\mathbb{R}$, that is the stereographic projection of $\mathbb{S}^1$). Since the proof is very smilar to the closed case, we will often drop portions of proof that do not differ to the closed case and we will emphasize on the main differences.

\begin{prop} \label{palaissmalegeneralsteklov}
Let $(\Sigma,g)$ be a Riemannian surface with a boundary and $(\beta_\eps, \Phi_\eps,g_\eps)$, be a $(PS)_{K}$ sequence as $\eps\to 0$.
Then, up to the extraction of a subsequence 
$e^{2u_\eps}dL_g$ and $\beta_\eps(1,.)$ MW$\star$--bubble tree converge to the measures $V_0 dL_g$ (possibly $0$ if $l\geq 1$) on $\partial\Sigma$ and $V_j dL_{\mathbb{S}^1}$ on $\left(\mathbb{S}^1\right)_j$ where $V_0,V_1,\cdots,V_l$ are $L^{\infty}$ densities.

If in addition $(\beta_\eps)$ and $(g_\eps)$ are minimizing sequences for $E$, then 
$$ V_0 = \Phi_0 \cdot \partial_\nu \Phi_0 \text{ and } V_j = \Phi_j \cdot \partial_r \Phi_j  $$  
where $\Phi_0 : (\Sigma,\partial\Sigma) \to (co\left(\mathcal{E}_\sigma\right),\mathcal{E}_\sigma)$ and $\Phi_j : \left(\mathbb{D},\mathbb{S}^1\right)_j \to (co\left(\mathcal{E}_\sigma\right),\mathcal{E}_\sigma)$ are free boundary harmonic maps harmonic maps into $co\left(\mathcal{E}_\sigma\right)$ and we have that
$$ I_F(\Sigma,[g]) = I_F(\widetilde{\Sigma},[\widetilde{g}]) $$
where $\widetilde{\Sigma} = \Sigma \sqcup \bigsqcup_{j=1}^l (\mathbb{D})_j$ endowed with $\widetilde{g}$ equal to $g$ on $\Sigma$ and the flat metric on the copies of $\mathbb{D}$.
\end{prop}

\begin{rem}
Notice that by a glueing method similar to \cite{ces} or \cite{fs3}, we always have 
$$ I_F(\Sigma,[g]) \leq  I_F(\widetilde{\Sigma},[\widetilde{g}]) $$
and if we know that the inequality is strict, then we automatically deduce that all the minimizing sequences for $I_F(\Sigma,[g])$ MW$\star$ converge to a measure absolutely continuous with respect to $dL_g$ with a smooth density ($l=0$ in the proposition)
\end{rem}

This proposition and the remark proves Theorem \ref{theo:main} in the case of Steklov eigenvalues, noticing that if $\tilde{V}$ is a positive extension of $V$ in $\widetilde{\Sigma}$, $V \tilde{g}$ is a smooth metric (contrary to the closed case, conical singularities are not possible on the boundary by a classical use of a Hopf lemma coming from the maximum principle)

\medskip

In this case, we have the following notations: if $\Omega \subset \Sigma$ is an open set, we denote the surface boundary of $\Omega$:
$$ \partial_s \Omega := \partial \overline{\Omega} \cap \partial \Sigma $$
and the domain boundary of $\Omega$
$$ \partial_d \Omega := \partial \overline{\Omega} \setminus \partial \Sigma $$
and $H^1_0(\Omega,g)$ is the set of $H^1$ functions of $\Omega$ equal to $0$ on the domain boundary of $\Omega$: $\partial_d \Omega$

\subsection{Tree of concentration points }

As in the closed case, we currently have the following general property that is similar to \eqref{prop:bubbletree}. The notation for $I$ is kept all along the proof.
\begin{prop} \label{prop:bubbletreesteklov}
Let $(\beta_\eps ,\Phi_\eps,e^{2u_\eps}g)$ be a $(PS)_K$ sequence. We assume in addition that $\sigma_\eps^I \to 0$ and that $\left(\sigma_\eps^{I+1}\right)$ is uniformly lower bounded as $\eps\to 0$. Then up to the extraction of a subsequence, $\beta_\eps$ and $e^{u_\eps}dL_g$ MW$\star$-bubble tree converge to the same measures $\mu_0$ on $\partial\Sigma$ and $\mu_j$ in $\left(\mathbb{S}^1\right)_j$ for $1\leq j\leq l$ where $l\leq I+1$. 
\end{prop}

Again, there is an abuse of notation with the use of sums on manifolds. Here, we work on an atlas of conformal charts such that if the chart intersects the boundary, $\Sigma$ is locally isometric to a portion of the half-space $\mathbb{R}^2_+ := \mathbb{R}\times \mathbb{R}_+ \cap U$ endowed with a metric conformal to the flat metric, such that $\mathbb{R}\times \{0\} \cap U$ corresponds to the boundary of $\Sigma$. Then, if $p \in \partial \Sigma$, we denote $\mathbb{D}_r^+(p)$ the Euclidean half balls in the charts centered at $p$.

We denote for $1\leq j\leq l$, $(x_j^\eps,\alpha_j^\eps)$ the associated points and scales. We denote for $0\leq j \leq l$, $\mu_j$ the pullback of the continuous part of $\nu_j$ (having the set of atoms $Z_j$) with repect to $\pi_{\mathbb{S}^1}$, the stereographic projection $\mathbb{S}^1 \to \mathbb{R}$ (restriction to $\mathbb{S}^1$ of a biholomorphism $\mathbb{D}\to \mathbb{R}^2_+$). The functions $f_\eps : \Sigma \to \R$ we consider are seen at the scale $(x_j^\eps,\alpha_j^\eps) \in \partial \Sigma \times \mathbb{R}_+^\star$ with the formula
$$ f_j^\eps := f_\eps( x_j^\eps + \alpha_j^\eps \pi_{\mathbb{S}^1}^{-1}\left(\cdot\right) ) $$
and in particular, we denote $\widetilde{\Phi_\eps}^j := \left(\Phi_\eps\right)_j$
while linear forms on continuous functions (measures) $\mu_\eps$ or bilinear forms on $H^1$ functions $\beta_\eps$ satisfy at the scale $(x_j^\eps,\alpha_j^\eps)$ for $\varphi,\psi \in \mathcal{C}_c^\infty\left( \mathbb{\pi}_{\mathbb{S}^1}^{-1}\left( \R \setminus Z_j \right) \right)$
$$ \langle \mu_j^\eps , \varphi \rangle := \left\langle \mu_\eps , \varphi\left( \frac{ \pi_{\mathbb{S}^1}\left(\cdot\right) - x_j^\eps }{\alpha_j^\eps} \right) \right\rangle $$
$$ \beta_j^\eps(\varphi,\psi):= \beta_\eps \left(\varphi\left( \frac{\pi_{\mathbb{S}^1}\left(\cdot\right) - x_j^\eps}{\alpha_j^\eps} \right),\psi\left( \frac{ \pi_{\mathbb{S}^1}\left(\cdot\right) - x_j^\eps }{\alpha_j^\eps} \right)\right)$$
and in particular, we denote $e^{u_j^\eps}dL_{\mathbb{S}^1} := \left( e^{u_\eps} dL_g \right)_j $. 

We say that the analysis in $\Sigma$ if $\mu_0 \neq 0$ and in $\left(\mathbb{D}\right)_j$ for $1\leq j \leq l$ (in this case $\mu_j \neq 0$) at the scale $(x_j^\eps,\alpha_j^\eps)$ of functions, measures, bilinear forms and sets we consider is an analysis in a "thick part" since the measure has a positive mass at this scale.

We also localize the space $\bar{X}$. We denote for an open set $\Omega$ of a smooth Riemannian surface with boundary $(\Sigma,g)$, $\overline{X}\left(\Omega,g\right)$ the closure of 
$$ X(\Omega,g) := \{ \left(\varphi,\psi\right)\mapsto \int_{\partial \Sigma} \varphi \psi e^{u}dL_g ; u\in \mathcal{C}^\infty(\partial\Sigma \cap \Omega) \} $$
in the set of symmetric bilinear forms on $H^{1}_0(\Omega,g)$ endowed with the norm
$$ \Vert \beta \Vert_{\overline{X}\left(\Omega,g\right)} := \sup_{\varphi,\psi \in H^1_0(\Omega) } \frac{\vert \beta(\varphi,\psi)\vert}{\Vert \varphi \Vert_{H^1_0(\Omega,g)} \Vert \psi \Vert_{H^1_0(\Omega,g)}}, $$
where
$$ \Vert \varphi \Vert_{H^1_0(\Omega,\partial,g)}^2 := \int_{\Omega} \vert \nabla \varphi \vert^2_g dA_g + \int_{\partial_s \Omega} \varphi^2 dL_g. $$

The goal in all the section is to prove that the limiting measures $\mu_0,\cdots,\mu_l$ are absolutely continous with respect to $dL_g$ or $dL_{\mathbb{S}^1}$ (the Lebesgue length measure of $\mathbb{S}^1$) with densities satisfying the conclusions of Proposition \ref{palaissmalegeneral}.

\subsection{Some convergence of $\omega_{\eps}$ to $1$ and first replacement of $\Phi_{\eps}$} 

We set $\omega_\eps$ the harmonic extension of the following map defined on $\partial \Sigma$
$$ \omega_{\eps} = \sqrt{\left\vert \Phi_{\eps} \right\vert_{\sigma_{\eps}}^2 + \theta_\eps^2 } 
\text{ in } \partial \Sigma \text{ and } \Delta_g \omega_\eps = 0 \text{ in } \Sigma$$ 
We first prove that $\nabla \omega_{\eps}$ converges to $0$ in $L^2$ and that $\sqrt{\sigma_\eps}\Phi_{\eps}$ has a similar $H^1$ behaviour as $\frac{\sqrt{\sigma_\eps}\Phi_{\eps}}{\omega_{\eps}}$

\begin{cl} \label{clomegaeps2steklov} We have that
\begin{equation}  \label{eqomegaepsto1steklov} \int_{\Sigma}(\omega_\eps^2-1)\left\vert \nabla \frac{\Phi_\eps}{\omega_\eps}  \right\vert_{\sigma_\eps}^2 dA_g +  \int_{\Sigma}  \left\vert \nabla \omega_\eps \right\vert^2 + \int_\Sigma \left\vert \nabla\left( \Phi_\eps - \frac{\Phi_\eps}{\omega_\eps} \right) \right\vert_{\sigma_\eps}^2  \leq O(\eps) \end{equation}
as $\eps\to 0$.
\end{cl}

The proof is similar to the proof of Claim  \ref{eqomegaepsto1steklov} but needs a particular attention because of the harmonic extension of $\omega_\eps$

\begin{proof} We first prove
\begin{equation} \label{eqomega_epsminus1steklov} L_\eps\left( \left\vert \sigma_\eps \Phi_\eps \right\vert^2 \left(1- \frac{1 }{\omega_\eps}\right) \right) \leq O(\eps) \end{equation}
as $\eps\to 0$. Since $\omega_\eps \geq 1$,  
and $\left\vert \Phi_\eps \right\vert_{\sigma_\eps}^2 \leq \omega_\eps^2$, we have that
\begin{equation*}
\begin{split}
 L_\eps\left( \left\vert \sigma_\eps \Phi_\eps \right\vert^2 \left(1- \frac{1 }{\omega_\eps} \right) \right) \leq & \sigma_K^\eps L_\eps\left(  \left(\omega_\eps^2-\omega_\eps\right) \right) \\
\leq & \sigma_K^\eps   \left( L_\eps\left(  \left\vert \Phi_\eps \right\vert_{\sigma_\eps}^2 \right) + L_\eps \left(  \theta_\eps^2 \right) - L_{\eps}(1) \right) \\
 \end{split} \end{equation*}
 so that
$$ L_\eps\left( \left\vert \sigma_\eps \Phi_\eps \right\vert^2 \left(1- \frac{1 }{\omega_\eps}\right) \right) \leq  \sigma_K^\eps   L_\eps( \theta_\eps^2) \leq  \sigma_K^\eps  \left\Vert \beta_\eps \right\Vert_{g_\eps} \left\Vert \theta_\eps \right\Vert^2_{H^1(\partial,g_\eps)} \leq O(\eps) $$
as $\eps \to 0$ since by assumption
$$ \left\Vert \beta_\eps \right\Vert_{g_\eps}  \leq \left\Vert e^{u_\eps}dL_g \right\Vert_{g_\eps} +\eps \leq \sup_{\varphi,\psi \in H^1} \frac{\vert \int_{\partial \Sigma}e^{u_\eps}\varphi\psi dL_g \vert}{\Vert \varphi \Vert_{H^1(g_\eps)} \Vert \psi \Vert_{H^1(g_\eps)}} +\eps \leq 1 +\eps. $$
by and we obtain \eqref{eqomega_epsminus1steklov}. We now prove \eqref{eqomegaepsto1steklov}:
\begin{equation*} 
\begin{split}
\int_{\Sigma} & \left\vert \nabla \frac{\Phi_\eps}{\omega_\eps} \right\vert_{\sigma_\eps}^2 - \int_\Sigma  \left\vert \nabla \Phi_\eps \right\vert_{\sigma_\eps}^2 - \int_\Sigma  \left\vert \nabla\left( \Phi_\eps - \frac{\Phi_\eps}{\omega_\eps} \right) \right\vert_{\sigma_\eps}^2 
\\ = & - 2 \int_\Sigma \left\langle \nabla \Phi_\eps, \nabla\left( \Phi_\eps - \frac{ \Phi_\eps}{\omega_\eps}\right) \right\rangle_{\sigma_\eps} 
= - 2 \int_\Sigma \Delta \Phi_\eps  {\sigma_\eps}.\left( \Phi_\eps - \frac{ \Phi_\eps}{\omega_\eps}\right)  \\
= & - 2 \beta_\eps\left( \sigma_\eps.\Phi_\eps , {\sigma_\eps}.\left( \Phi_\eps - \frac{ \Phi_\eps}{\omega_\eps}\right) \right)
= - 2 L_\eps\left(   \left\vert \sigma_\eps \Phi_\eps \right\vert^2 \left(1- \frac{1 }{\omega_\eps} \right) \right) = O(\eps)
\end{split}
\end{equation*}
where we tested $ \Delta \Phi_{\eps} =  \beta_\eps( \sigma_\eps\Phi_\eps,.) $ in $\Sigma$ against $ \sigma_\eps. \left( \Phi_\eps-\frac{ \Phi_\eps}{\omega_\eps}\right)$, and we used \eqref{eqomega_epsminus1steklov}.

In particular, we have
$$ 0 \leq \int_\Sigma \left\vert \nabla\left( \Phi_\eps - \frac{\Phi_\eps}{\omega_\eps} \right) \right\vert_{\sigma_\eps}^2 \leq \int_\Sigma  \left( \left\vert \nabla \frac{\Phi_\eps}{\omega_\eps} \right\vert_{\sigma_\eps}^2 - \left\vert \nabla \Phi_\eps \right\vert_{\sigma_\eps}^2\right) + O(\eps) $$
as $\eps\to 0$ and knowing that with the straightforward computations we have
\begin{equation*}
\begin{split} \left\vert \nabla \frac{\Phi_\eps}{\omega_\eps} \right\vert_{\sigma_\eps}^2 - \left\vert \nabla \Phi_\eps \right\vert_{\sigma_\eps}^2 = & \left( 1 - \omega_\eps^2  \right)\left\vert \nabla \frac{\Phi_\eps}{\omega_\eps} \right\vert_{\sigma_\eps}^2 - \left( \left\vert \nabla \omega_\eps \right\vert^2 \frac{\left\vert \Phi_\eps \right\vert^2_{\sigma_\eps}}{\omega_\eps^2}  + \omega_\eps \nabla \omega_\eps \nabla \frac{\left\vert \Phi_\eps \right\vert^2_{\sigma_\eps}}{\omega_\eps^2} \right) \\
=  & \left( 1 - \omega_\eps^2  \right)\left\vert \nabla \frac{\Phi_\eps}{\omega_\eps} \right\vert_{\sigma_\eps}^2 - \nabla \omega_\eps \nabla  \frac{\left\vert \Phi_\eps \right\vert^2_{\sigma_\eps}}{\omega_\eps} 
\end{split} \end{equation*}
Computing that
\begin{equation*}
\begin{split} \int_\Sigma  \nabla \omega_\eps, \nabla  \frac{\left\vert \Phi_\eps \right\vert^2_{\sigma_\eps}}{\omega_\eps}  & = \int_{\partial\Sigma} \partial_\nu \omega_\eps \left(\omega_\eps - \frac{\theta_\eps^2}{\omega_\eps}\right)  \\
& = - \int_{\Sigma} \Delta \frac{\omega_\eps^2}{2} - \int_{\Sigma} \nabla \omega_\eps \nabla \frac{\theta_\eps^2}{\omega_\eps} \\
& = \int_{\Sigma} \left\vert \nabla \omega_\eps \right\vert^2 + \int_{\Sigma} \frac{\theta_\eps^2}{\omega_\eps^2} \left\vert \nabla \omega_\eps \right\vert^2 - 2 \int_{\Sigma} \frac{\theta_\eps}{\omega_\eps} \nabla \theta_\eps \nabla \omega_\eps \\
& \geq \int_{\Sigma} \left\vert \nabla \omega_\eps \right\vert^2 -  \int_{\Sigma} \left\vert \nabla \theta_\eps \right\vert^2
\end{split} \end{equation*}
and we obtain since $\frac{\theta_\eps}{\omega_\eps}$ is uniformly bounded by $1$ that
\begin{equation*}
\begin{split} \int_\Sigma \left\vert \nabla\left( \Phi_\eps - \frac{\Phi_\eps}{\omega_\eps} \right) \right\vert_{\sigma_\eps}^2 + \int_{\Sigma} \left( \omega_\eps^2 - 1 \right)\left\vert \nabla \frac{\Phi_\eps}{\omega_\eps} \right\vert_{\sigma_\eps}^2 + \int_{\Sigma} \left\vert \nabla \omega_\eps \right\vert^2 \leq O\left(\eps \right)  
\end{split} \end{equation*}
as $\eps \to 0$,


\end{proof}

\subsection{Good/bad points in thick parts and immediate consequences}

\subsubsection{Construction of a finite number of bad points}

In the following, we perform local regularity estimates on $(\Phi_\eps)$. These estimates can only be done far from "bad points" we select in Claim \ref{clbadpointssteklov}.
For $\Omega \subset \Sigma$ a domain of $\Sigma$, we recall that
$$ \sigma_\star(\Omega, g, \beta_\eps) = \inf_{\varphi \in \mathcal{C}_c^{\infty}\left(\Omega\right) } \frac{\int_\Omega  \left\vert \nabla \varphi \right\vert_g^2 dA_g }{\beta_\eps(\varphi,\varphi) }. $$
We recall that $\sigma_K^\eps := \max_{i\in\{1,\cdots,m_\eps\}}\sigma_{i}^\eps$ where $\sigma_i^\eps$ is a $i$-th Steklov eigenvalue on $(\Sigma,g,\beta_\eps)$. Denoting $g_j = g$ if $j=0$ and $g_j = g_{\mathbb{D}}$ if $j\geq 1$, The proof of the following claim exactly follows the proof of Claim \ref{clbadpoints}:
\begin{cl} \label{clbadpointssteklov} Up to a subsequence, there is $0< r_{\star} <1$ and a set of at most $K+1$ bad points $P_j  \subset \partial\Sigma$ and  such that for any $p \in \partial\Sigma \setminus P_0$ and any $r < \min\left( r_\star, d_g(p,P_j ) \right) $, then for $\eps$ small enough,
$$ \sigma_\star\left( \mathbb{D}^+_r(p), g_j, \beta_\eps^j \right) \geq \sigma_K^\eps. $$
\end{cl}
In the following, for $\rho >0$, we denote 
$$ \Omega_\rho^0 = \Sigma \setminus \bigcup_{p\in P_0} \mathbb{D}^+_{\rho}(p) \text{ and } \Omega_\rho^j = \mathbb{D} \setminus \bigcup_{p\in P_j} \mathbb{D}^+_{\rho}(p). $$

\subsubsection{Smallness of $\omega_\eps - 1$ and $\theta_\eps$ near good points of thick parts}

We have the following convergence of $\omega_\eps$ to $1$ and $\theta_\eps$ to $0$ in thick parts, and if $\sigma_i^\eps \to 0$, then $\int_{\partial_s \Omega_\rho^j} \left(\sqrt{\sigma_i^\eps}  \widetilde{\phi_i^\eps}^j  \right)^2 \to 0$ (see Claim \ref{clsmallthetaomega-1} for the proof in the closed case)
\begin{cl}  We have for any $0 < \rho \leq \rho_0$ that for $1\leq j\leq l$ and for $j = 0$ if $\mu_0\neq 0$ that
\begin{equation}  \label{eqomegaepsto1general2steklov}   \int_{\partial_s \Omega_\rho^j}(\omega_\eps^j -1)^2 +\int_{\partial_s \Omega_\rho^j}\left(\theta_\eps^j\right)^2 \leq O\left( \eps\right) \end{equation}
\begin{equation} \label{eq:poincareeigenfunctionssteklov}
\int_{\partial_s \Omega_\rho^j} \left(\sqrt{\sigma_i^\eps}  \widetilde{\phi_i^\eps}^j  \right)^2  \leq O\left( \sigma_i^\eps t_i^\eps \right) 
\end{equation}
as $\eps\to 0$ and then $\rho \to 0$ (where the integrals are computed with respect to $dL_g$ if $j=0$ and the measure $dL_{\mathbb{S}^1}$ if $j\geq 1$).
\end{cl}

\subsubsection{Good annuli close to bad points}
We denote for a point $p$ and $r_2 < r_1$
$$ \mathbb{A}_{r_1,r_2}(p)^+ := \mathbb{D}^+_{r_1}(p) \setminus \mathbb{D}^+_{r_2}(p) $$
Following Claim \ref{clannulusbadpoint} in the closed case,
\begin{cl} \label{clannulusbadpointsteklov}
Let $j \in \{0,\cdots,l\}$ and let $p \in P_j$, then, up to the extraction of a subsequence there is $r>0$ and $s_\eps \to 0$ as $\eps \to 0$ such that
$$ \sigma_\star\left( \mathbb{A}^+_{r,s_\eps}(p), g_j, \beta_\eps^j \right) \geq \sigma_K^\eps $$
\end{cl}

\subsubsection{Non concentration of energies near good points and arbitrarily close to bad points}
The following proof of non-concentration is fairly left to the reader following the proof of Claim \ref{propnonconcentgeneral}.
\begin{cl} \label{propnonconcentgeneralsteklov}
Let $p\in \Sigma \setminus P_0 $ or $\mathbb{D} \setminus P_j$, be a good point
then for any $r$ such that $\sqrt{r} < r_\star(p):= \min\left(r_\star,\frac{d(p,P_j)}{2}\right)$ and any function $\zeta \in \mathcal{C}_c^\infty(\mathbb{D}^+_r(p))$ such that $0 \leq \zeta \leq 1$
\begin{equation} \label{noconcentsteklov}
\lim_{r\to 0} \limsup_{\eps\to 0} L_j^\eps( \zeta ) = \lim_{r\to 0} \limsup_{\eps\to 0}\int_{\mathbb{D}^+_r(p)} \left\vert \nabla \widetilde{\Phi_{\eps}}^j \right\vert_{\sigma_\eps}^2 =  0
\end{equation}
In addition, we have that for a bad point $p \in P_j$ and $r \leq r_\star$, and any function $\zeta \in \mathcal{C}_c^\infty(\mathbb{A}^+_{r,\sqrt{s_\eps}}(p))$ such that $0 \leq \zeta \leq 1$
\begin{equation} \label{noconcentannulussteklov}
\lim_{r\to 0} \limsup_{\eps\to 0} L_j^\eps( \zeta ) = \lim_{r\to 0} \limsup_{\eps\to 0}\int_{\mathbb{A}^+_{r,\sqrt{s_\eps}}(p)} \left\vert \nabla \widetilde{\Phi_{\eps}}^j \right\vert_{\sigma_\eps}^2 =  0
\end{equation}
\end{cl}

\subsection{Construction of local harmonic replacements}
We set 
\begin{equation}\label{eq:defvarphithetahatsigma}\hat{\theta}_\eps^j := (\theta_\eps^j , \sqrt{\sigma_1^\eps} \widetilde{\phi_1^\eps}^{j} ,\cdots,\sqrt{\sigma_I^\eps}\widetilde{\phi_I^\eps}^{j}) \text{ and } \varphi_\eps^j := \left(\widetilde{\phi_1^\eps}^j,\cdots,\widetilde{\phi_{m_\eps}^\eps}^{j}\right) \text{ and } \hat{\sigma}_\eps := ( \sigma_{I+1}^\eps,\cdots,\sigma_{m_\eps}^\eps ). \end{equation} 
First we build a local replacement of $\widetilde{\Phi_\eps}^j$ which will be written $ \sqrt{\omega_\eps^2 - \vert \tau_\eps\vert^2} \Psi_\eps $ where $\vert \tau_\eps \vert$ is a local free boundary harmonic replacement into $\R^{I+1}$ of $\hat{\theta}_\eps^j$ and $\Psi_\eps$ is a local free boundary harmonic replacement into an Euclidean ellipsoid of parameter $\hat{\sigma}_\eps$ of $\frac{\varphi_\eps^j}{\left\vert \varphi_\eps^j \right\vert_{\hat{\sigma}_\eps}}$. In particular, in the following claim, we give a sense to the replacement $\Psi_\eps$ and prove that it can have an arbitrary small energy. We choose $\eps_0 := \eps_{\alpha}'$  in order to have \ref{eq:energyconvexitysteklov} with $\alpha$ an upper bound for $\max\left\{ \sigma_K^\eps, \left(\sigma_{I+1}^\eps\right)^{-1} \right\}$. This implies the uniqueness of the free boundary harmonic replacement.

\begin{cl} \label{clharmreplacementsteklov} There is $\eta>0$ such that for any $p\in \mathbb{S}^1$ ($\partial\Sigma$ if $j=0$ and $\mu_0\neq 0$) there is $r(p)$ and $r(p)^2 \leq  r_\eps(p) \leq r(p)$ such that there are unique maps $\tau_\eps$ and $\Psi_\eps$ satisfying
$$ \tau_\eps = \hat{\theta}_\eps^j \text{ and } \left\vert \varphi_\eps^j \right\vert_{\hat{\sigma}_\eps} \geq \frac{1}{2} \text{ and } \Psi_\eps = \frac{\varphi_\eps^j}{\left\vert \varphi_\eps^j \right\vert_{\hat{\sigma}_\eps}} $$
almost everywhere on $ \partial_d \mathbb{D}_{r_\eps}^+(p)$ and $\vert \psi_\eps \vert_{\hat{\sigma}_\eps} = 1$ on $\partial_s \mathbb{D}_{r_\eps}^+(p)$
\begin{equation*} \int_{\mathbb{D}^+_{r_\eps}(p)} \left\vert \nabla \Psi_\eps \right\vert^2 = \inf \left\{ \int_{\mathbb{D}^+_{r_\eps}(p)} \left\vert \nabla \Psi \right\vert^2 ; \Psi \in H^1  ; \begin{cases} \vert \Psi \vert_{\hat{\sigma}_\eps} =_{a.e} 1 \text{ on } \partial_{s}\mathbb{D}^+_{r_\eps}(p) \\ \Psi =_{a.e} \frac{\varphi_\eps^j}{\left\vert \varphi_\eps^j \right\vert_{\hat{\sigma}_\eps}} \text{ on }  \partial_{d}\mathbb{D}^+_{r_\eps}(p) \end{cases} \right\} \leq \eps_0. \end{equation*}
In particular $\Psi_\eps$ is a free boundary harmonic map into the ellipsoid $\{ \left\vert x \right\vert_{\hat{\sigma}_\eps} = 1 \}$ and
$$ \Delta \Psi_\eps = 0 \text{ in } \mathbb{D}_{r_\eps}^+(p)  \text{ and } \partial_\nu \Psi_\eps = \left(\Psi_\eps\cdot \partial_\nu \Psi_\eps\right) \hat{\sigma}_\eps \Psi_\eps \text{ on } \partial_s \mathbb{D}_{r_\eps}^+(p) $$
$$ \Delta \tau_\eps = 0 \text{ in } \mathbb{D}_{r_\eps}^+(p)  \text{ and } \partial_\nu \tau_\eps = 0 \text{ on } \partial_s \mathbb{D}_{r_\eps}^+(p)  $$
and $\vert \tau_\eps \vert^2 \leq \frac{1}{4}$ in $\mathbb{D}_{r_\eps}^+(p)$.
\end{cl}

\begin{proof}
During all the proof, we drop the indices or exponents $j$ of all the functions because the argument is similar in all the thick parts. Thanks to \eqref{noconcentsteklov}, let $p \in \Sigma \setminus P_0$ or $\mathbb{S}^2 \setminus P_j$, let $r_0(p) \leq r_{\star}$ be such that any small $\eps$,
$$ \int_{\mathbb{D}^+_{r_0(p)}(p)} \left\vert \nabla \varphi_\eps \right\vert^2 \leq \frac{1}{2} \eps_\alpha. $$
for a constant $0<\delta\leq 1$ we will choose later. If $p \in P_j$, with the use of \eqref{noconcentannulussteklov}, we choose $r_0(p)$ such that, 
$$ \int_{\mathbb{A}^+_{\delta_0^{-1} r_0(p),\delta_0 \frac{r_0(p)^2}{4}}(p)} \left\vert \nabla \varphi_\eps \right\vert^2 \leq \frac{1}{2} \eps_\alpha'. $$

Then, by Claim \ref{cluniform} and Claim \ref{cluniformannulus}, knowing that $\vert \varphi_\eps \vert^2_{\hat{\sigma}_\eps}  + \vert \hat{\theta}_\eps \vert^2 \geq 1$ in $\partial_s \mathbb{D}^+_{2r_0(p)}(p)$ if $p\notin P_j$ or $\partial_s \mathbb{A}^+_{\delta_0^{-1} r_0(p),\delta_0 \frac{r_0(p)^2}{4}}(p)$ if $p\in P_j$, that $(\sqrt{\hat{\sigma}_\eps}\varphi_\eps,\hat{\theta}_\eps)$ is a Euclidean harmonic map and that $\int \vert \nabla \hat{\theta}_\eps \vert^2 \to 0 $ as $\eps \to 0$, we obtain for $\alpha = \frac{1}{4}$ that
\begin{equation} \label{eqphiepsthetaepssupsteklov} \vert \varphi_\eps \vert^2_{\hat{\sigma}_\eps}  + \vert \hat{\theta}_\eps \vert^2 \geq \frac{3}{4}  \end{equation}
in $ \mathbb{D}^+_{r_0(p)}(p) $ if $p\notin P_j$ or $\mathbb{A}^+_{r_0(p),\frac{ r_0(p)^2}{4}}(p)$ if $p \in P_j$.

Up to reduce $r_0(p)$, thanks to \eqref{noconcentsteklov} again, we assume
$$ \int_{\mathbb{D}_{r_0}(p)} \left\vert \nabla \varphi_\eps \right\vert^2 \leq \frac{\delta}{2} \eps_0. $$
for a constant $0<\delta\leq 1$ we will choose later. If $p \in P_j$, with the use of \eqref{noconcentannulussteklov} again, we choose $r_0(p)$ such that , 
$$ \int_{\mathbb{A}_{r_0(p),\frac{r_0(p)^2}{4}}(p)} \left\vert \nabla \varphi_\eps \right\vert^2 \leq \frac{\delta}{2} \eps_0. $$
Now we use a symmetrization $(x,y)\mapsto (x,-y)$ to extend $\varphi_\eps$ on $\mathbb{D}_{r_0(p)}$ if $p \notin P_j$ or $\mathbb{A}_{r_0(p),\frac{r_0(p)^2}{2}}(p) $ if $p\in P_j$. By the Courant-Lebesgue lemma with $\frac{r_0(p)}{2} <r< r_0(p)$, let $r^2 \leq r_{\eps} \leq r$ be a radius such that  
\begin{equation} \label{eqcourantlebesguesteklov}
\begin{split} 
\int_{\partial\mathbb{D}_{r_\eps(p)}(p)} & \left\vert \partial_\theta \hat{\theta}_\eps \right\vert ^2 d\theta +
 \int_{\partial\mathbb{D}_{r_\eps(p)}(p)} \left\vert \partial_\theta \varphi_\eps \right\vert ^2 d\theta \\ \leq & 
\frac{1}{\ln 2} 
\left( \int_{\mathbb{A}_{r,r^2}(p)}  \left\vert \nabla \hat{\theta}_\eps \right\vert^2 + 
\int_{\mathbb{A}_{r,r^2}(p)}  \left\vert \nabla \varphi_\eps \right\vert^2 
\right) 
 \leq  \frac{2}{\ln 2} \delta \eps_0 .
\end{split}  \end{equation}
As a consequence, a vector-valued Morrey embedding theorem yields
\begin{equation} \label{eqconsequencecourantlebesguesteklov} 
\max_{q,q' \in \partial\mathbb{D}_{r_\eps(p)}(p)} \left\vert \hat{\theta}_\eps(q) - \hat{\theta}_\eps(q') \right\vert^2 + 
\max_{q,q' \in \partial\mathbb{D}_{r_\eps(p)}(p)} \sum_{i=1}^{n_\eps}  \left\vert \varphi_i^\eps(q) - \varphi_i^\eps(q') \right\vert^2 \leq \frac{2\pi}{\ln 2}  \delta \eps_0. \end{equation}
By the equivalence of the norms
$$ \left( \int_{\partial_s \mathbb{A}^+_{r,r^2}(p)} \varphi^2 + \int_{\mathbb{A}^+_{r,r^2}(p)} \vert \nabla \varphi \vert^2 \right)^{\frac{1}{2}} \text{ and } \left( \int_{\partial \mathbb{D}_{r_\eps(p)}(p)} \varphi^2 + \int_{\mathbb{A}^+_{r,r^2}(p)} \vert \nabla \varphi \vert^2 \right)^{\frac{1}{2}} $$
and by \eqref{eqomegaepsto1general2steklov} we have that
$$  \int_{\partial \mathbb{D}_{r_\eps}(p)} \vert \hat{\theta}_\eps\vert ^2  \leq o(1) $$
as $\eps\to 0$. 
Using \eqref{eqconsequencecourantlebesguesteklov},
$$ \sup_{q\in \partial\mathbb{D}_{r_\eps}(p)} \vert \hat{\theta}_\eps(q)\vert \leq o(1) + \sqrt{\frac{2\pi}{\ln 2}\delta\eps_0} $$
so that choosing $\delta \leq \frac{1}{64} \sqrt{\frac{\ln 2}{\pi\eps_0}} $, by \eqref{eqphiepsthetaepssupsteklov} and by symmetry, $\left\vert \hat{\theta}_\eps \right\vert^2 \leq \frac{1}{4}$ and $\left\vert \varphi_\eps \right\vert_{\hat{\sigma}_\eps}^2 \geq \frac{1}{2}$ on $\partial\mathbb{D}_{r_\eps(p)}(p) $ for $\eps$ small enough. By the maximum principle, $\left\vert \tau_\eps \right\vert^2 \leq \frac{1}{4}$ in $\partial\mathbb{D}_{r_\eps(p)}^+(p) $.

We let $\Psi_\eps : (\mathbb{D}_{r_\eps(p)}(p),\partial \mathbb{D}_{r_\eps(p)}(p)) \to \left(co(\mathcal{E}_{\hat{\sigma}_\eps}),\mathcal{E}_{\hat{\sigma}_\eps}\right)$ be a harmonic extension of $\frac{\varphi_\eps}{\left\vert \varphi_\eps \right\vert_{\hat{\sigma}_\eps}}$ (that is a minimizer of the energy on maps $\Psi$ satisfying $\left\vert \Psi \right\vert_{\hat{\Lambda}_\eps} = 1$ on $\partial \mathbb{D}_{r_\eps(p)}(p)$)
In order to prove uniqueness of $\Psi_\eps$, we have to prove that its energy is small enough. 

Let $\eta \in \mathcal{C}^\infty_c\left( \mathbb{D}_{r^2}(p) \right)$ be a cut-off function such that $\eta \geq 1$ in $\mathbb{D}_{\frac{r^2}{2}}(p)$ and $\vert \nabla \eta \vert \leq \frac{1}{r}$. We set $T_\eps(x) := \left( 1 - \eta \right) \varphi_\eps\left( r_\eps \frac{x}{\vert x \vert} \right) + \eta \varphi_\eps(q_\eps)$ and we compute the energy of $\frac{T_\eps}{\vert T_\eps \vert}$ knowing that
$$ \int_{\mathbb{D}_{r_\eps(p)}(p)} \left\vert \nabla \Psi_\eps \right\vert^2_g dA_g \leq \int_{\mathbb{D}_{r_\eps(p)}(p)} \left\vert \nabla \frac{T_\eps}{\vert T_\eps \vert} \right\vert^2_g dA_g $$
We have that 
$$\left\vert \nabla \frac{T_\eps}{\vert T_\eps \vert} \right\vert^2 \leq  \frac{\vert \nabla T_\eps \vert^2}{\vert T_\eps \vert^2} \leq \frac{2\left( 1-\eta \right)^2 \frac{\vert \nabla_\tau \varphi_\eps \vert^2}{r^2} + 2\vert \nabla \eta \vert^2 \max_{q \in \partial\mathbb{D}_{r_\eps}(p)}\vert \varphi_\eps(q)-\varphi_\eps(q_\eps)  \vert^2}{\left( \vert \varphi_\eps(q_\eps) \vert - \max_{q \in \partial\mathbb{D}_{r_\eps}(p)}\vert \varphi_\eps(q)-\varphi_\eps(q_\eps)  \vert \right)^2} $$
so that using the previous smallness estimates coming from \eqref{eqconsequencecourantlebesguesteklov} 
and up to reduce $\delta$, we complete the proof of the Claim.
\end{proof}

\subsection{Local $H^1$ comparison of eigenfunctions to the harmonic replacements}

\begin{cl} \label{clcomparegeneralsteklov} We have for all $p\in \partial\Sigma$ and $r_\eps(p)$ given by Claim \ref{clharmreplacementsteklov}
$$  \int_{\mathbb{D}_{r_\eps(p)}^+(p)}   \left\vert \nabla \left( \Psi_\eps- \hat{\varphi}_\eps^j \right) \right\vert^2 
= o\left(1\right) $$
as $\eps \to 0$ where with the notations of Claim \ref{clharmreplacementsteklov}
\begin{equation*}
 \hat{\varphi}_\eps^j = \begin{cases} \frac{\varphi_\eps^j}{\rho_\eps^j} \text{ if } p \in \partial\Sigma \setminus P_0 \text{ if } j=0 \text{ or } \mathbb{S}^1 \setminus P_j \text{ if } j \geq 1 \\ \left(1-\eta_\eps\right)  \frac{\varphi_\eps^j}{\rho_\eps^j} +  \eta_\eps \Psi_\eps \text{ if } p \in P_j
 \end{cases} \end{equation*}
and $\rho_\eps^j := \sqrt{\left(\omega_\eps^j\right)^2 - \vert \tau_\eps \vert^2}$ and $\eta_\eps \in \mathcal{C}_c^\infty\left(\mathbb{D}^+_{\sqrt{s_\eps}}(p)\right)$ such that $\eta_\eps = 1$ in $\mathbb{D}^+_{s_\eps}(p)$, $0\leq \eta_\eps \leq 1$ and
\begin{equation}\label{eq:smallenergiesetarhosteklov} \int_{\mathbb{D}_{r_\eps(p)}(p)} \vert \nabla \eta_\eps \vert^2 = O\left(\frac{1}{\ln \frac{1}{s_\eps}}\right) \text{ and } \int_{\mathbb{D}_{r_\eps(p)}(p)} \vert \nabla \rho_\eps^j \vert^2 = O(\eps)\end{equation}
\end{cl}

\begin{proof}
We only write the proof of the claim for $p \in P_j$ since the other case exactly follows the same proof with $\eta_\eps = 0$ and $\mathbb{D}_{r_\eps(p)}(p)$ instead of $\mathbb{A}_{r_\eps(p),s_\eps}(p)$. We drop the index/exponent $j$ in all the proof since it works the same way in every thick part. We let $r_\eps(p)$, $\Psi_\eps$, $\tau_\eps$ be given by Claim \ref{clharmreplacementsteklov}. 
Notice that \eqref{eq:smallenergiesetarhosteklov} on $\rho_\eps^j$ is a simple consequence of Claim \eqref{eqomegaepsto1steklov} and the $(PS)_K$ that gives $\int_\Sigma \vert \nabla \tau_\eps \vert^2 = O(\eps)$.
Notice also that $\rho_\eps$ is chosen so that  $\hat{\varphi}_{\eps}^i - \Psi_{\eps}^i$ is equal to $0$ on $\partial\mathbb{D}_{r_\eps(p)}(p)$. With the choice of $\eta_\eps$, it is equal to $0$ on $\partial\mathbb{A}_{r_\eps(p),s_\eps}(p)$. We will use this property in Step 1 and Step 2. Using both steps will complete the proof of the Claim.

\medskip
\noindent\textbf{Step 1:}
\medskip

\begin{equation} \label{eqsmalldiffenergysteklov}
 \int_{\mathbb{D}^+_{r_\eps(p)}(p)}  \left\vert \nabla \hat{\varphi}_\eps \right\vert^2 - \int_{\mathbb{D}^+_{r_\eps(p)}(p)} \left\vert \nabla  \Psi_\eps \right\vert^2  \leq  o(1)
\end{equation}
as $\eps \to 0$

\medskip
\noindent \textbf{Proof of Step 1:}
We test the function $ \hat{\varphi}_{\eps}^i - \Psi_{\eps}^i$ in the variational characterization of $\sigma_{\star}:= \sigma_{\star}\left(\mathbb{A}^+_{r_\eps(p),s_\eps}(p),\beta_\eps \right)$ knowing Claim \ref{clbadpointssteklov}:
$$ \sigma_{i}^\eps L_\eps \left( \left(\hat{\varphi}_{\eps}^i - \Psi_{\eps}^i \right)^2 \right) \leq  \sigma_{\star} L_\eps \left( \left(\hat{\varphi}_{\eps}^i - \Psi_{\eps}^i \right)^2 \right) \leq \int_{ \mathbb{D}^+_{r_\eps(p)}(p) } \left\vert \nabla\left(\hat{\varphi}_{\eps}^i - \Psi_{\eps}^i\right)\right\vert^2$$
and we sum on $i$ to get
\begin{equation} \label{eqtestsigmastar}
 L_\eps\left( \left\vert \hat{\varphi}_{\eps} - \Psi_{\eps}\right\vert_{\hat{\sigma}_\eps}^2\right) \leq   \int_{ \mathbb{D}^+_{r_\eps(p)}(p) } \left\vert \nabla \hat{\varphi}_{\eps}  \right\vert^2 +  \int_{ \mathbb{D}^+_{r_\eps(p)}(p) } \left\vert \nabla  \Psi_{\eps}  \right\vert^2 
- 2\int_{ \mathbb{D}^+_{r_\eps(p)}(p) } \nabla \hat{\varphi}_{\eps} \nabla \Psi_{\eps} 
 \end{equation}
 Now, we test the equation on $\Phi_\eps$: $\Delta_g \Phi_\eps = \beta_\eps(\sigma_\eps\Phi_\eps,.) $ against $\frac{1-\eta_\eps}{\rho_\eps}\left(\hat{\varphi}_\eps - \Psi_\eps\right)$ and we multiply by $2$:
\begin{equation*}
\begin{split}2 \int_{\mathbb{D}^+_{r_\eps(p)}(p)} \nabla \varphi_\eps \nabla\left( \frac{1-\eta_\eps}{\rho_\eps}\left(\hat{\varphi}_\eps - \Psi_\eps\right)  \right) & = 2 L_{\eps}\left(  \left\langle \varphi_\eps , \frac{1-\eta_\eps}{\rho_\eps}\left(\hat{\varphi}_\eps - \Psi_\eps\right) \right\rangle_{\hat{\sigma}_\eps} \right) \\
& = L_\eps\left( \left\vert \hat{\varphi}_\eps - \Psi_\eps  \right\vert_{\hat{\sigma}_\eps}^2\right) + L_\eps\left((1-\eta_\eps)^2 \frac{ \vert \tau_\eps \vert^2 - \vert \hat{\theta}_\eps \vert^2 }{\omega_\eps^2 - \vert \tau_\eps \vert^2 }  \right)
\end{split}
\end{equation*}
where for the last equality, we used that $\langle X,(X-Y)\rangle_\sigma = \frac{1}{2} \left\vert X-Y \right\vert^2_\sigma + \frac{1}{2} \left(\left\vert X \right\vert^2_\sigma - \left\vert Y \right\vert^2_\sigma \right)$ with $X = (1-\eta_\eps)  \frac{\varphi_\eps}{\rho_\eps}$, $Y = (1-\eta_\eps) \Psi_\eps$ and the equality
$$ \hat{\varphi}_\eps - \Psi_\eps = (1-\eta_\eps) \left( \frac{\varphi_\eps}{\rho_\eps}-\Psi_\eps \right). $$

We obtain that 
\begin{equation*}
\begin{split}  \int_{\mathbb{D}^+_{r_\eps(p)}(p)}  \left\vert \nabla \hat{\varphi}_\eps \right\vert^2 - \int_{\mathbb{D}^+_{r_\eps(p)}(p)} \left\vert \nabla  \Psi_\eps \right\vert^2  \leq  L_\eps\left((1-\eta_\eps)^2 \frac{ \vert \tau_\eps \vert^2 - \vert \hat{\theta}_\eps \vert^2 }{\omega_\eps^2 - \vert \tau_\eps \vert^2 }  \right). \\
+ 2 \int_{\mathbb{D}^+_{r_\eps(p)}(p)} \left(  \nabla \hat{\varphi}_\eps  \nabla\left( \hat{\varphi}_\eps - \Psi_\eps  \right) -\nabla \varphi_\eps \nabla \left( \frac{1-\eta_\eps}{\rho_\eps}\left( \hat{\varphi}_\eps - \Psi_\eps  \right) \right) \right) = I + II
\end{split}
\end{equation*}
The first right-hand term satisfies by a Cauchy-Schwarz inequality and properties of $\sigma_{\star}:= \sigma_{\star}\left(\mathbb{A}_{r_\eps(p),s_\eps}(p),\beta_\eps \right)$
\begin{equation*}
\begin{split} I^2 \leq & 4 L_\eps\left( \left\vert (1-\eta_\eps)(\tau_\eps - \hat{\theta}_\eps) \right\vert^2 \right) L_\eps\left( \left\vert (1-\eta_\eps)(\tau_\eps + \hat{\theta}_\eps) \right\vert^2 \right) \\
\leq & C \frac{1}{\sigma_\star} \int_{\mathbb{D}^+_{r_\eps(p)}(p)} \left\vert \nabla \left((1-\eta_\eps)(\tau_\eps - \hat{\theta}_\eps)   \right) \right\vert^2 \leq o(1) 
 \end{split} \end{equation*}
as $\eps \to 0$ since the energies of $\hat{\theta}_\eps$, $\tau_\eps$ and $\eta_\eps$ go to $0$ as $\eps\to 0$.

The second right-hand term satisfies
\begin{equation*}
\begin{split}  II  = &  2 \int_{\mathbb{D}^+_{r_\eps(p)}(p)} \nabla\left( \hat{\varphi}_\eps - \psi_\eps\right)\nabla \left( \hat{\varphi}_\eps - \varphi_\eps\frac{(1-\eta_\eps)}{\rho_\eps} \right) \\
 & + 2 \int_{\mathbb{D}^+_{r_\eps(p)}(p)} \nabla\frac{\eta_\eps}{\rho_\eps} \left( \left(\hat{\varphi}_\eps-\psi_\eps\right) \nabla \varphi_\eps - \varphi_\eps \nabla \left(\hat{\varphi}_\eps - \Psi_\eps \right) \right)  \\
 = & 2\int_{\mathbb{D}^+_{r_\eps(p)}(p)} \nabla\left( \hat{\varphi}_\eps - \psi_\eps\right)\nabla \left( \eta_\eps \Psi_\eps \right) \\
 & +  2\int_{\mathbb{D}^+_{r_\eps(p)}(p)} \left(\nabla \eta_\eps - \eta_\eps \frac{\nabla \rho_\eps}{\rho_\eps}\right) \left( \left(\hat{\varphi}_\eps-\psi_\eps\right) \frac{\nabla \varphi_\eps}{\rho_\eps} - \frac{\varphi_\eps}{\rho_\eps} \nabla \left(\hat{\varphi}_\eps - \Psi_\eps \right) \right)  \\
 \leq & C \left(\int_{\mathbb{D}^+_{r_\eps(p)}(p)} \eta_\eps^2 \vert \nabla \Psi_\eps \vert^2 + \vert \nabla \eta_\eps \vert^2  \right)^{\frac{1}{2}} + C \left( \int_{\mathbb{D}^+_{r_\eps(p)}(p)}  \vert \nabla \rho_\eps \vert^2 + \vert \nabla \eta_\eps \vert^2 \right)^{\frac{1}{2}} = o(1)
\end{split}\end{equation*}
as $\eps\to 0$ where we used for the inequality that the energy of $\varphi_\eps$ and $\hat{\varphi}_\eps-\Psi_\eps$ is uniformly bounded, that $\rho_\eps^{-1}$, $\frac{\varphi_\eps}{\rho_\eps}$ and $\hat{\varphi}_\eps-\Psi_\eps$ are uniformly bounded in $L^{\infty}$ as $\eps \to 0$. For the last equality, we use that the energy of $\rho_\eps$ and $\eta_\eps$ converges to $0$, and that the $L^\infty$ norm of $\vert \nabla \Psi_\eps \vert^2 $ is uniformly bounded in $\mathbb{D}^+_{\frac{r_\eps(p)}{2}}(p)$ by $\eps$-regularity on free boundary harmonic maps (see Claim \ref{cl:epsregsteklovcase}). Finally we obtain \eqref{eqsmalldiffenergysteklov}

\medskip
\noindent \textbf{Step 2:}
\medskip

$$ \int_{\mathbb{D}^+_{r_\eps(p)}(p)}   \left\vert \nabla \left( \Psi_\eps- \hat{\varphi}_\eps^j \right) \right\vert^2 \leq  \int_{\mathbb{D}^+_{r_\eps(p)}(p)}  \left\vert \nabla \hat{\varphi}_\eps \right\vert^2 - \int_{\mathbb{D}^+_{r_\eps(p)}(p)} \left\vert \nabla  \Psi_\eps \right\vert^2 + o(1) $$
as $\eps \to 0$.

\medskip
\noindent \textbf{Proof of Step 2:}
\medskip

We test the equation on $\Psi_\eps$: $\Delta \Psi_\eps = 0$ and $\partial_\nu \Psi_\eps = \left(\Psi_\eps \cdot \partial_\nu \Psi_\eps \right) \sigma_\eps \Psi_\eps $ against $\Psi_\eps - \hat{\varphi}_\eps$ and we multiply by $2$ to obtain
\begin{equation*}
\begin{split} 2 & \int_{\mathbb{D}^+_{r_\eps(p)}(p)}  \nabla \Psi_\eps \nabla \left( \Psi_\eps - \hat{\varphi}_\eps  \right) = 2 \int_{\partial_s \mathbb{D}^+_{r_\eps(p)}(p)}  \left(\Psi_\eps \cdot \partial_\nu \Psi_\eps \right) \langle \Psi_\eps, \Psi_\eps - \hat{\varphi}_\eps \rangle_{\hat{\sigma}_\eps}  \\ 
= &  \int_{\partial_s \mathbb{D}^+_{r_\eps(p)}(p)} \left(\Psi_\eps \cdot \partial_\nu \Psi_\eps \right) \left( \left\vert \Psi_\eps - \hat{\varphi}_\eps \right\vert_{\sigma_\eps}^2 + (1-\eta_\eps) \frac{  \vert \hat{\theta}_\eps \vert^2 - \vert \tau_\eps \vert^2  }{\omega_\eps^2 - \vert \tau_\eps \vert^2 }  \right)   \\
\leq & C \left(\frac{ \sigma_K^\eps}{ \sigma_{I+1}^\eps}\right)^2\eps_0 \left(  \int_{\mathbb{D}^+_{r_\eps(p)}(p)}  \left\vert \nabla \left( \Psi_\eps - \hat{\varphi}_\eps\right) \right\vert^2 + \left(\int_{\mathbb{D}^+_{r_\eps(p)}(p)} \left\vert \nabla \left( \hat{\theta}_\eps-\tau_\eps \right) \right\vert^2 + \vert \nabla \eta_\eps \vert^2\right)^{\frac{1}{2}}  \right)
\end{split} \end{equation*}
where we used again that $\langle X,(X-Y)\rangle_\sigma = \frac{1}{2} \left\vert X-Y \right\vert^2_\sigma + \frac{1}{2} \left(\left\vert X \right\vert^2_\sigma - \left\vert Y \right\vert^2_\sigma \right)$ with $X = \Psi_\eps$ and $Y = \Psi_\eps - \hat{\varphi}_\eps$ for the second equality. The first inequality is a consequence of the rescaling on $\mathbb{D}^+_{r_\eps(p)}(p)$ of the following Hardy inequality \cite{laurainpetrides} Theorem 3.2
$$\forall u \in H_0^1(\mathbb{D}^+),  \int_{[-1,1]\times\{0\}} \frac{u^2}{1-\vert x \vert} \leq \frac{\pi}{2} \int_{\mathbb{D}^+} \vert \nabla u \vert^2 $$
using the $\eps$-regularity of the energy of harmonic maps (see Claim \ref{cl:epsregsteklovcase}), we have
$$ \vert \Psi_\eps \cdot \partial_\nu \Psi_\eps \vert(x) \leq \vert \nabla \Psi_\eps \vert(x) \leq  \frac{C}{\left(r_\eps(p)-\vert x - p \vert\right)} \sqrt{\int_{\mathbb{D}^+_{r_\eps(p)}(p)}  \vert \nabla \Psi_\eps \vert^2}.$$
Then, we have that
\begin{equation*}
\begin{split} \int_{\mathbb{D}^+_{r_\eps(p)}(p)} &  \left\vert \nabla \left( \Psi_\eps- \hat{\varphi}_\eps^j \right) \right\vert^2 = \int_{\mathbb{D}^+_{r_\eps(p)}(p)} \left( \left\vert \nabla \hat{\varphi}_\eps \right\vert^2 - \left\vert \nabla  \Psi_\eps \right\vert^2 + 2   \nabla \Psi_\eps \nabla \left( \Psi_\eps - \hat{\varphi}_\eps  \right) \right) \\ 
& \leq \int_{\mathbb{D}^+_{r_\eps(p)}(p)}  \left\vert \nabla \hat{\varphi}_\eps \right\vert^2 - \int_{\mathbb{D}^+_{r_\eps(p)}(p)} \left\vert \nabla  \Psi_\eps \right\vert^2 + C' \eps_0 \int_{\mathbb{D}^+_{r_\eps(p)}(p)}   \left\vert \nabla \left( \Psi_\eps- \hat{\varphi}_\eps^j \right) \right\vert^2 + o(1)    \end{split} \end{equation*}
as $\eps \to 0$. Choosing $\eps_0 \leq \left(2C'\right)^{-1}$, we obtain Step 2 and the Claim.
\end{proof}

\subsection{Convergence results on the Palais-Smale sequence}

We consider $\widetilde{\Sigma} := \Sigma \sqcup \bigsqcup_{j=1}^l \left(\mathbb{D}\right)_j$ endowed with the metric $\tilde{g}$ equal to $g$ on $\Sigma$ and the flat metric $g_{\mathbb{D}}$ on $(\mathbb{D})_j$ for $1\leq j\leq l$.
Thanks to the previous claims, we can construct a covering of $\widetilde{\Sigma}$ of disks $\{ \mathbb{D}_{r_\eps(p)}(p) \}_{p \in Q}$ where $Q$ is a finite set independent of $\eps$
such that the conclusions of Claim \ref{clcomparegeneralsteklov} hold on any $\mathbb{D}_{r_\eps(p)}(p)$. We use this property to localize and prove the following:
\begin{cl}
There is $V_0 \in L^\infty_+(\partial\Sigma)$ and $V_1,\cdots,V_j \in L^\infty_+(\mathbb{S}^1)$ such that for any $\eta_0 \in \mathcal{C}_c^\infty\left( \Sigma \setminus P_0 \right)$ and $\eta_j \in \mathcal{C}^\infty_c\left( \mathbb{D} \setminus P_j \right)$ 
for $0 \leq j \leq l$,
\begin{equation}\label{eqconvbetaepstoVsteklov} \beta_j^\eps(\eta_j, 1) - \int \eta_j V_j \leq o(1) \left(\Vert \nabla \eta \Vert_{L^2} + \Vert \eta \Vert_{L^\infty} \right).\end{equation}
as $\eps\to 0$. In particular $\mu_0 = V_0 dL_g$ and $\mu_j = V_j dL_{\mathbb{S}^1}$ for $1 \leq j\leq l$
\end{cl}

\begin{proof}
We prove the result for a given $0 \leq j\leq l$ and we drop the use of $j$ in the indices/exponents of functions. We localize the result: let $\eta$ be a cut-off function at the neighborhood of a good point such that a harmonic replacement given by Claim \ref{clharmreplacement} is well-defined on $K = supp(\eta)$ for any large $\eps$, and such that for any large $\eps$,
$$ \Vert \left\vert \nabla \Psi_\eps \right\vert^2 \Vert_{L^{\infty}(K)} \leq A $$
for some constant $A$ by $\eps$-regularity of free boundary harmonic maps (see Claim \ref{cl:epsregsteklovcase}). Then, $ \Psi_\eps \cdot \partial_\nu \Psi_\eps $ converges to some function $V_j \in L^\infty(\partial_s K)$ strongly in $L^p(\partial_s K)$ for $1\leq p < +\infty$. 
We test the function $\frac{\eta \varphi_\eps}{\rho_\eps^2}$ against the equation on $\varphi_\eps$: $\Delta \varphi_\eps = \hat{\sigma}_\eps\beta_\eps(\varphi_\eps,\cdot)$. We obtain
\begin{equation*}\begin{split} \beta_\eps(1,\eta)  = & \beta_\eps\left(\frac{\vert\varphi_\eps\vert_{\hat{\sigma}_\eps}^2}{\rho_\eps^2},\eta\right) = \hat{\sigma}_\eps\beta_\eps\left(\varphi_\eps,\frac{\varphi_\eps \eta}{\rho_\eps^2}\right) = \int_{K} \nabla \varphi_\eps \nabla \frac{\varphi_\eps \eta}{\rho_\eps^2} \\
 = & \int_K  \frac{\varphi_\eps}{\rho_\eps}  \nabla \frac{\varphi_\eps}{\rho_\eps}  \nabla \eta - \int_K  \frac{\vert\varphi_\eps\vert^2}{\rho_\eps} \nabla \frac{1}{\rho_\eps} \nabla \eta  \\ 
& + \int_K \eta \left\vert \nabla \frac{\varphi_\eps}{\rho_\eps} \right\vert^2 +  \int_K \eta \nabla\frac{1}{\rho_\eps}\left( \frac{\varphi_\eps}{\rho_\eps}  \nabla \varphi_\eps - \varphi_\eps \nabla\frac{\varphi_\eps}{\rho_\eps}\right)    \\
 = & \int_K \left( \eta \vert \nabla \Psi_\eps \vert^2 + \Psi_\eps \nabla \Psi_\eps  \nabla \eta\right) + \int_K  \eta\left( \left\vert \nabla \frac{\varphi_\eps}{\rho_\eps} \right\vert^2 - \vert \nabla \Psi_\eps \vert^2\right)  \\
& + \int_K \left( \Psi_\eps \nabla \Psi_\eps - \frac{\varphi_\eps}{\rho_\eps}\nabla \frac{\varphi_\eps}{\rho_\eps}\right)  \nabla \eta  + \int_K \nabla \frac{1}{\rho_\eps}\left( \frac{\vert\varphi_\eps\vert^2}{\rho_\eps}\nabla \eta + \eta\left( \frac{\varphi_\eps}{\rho_\eps}  \nabla \varphi_\eps - \varphi_\eps \nabla\frac{\varphi_\eps}{\rho_\eps}  \right)  \right) \\
= & \int_K \nabla (\eta\Psi_\eps)\nabla \Psi_\eps + o(1) \left(\Vert \nabla \eta \Vert_{L^2} + \Vert \eta \Vert_{L^\infty} \right) \\
= & \int_{\partial_s K} \eta \Psi_\eps\cdot \partial_\nu \Psi_\eps +o(1)\left(\Vert \nabla \eta \Vert_{L^2} + \Vert \eta \Vert_{L^\infty} \right)
\end{split}\end{equation*}
where the penultimate equality comes from Claim \ref{clcomparegeneralsteklov}. We completed the proof.
\end{proof}

We recall that for a Riemannian surface $(\Sigma,g)$,
$$ I_F(\Sigma,g)  = \inf_{\beta \in \bar{X}} F(\bar{\sigma}_1(\Sigma,g,\beta),\cdots,\bar{\sigma}_m(\Sigma,g,\beta)) $$
From the previous claim, we obtain a measure $VdL_{\tilde{g}}$ equal to $V_0dL_g$ on $\partial \Sigma$ and $V_j dL_{\mathbb{S}^1}$ on $(\mathbb{S}^1)_j$ for $1\leq j\leq l$.
By upper semi-continuity of eigenvalues with respect to bubble tree convergence, and then lower semi-continuity of $f(\Sigma,g,\beta) := F(\bar{\sigma}_1(\Sigma,g,\beta),\cdots,\bar{\sigma}_m(\Sigma,g,\beta))$ with respect to bubble tree convergence, we obtain that
$$ I_F(\Sigma,g) = \liminf_{\eps\to 0} E(\Sigma,g,\beta_\eps) \geq E(\widetilde{\Sigma},\widetilde{g}, V dA_{\tilde{g}} )  ) \geq I_{F}(\widetilde{\Sigma}, \widetilde{g})  $$
In addition, we know by glueing methods that 
$I_{F}(\widetilde{\Sigma}, \widetilde{g}) \geq  I_F(\Sigma,g) $ (see \cite{ces}, \cite{}). Therefore, all the inequalities are equalities and $V dL_{\tilde{g}}$ is a minimizer for $E$ on $\left(\widetilde{\Sigma},\widetilde{g}\right)$. 

By Euler-Lagrange equation applied to the minimizer $V dL_{\tilde{g}}$, we obtain the existence of $\Phi : \left(\widetilde{\Sigma},\partial \widetilde{\Sigma}\right) \to \mathbb{R}^n$ such that setting $\lambda_k := \lambda_k(\widetilde{\Sigma},\widetilde{g}, VdL_{\widetilde{g}})$, and $\sigma := (\sigma_1,\cdots,\sigma_n)$
\begin{itemize}
\item $\Delta_{\widetilde{g}} \Phi = 0$ in $\widetilde{\Sigma}$ and $\partial_{\tilde{\nu}} \Phi = \sigma V \Phi$ on $\partial\widetilde{\Sigma}$
\item $\vert \Phi \vert^2_\sigma \geq 1$ and $\int_{\partial\widetilde{\Sigma}} \vert \Phi \vert^2_\sigma V dL_{\widetilde{g}} = 1 $.
\end{itemize}
Applying Claim \eqref{eqomegaepsto1steklov} with $\theta_\eps = 0$, we obtain that $\vert \Phi \vert^2_\sigma = 1$ on $\partial \widetilde{\Sigma}$, so that $\Phi : \left(\widetilde{\Sigma},\partial \widetilde{\Sigma}\right) \to \left(co(\mathcal{E}_{\sigma}),\mathcal{E}_{\sigma}\right)$ is a free boundary harmonic map. In addition, we have that
$$ V = \Phi\cdot \partial_{\tilde{\nu}} \Phi $$
and since a free boundary harmonic map has to be smooth, $V$ is a smooth function. 

Finally, by the Hopf lemma  $\Phi\cdot \partial_{\tilde{\nu}} \Phi(x) >0$ for any $x$. Indeed, setting $\psi(y) = \langle \sigma \Phi(x),\Phi(y) \rangle$, we have that
$$  \psi(y) \leq \vert\Phi(x)\vert_\sigma \vert \Phi(y) \vert_\sigma \leq 1 $$
for any $y\in \widetilde{\Sigma}$ where we used the maximum principle on the subharmonic map $\vert \Phi \vert_{\sigma}^2$ that is equal to $1$ on $\partial \widetilde{\Sigma}$.  This inequality is an equality if $y = x$. Since $\psi$ is harmonic, we have that by the Hopf lemma that $\partial_\nu \psi(x) >0$. Therefore, noting that $\partial_{\tilde{\nu}}\Psi$ is parallel to $\psi$ in $\partial \Sigma$, $V(x) = \Phi\cdot \partial_{\tilde{\nu}} \Phi(x) = \partial_\nu \psi(x) >0 $.

\section{Regularity estimates for harmonic maps independent of the dimension of the target ellipsoid}
\label{sec4}

\begin{cl} \label{cl:epsregclosedcase}
For any $\alpha>1$, there is $C_\alpha > 0$ and $\eps_\alpha>0$ such that for every $n\in \mathbb{N}$ and $\Lambda = (\lambda_1,\cdots,\lambda_n)$ with
$$ \max_{1\leq i \leq n} \lambda_i \leq \alpha \text{ and } \min_{1\leq i \leq n} \lambda_i \geq \alpha^{-1}, $$
such that $\Phi : \mathbb{D} \to \mathcal{E}_{\Lambda}$ is a harmonic map satisfying
$$ \int_{\mathbb{D}} \left\vert \nabla \Phi \right\vert^2 \leq \eps_\alpha $$
Then
$$ \Vert \nabla \Phi \Vert_{L^{\infty}\left(\mathbb{D}_{\frac{1}{2}}\right)}^2 \leq C_{\alpha}  \int_{\mathbb{D}} \left\vert \nabla \Phi \right\vert^2 $$
\end{cl}

\begin{cor}[Energy convexity of harmonic maps \cite{coldingminicozzi}\cite{laurainpetrides}]
For any $\alpha>1$, there is $0<\eps_\alpha'< \eps_\alpha$ such that for every $n\in \mathbb{N}$ and $\Lambda = (\lambda_1,\cdots,\lambda_n)$ with
$$ \max_{1\leq i \leq n} \lambda_i \leq \alpha \text{ and } \min_{1\leq i \leq n} \lambda_i \geq \alpha^{-1}, $$
such that $\Psi : \mathbb{D} \to \mathcal{E}_{\Lambda}$ is a harmonic map satisfying
$$ \int_{\mathbb{D}} \left\vert \nabla \Psi \right\vert^2 \leq \eps_\alpha' $$
Then, for any map $\Phi \in H^1(\mathbb{D},\mathcal{E}_\Lambda)$ such that $\Phi =_{a.e} \Psi$ on $\partial \mathbb{D}$, then
\begin{equation} \label{eq:energyconvexity} \frac{1}{2} \int_\mathbb{D}\vert \nabla \left(\Phi-\Psi\right) \vert^2 \leq \int_\mathbb{D}\vert \nabla \Phi \vert^2 - \int_\mathbb{D}\vert \nabla \Psi\vert^2  \end{equation}
\end{cor}

\begin{cl}{\cite{scheven}, lemma 3.1} \label{cluniform}
for any $\alpha>0$ there is $\eps_{\alpha} >0$  
such that for any $n\in \mathbb{N}$ and $\Phi : \mathbb{D}^+ \to \mathbb{R}^n$ a Euclidean harmonic map such that $\left\vert \Phi \right\vert^2 \geq 1$ on $[-1,1]\times \{0\}$ and such that
$$ \int_{\mathbb{D}_+} \left\vert \nabla \Phi \right\vert^2 \leq \eps_{\alpha}, $$
we have $\left\vert \Phi \right\vert^2 \geq 1 - \alpha $ on $ \mathbb{D}_{\frac{1}{2}}^+$.
\end{cl}

\begin{cl} \label{cluniformannulus}
There is a small $0<\delta_0 < 1$ such that for any $\alpha>0$ there is $\eps_{\alpha}' >0$  
such that for any $n\in \mathbb{N}$ any $r>0$, and $\Phi : \mathbb{A}^+_{\delta_0^{-1},\delta_0 r} \to \mathbb{R}^n$ a Euclidean harmonic map such that $\left\vert \Phi \right\vert^2 \geq 1$ on $\left([-\delta_0^{-1},-\delta_0 r]\cup[\delta_0 r,\delta_0^{-1}]\right)\times \{0\}$ and such that
$$ \int_{\mathbb{A}^+_{\delta_0^{-1},\delta_0 r}} \left\vert \nabla \Phi \right\vert^2 \leq \eps_{\alpha}', $$
we have $\left\vert \Phi \right\vert^2 \geq 1 - \alpha $ on $ \mathbb{A}_{1,r}^+$.
\end{cl}

\begin{cl} \label{cl:epsregsteklovcase}
For any $\alpha>1$, there is $C_\alpha > 0$ and $\eps_\alpha>0$ such that for every $n \in \mathbb{N}$ and $\sigma = (\sigma_1,\cdots,\sigma_n)$ with
$$ \max_{1\leq i \leq n} \sigma_i \leq \alpha \text{ and } \min_{1\leq i \leq n} \sigma_i \geq \alpha^{-1}, $$
such that $\Phi : (\mathbb{D}^+,[-1,1]) \to \left(co\left(\mathcal{E}_{\sigma}\right),\mathcal{E}_{\sigma} \right)$ is a free boundary harmonic map satisfying
$$ \int_{\mathbb{D}^+} \left\vert \nabla \Phi \right\vert^2 \leq \eps_\alpha $$
Then
$$ \Vert \nabla \Phi \Vert_{L^{\infty}\left(\mathbb{D}_{\frac{1}{2}}^+\right)}^2 \leq C_{\alpha}  \int_{\mathbb{D}^+} \left\vert \nabla \Phi \right\vert^2 $$
\end{cl}

\begin{cor}[energy convexity of free boundary harmonic maps \cite{laurainpetrides}]
For any $\alpha>1$, there is $0<\eps_\alpha'< \eps_\alpha$ such that for every $n\in \mathbb{N}$ and $\sigma = (\sigma_1,\cdots,\sigma_n)$ with
$$ \max_{1\leq i \leq n} \sigma_i \leq \alpha \text{ and } \min_{1\leq i \leq n} \sigma_i \geq \alpha^{-1}, $$
such that $\Psi : \left(\mathbb{D}_+, [-1,1]\right) \to \left(co\left(\mathcal{E}_{\sigma}\right),\mathcal{E}_{\sigma}\right)$ is a harmonic map satisfying
$$ \int_{\mathbb{D}_+} \left\vert \nabla \Psi \right\vert^2 \leq \eps_\alpha' $$
Then, for any map $\Phi \in H^1(\mathbb{D}_+,\R^n)$ such that $\Phi =_{a.e} \Psi$ on $\mathbb{D}_+ \cap \partial \mathbb{D}$ and $\vert\Phi\vert_\sigma =_{a.e} 1$ on $[-1,1]$, then
\begin{equation} \label{eq:energyconvexitysteklov} \frac{1}{2} \int_{\mathbb{D}_+}\vert \nabla \left(\Phi-\Psi\right) \vert^2 \leq \int_{\mathbb{D}_+}\vert \nabla \Phi \vert^2 - \int_{\mathbb{D}_+} \vert \nabla \Psi\vert^2  \end{equation}
\end{cor}

\end{document}